\documentclass{tac}
\usepackage{amssymb}
\usepackage{tikz}
\usetikzlibrary{calc}
\usepackage[mathscr]{eucal}
\usepackage{float}
\usepackage{afterpage}

\title[Semantic Proof of Confluence]%
{Semantic Proof of Confluence of the Categorical Reduction System
 for Linear Logic}
\author{Ryu Hasegawa}
\address{Graduate School
 of Mathematical Sciences, The University of Tokyo, Komaba 3-8-1,
 Meguro-ku, Tokyo 153-8914, Japan}


\makeatletter

\newdimen\theight
\def\rightref#1{%
             \vadjust{\setbox0=\hbox{\quad\quad
             \vbox{\leftskip=0pt\hsize=10em\raggedright\small\sf{#1}}}%
             \theight=10pt
             \kern -\theight \vbox to \theight{\rightline{\rlap{\box0}}%
             \vss}%
             }}%

\newif\ifcomment
\commentfalse

\def\mor#1{\mathrel{\mathop{\longrightarrow}\limits^{\vbox
 to0pt{\vss \hbox{$\scriptstyle #1$}\kern-2pt}}}}

\newdimen\latticeUnit
\newdimen\boxMargin \newdimen\layerMargin

\newbox\tbox \newbox\ttbox
\chardef\tc=8 \chardef\tcc=9
\newdimen\tdimen \newdimen\ttdimen
\newdimen\sdimen \newdimen\ssdimen
\newtoks\ttoks \newtoks\tttoks
\newcount\tcount \newcount\ttcount
\newcount\tttcount \newcount\ttttcount

\newbox\diagbox
\newdimen\rightlim \newdimen\leftlim
\newdimen\upperlim \newdimen\lowerlim
\newdimen\dimA \newdimen\dimB
\newdimen\dimC \newdimen\dimD
\newdimen\dimE \newdimen\dimF \newdimen\dimG

\def\clist{\\}

\def\rightappend#1\to#2{\ttoks={#1\\}\tttoks=\expandafter{#2}%
 \edef#2{\the\tttoks\the\ttoks}}

{\catcode`p=12 \catcode`t=12 \gdef\gazonc#1pt{#1}}
\let\getfactor=\gazonc
\def\dimensionToNumber#1{\expandafter\getfactor\the#1}
\let\ptless=\dimensionToNumber

\def\nodebox#1{%
 \futurelet\com\nnodeboxx#1\with}

\def\nnodeboxx{%
 \ifx\com\invisible
  \let\next=\invisibleNodeBox
 \else
  \let\next=\visibleNodeBox
 \fi\next}

\def\visibleNodeBox#1\with#2#3{%
 \def\bm{#3}
 \setbox#2=\hbox{\kern\bm\vbox{\kern\bm\hbox{$\displaystyle{#1}$%
  }\kern\bm}\kern\bm}}

\def\invisible{}

\def\invisibleNodeBox#1\with#2#3{%
 \def\bm{#3}
 \setbox#2=\hbox{\kern\bm\vbox{\kern\bm\hbox{$\displaystyle{#1}$%
  }\kern\bm}\kern\bm}
 \setbox#2=\hbox{\vrule width0pt height\ht#2 depth\dp#2%
  \vbox to0pt{\hrule height0pt depth0pt width\wd#2}}}

\gdef\object(#1,#2)=#3{%
 \nodebox{#3}\tbox\boxMargin
 \tdimen=.5\wd\tbox \ttdimen=.5\ht\tbox \advance\ttdimen by.5\dp\tbox
 \edef\bbbox{(#1,#2)(\ptless{\tdimen}pt,\ptless{\ttdimen}pt)}%
 \expandafter\rightappend\bbbox\to\clist
 \advance\tdimen by-#1\latticeUnit \advance\ttdimen by-#2\latticeUnit
 \tdimen=-\tdimen \ttdimen=-\ttdimen
 \sdimen=\tdimen
  \ifdim\sdimen<\leftlim \global\leftlim=\sdimen \fi
  \advance\sdimen by\wd\tbox
  \ifdim\sdimen>\rightlim \global\rightlim=\sdimen \fi
 \sdimen=\ttdimen
  \ifdim\sdimen<\lowerlim \global\lowerlim=\sdimen \fi
  \advance\sdimen by\ht\tbox
  \ifdim\sdimen>\upperlim \global\upperlim=\sdimen \fi
 \put(\ptless{\tdimen},\ptless{\ttdimen}){\unhbox\tbox}
}

\def\refer(#1,#2)\to(#3,#4){\def\a{#1}\def\b{#2}%
 \tcount=#3 \ttcount=#4
 \expandafter\rreferr\clist\empty}
\def\rreferr\\#1{\ifx#1\empty \let\next=\relax
 \else \let\next=\rrreferrr \fi \next}
\def\rrreferrr#1,#2)(#3,#4){\def\aa{#1}\def\bb{#2}%
 \ifdim\a pt=\aa pt \ifdim\b pt=\bb pt
  \dimen\tcount=#3 \dimen\ttcount=#4
 \fi\fi \rreferr}

\gdef\edge{\let\lcommand=\line
 \morphismSwitch}

\gdef\morphism{\let\lcommand=\vector
 \morphismSwitch}

\def\morphismSwitch(#1,#2)to(#3,#4){%
 \def\mNext{\ifx\farg[\morphismBody(#1,#2)(#3,#4)%
  \else \slantMorphism(#1,#2)(#3,#4)\fi}
 \futurelet\farg\mNext}

\def\morphismBody(#1,#2)(#3,#4){%
 \slantMorphism(#1,#2)(#3,#4)%
 \expandafter\attachSwitch}

\def\abs#1{\ifdim#1<0pt-#1\else#1\fi}

\def\attachSwitch#1#2]{
 \nodebox{{\scriptstyle #2}}\tbox\boxMargin
 \rightsidefalse
 \diagonalfalse \horizontalfalse \verticalfalse
 \def\tNext{\ifx\targ[\expandafter\attachOption
  \else \attachBody\fi}
 \futurelet\targ\tNext}

\newif\ifrightside
\newif\ifhorizontal
\newif\ifvertical
\newif\ifdiagonal

\def\attachOption#1{%
 \attachOptionLoop}
\def\attachOptionLoop#1{%
 \ifx#1]
  \let\next\attachBody
 \else
  \ifx#1R \rightsidetrue \fi
  \ifx#1H \horizontaltrue \fi
  \ifx#1V \verticaltrue \fi
  \ifx#1D \diagonaltrue \fi
  \let\next\attachOptionLoop
 \fi
 \next
}

\newtoks\labelPosition
\labelPosition={.5}

\def\attachBody{%
 \dimC=\xC \advance\dimC by-\xB
 \dimD=\yC \advance\dimD by-\yB
 \dimA=\xB \advance\dimA by\the\labelPosition\dimC
 \dimB=\yB \advance\dimB by\the\labelPosition\dimD
 \ifhorizontal \ifdim0pt<\abs\dimD \attachBodyHorizontal \fi
 \else \ifvertical \ifdim0pt<\abs\dimC \attachBodyVertical \fi
 \else \ifdiagonal \attachBodyDiagonal
 \else
  \tdimen.2\dimD \tdimen=\abs\tdimen
  \ifdim\abs\dimC<\tdimen
   \attachBodyHorizontal
  \else \tdimen.2\dimC \tdimen=\abs\tdimen
  \ifdim\abs\dimD<\tdimen
   \attachBodyVertical
  \else
   \attachBodyDiagonal
  \fi\fi
 \fi\fi\fi}

\def\attachBodyHorizontal{%
 \dimG=.5\ht\tbox \advance\dimG by.5\dp\tbox
 \ifrightside \dimE=\dimD \else \dimE=-\dimD \fi
 \dimF=\dimC
 \ifdim\dimD<0.0pt
  \ifrightside \advance\dimA by-.5\wd\tbox
  \else \advance\dimA by.5\wd\tbox \fi
 \else
  \ifrightside \advance\dimA by.5\wd\tbox
  \else \advance\dimA by-.5\wd\tbox \fi
 \fi
 \displaceLabel(\ptless{\dimA},\ptless{\dimB})(\ptless{\dimE},\ptless{\dimF},\ptless{\dimG},1)
}

\def\attachBodyVertical{%
 \dimG=.5\wd\tbox
 \ifrightside \dimE=\dimC \else \dimE=-\dimC \fi
 \dimF=\dimD
 \ifdim\dimC<0.0pt
  \ifrightside \advance\dimB by.5\ht\tbox \advance\dimB by.5\dp\tbox
  \else \advance\dimB by-.5\ht\tbox \advance\dimB by-.5\dp\tbox \fi
 \else
  \ifrightside \advance\dimB by-.5\ht\tbox \advance\dimB by-.5\dp\tbox
  \else \advance\dimB by.5\ht\tbox \advance\dimB by.5\dp\tbox \fi
 \fi
 \displaceLabel(\ptless{\dimA},\ptless{\dimB})(\ptless{\dimE},\ptless{\dimF},\ptless{\dimG},-1)
}

\def\attachBodyDiagonal{%
 \ifdim\dimD<0.0pt
  \ifrightside
   \advance\dimA by-.5\wd\tbox \advance\dimA by.5\boxMargin
  \else
   \advance\dimA by.5\wd\tbox \advance\dimA by-.5\boxMargin
  \fi
 \else
  \ifrightside
   \advance\dimA by.5\wd\tbox \advance\dimA by-.5\boxMargin
  \else
   \advance\dimA by-.5\wd\tbox \advance\dimA by.5\boxMargin
  \fi
 \fi
 \ifdim\dimC<0.0pt
  \ifrightside
   \advance\dimB by.5\ht\tbox \advance\dimB by.5\dp\tbox
   \advance\dimB by-.5\boxMargin
  \else
   \advance\dimB by-.5\ht\tbox \advance\dimB by-.5\dp\tbox
   \advance\dimB by.5\boxMargin
  \fi
 \else
  \ifrightside
   \advance\dimB by-.5\ht\tbox \advance\dimB by-.5\dp\tbox
   \advance\dimB by.5\boxMargin
  \else
   \advance\dimB by.5\ht\tbox \advance\dimB by.5\dp\tbox
   \advance\dimB by-.5\boxMargin
  \fi
 \fi
 \displaceLabel(\ptless{\dimA},\ptless{\dimB})(1.0,0.0,0.0,0)
}

\def\displaceLabel(#1,#2)(#3,#4,#5,#6){%
 \put(#1,#2){%
  \dimA=#4pt
  \ifdim\dimA<0.0pt \dimA=-\dimA\fi
  \ifnum#6>0
   \begin{tikzpicture}[xscale=0.0352778, yscale=0.0352778, thin, inner sep=0]
     \path[use as bounding box] (-#5*\ptless\dimA/#3,0);
     \node at (0,0) {%
       \hbox to0pt{\hss \vbox to0pt{\vss
       \hbox{\copy\tbox}%
       \vss}\hss}};
   \end{tikzpicture}
  \else
   \begin{tikzpicture}[xscale=0.0352778, yscale=0.0352778, thin, inner sep=0]
     \path[use as bounding box] (0,#5*\ptless\dimA/#3);
     \node at (0,0) {%
       \hbox to0pt{\hss \vbox to0pt{\vss
       \hbox{\copy\tbox}%
       \vss}\hss}};
   \end{tikzpicture}
  \fi}
}

\newdimen\hd \newdimen\vd \newdimen\cd \newdimen\md
\newdimen\xB \newdimen\yB
\newdimen\xC \newdimen\yC

\def\slantMorphism(#1,#2)(#3,#4){%
 \hd=#3\latticeUnit \advance\hd by-#1\latticeUnit
 \vd=#4\latticeUnit \advance\vd by-#2\latticeUnit
 \vectorPosition(#1,#2)(#3,#4)%
 \edef\arg{{\ptless{\xB}}{\ptless{\yB}}{\ptless{\xC}}{\ptless{\yC}}}
  \ifx\lcommand\vector
   \expandafter\drawVect\arg
  \else
   \expandafter\drawLine\arg
  \fi
}

\newif\iflayer
\layerfalse

\def\drawLine#1#2#3#4{%
 \ifdim#1pt<#3pt \dimA=#1pt \else \dimA=#3pt\fi
 \ifdim#2pt<#4pt \dimB=#2pt \else \dimB=#4pt\fi
 \iflayer
  \put(\ptless\dimA,\ptless\dimB){%
    \tikz[xscale=0.0352778, yscale=0.0352778]
      \draw[line width=2..5,white] (#1,#2) -- (#3,#4);
    \tikz[xscale=0.0352778, yscale=0.0352778]
      \draw[line width=.5,black] (#1,#2) -- (#3,#4);}
 \else
  \put(\ptless\dimA,\ptless\dimB){%
    \tikz[xscale=0.0352778, yscale=0.0352778]
      \draw[line width=.5,black] (#1,#2) -- (#3,#4);}
 \fi
}
\def\drawVect#1#2#3#4{%
 \ifdim#1pt<#3pt \dimA=#1pt \else \dimA=#3pt\fi
 \ifdim#2pt<#4pt \dimB=#2pt \else \dimB=#4pt\fi
 \iflayer
  \put(\ptless\dimA,\ptless\dimB){%
    \tikz[xscale=0.0352778, yscale=0.0352778]
      \draw[>=latex,->,line width=2..5,white] (#1,#2) -- (#3,#4);
    \tikz[xscale=0.0352778, yscale=0.0352778]
      \draw[>=latex,->,line width=.5,black] (#1,#2) -- (#3,#4);}
 \else
  \put(\ptless\dimA,\ptless\dimB){%
    \tikz[xscale=0.0352778, yscale=0.0352778]
      \draw[>=latex,->,line width=.5,black] (#1,#2) -- (#3,#4);}
 \fi
}

{\catcode`p=12 \catcode`t=12 \gdef\hazonc#1.#2pt{#1}}
\def\toJnt#1{\expandafter\hazonc\the#1}
\let\getInt=\toJnt

\newdimen\xd \newdimen\yd \newdimen\xe \newdimen\ye

\def\vectorPosition(#1,#2)(#3,#4){%
 \ifdim \hd>0.0pt
  \ifdim \vd>0.0pt
   \obtainDelta(#1,#2)(#3,#4)
  \else
   \vd=-\vd
   \obtainDelta(#1,#2)(#3,#4)
   \vd=-\vd \yd=-\yd \ye=-\ye
  \fi
 \else
  \hd=-\hd
  \ifdim \vd>0.0pt
   \obtainDelta(#1,#2)(#3,#4)
   \xd=-\xd \xe=-\xe
  \else
   \vd=-\vd
   \obtainDelta(#1,#2)(#3,#4)
   \vd=-\vd \xd=-\xd \xe=-\xe \yd=-\yd \ye=-\ye
  \fi
  \hd=-\hd
 \fi
 \xB=#1\latticeUnit \advance\xB by\xd
 \yB=#2\latticeUnit \advance\yB by\yd
 \xC=#3\latticeUnit \advance\xC by\xe
 \yC=#4\latticeUnit \advance\yC by\ye
}

\def\obtainDelta(#1,#2)(#3,#4){%
 \refer(#1,#2)\to(\tc,\tcc)%
 \ifdim\hd<\vd 
  \dimA=\dimen\tcc
  \multiply\dimA by\getInt\hd
  \divide\dimA by\getInt\vd    
  \dimB=\dimen\tc
  \ifdim \dimB<\dimA
   \dimA=\dimen\tc
   \multiply\dimA by\getInt\vd
   \divide\dimA by\getInt\hd   
   \xd=\dimen\tc \yd=\dimA
  \else
   \xd=\dimA \yd=\dimen\tcc
  \fi
 \else 
  \dimA=\dimen\tc
  \multiply\dimA by\getInt\vd
  \divide\dimA by\getInt\hd    
  \dimB=\dimen\tcc
  \ifdim \dimB<\dimA
   \dimA=\dimen\tcc
   \multiply\dimA by\getInt\hd
   \divide\dimA by\getInt\vd   
   \xd=\dimA \yd=\dimen\tcc
  \else
   \xd=\dimen\tc \yd=\dimA
  \fi
 \fi
 \refer(#3,#4)\to(\tc,\tcc)%
 \ifdim\hd<\vd 
  \dimA=\dimen\tcc
  \multiply\dimA by\getInt\hd
  \divide\dimA by\getInt\vd    
  \dimB=\dimen\tc
  \ifdim \dimB<\dimA
   \dimA=\dimen\tc
   \multiply\dimA by\getInt\vd
   \divide\dimA by\getInt\hd   
   \xe=-\dimen\tc \ye=-\dimA
  \else
   \xe=-\dimA \ye=-\dimen\tcc
  \fi
 \else 
  \dimA=\dimen\tc
  \multiply\dimA by\getInt\vd
  \divide\dimA by\getInt\hd    
  \dimB=\dimen\tcc
  \ifdim \dimB<\dimA
   \dimA=\dimen\tcc
   \multiply\dimA by\getInt\hd
   \divide\dimA by\getInt\vd   
   \xe=-\dimA \ye=-\dimen\tcc
  \else
   \xe=-\dimen\tc \ye=-\dimA
  \fi
 \fi
}

\def\ptToCoord#1{%
 \tcount=\expandafter\toInt\number\ptless\latticeUnit
 \divide#1 by\tcount}

\newenvironment{diagramme}{%
  \latticeUnit=100pt \boxMargin=3pt \layerMargin=3pt
  \rightlim=0pt \leftlim=0pt \upperlim=0pt \lowerlim=0pt
  \setbox\diagbox=\hbox\bgroup
  \begin{picture}(0,0)}%
 {\end{picture}\egroup
  \tdimen=\rightlim \advance\tdimen by-\leftlim
  \ttdimen=\upperlim \advance\ttdimen by-\lowerlim
  \begin{picture}(\ptless{\tdimen},\ptless{\ttdimen})%
                 (\ptless{\leftlim},\ptless{\lowerlim})
   \put(0,0){\unhbox\diagbox}
  \end{picture}%
}

\newcount\tabsc

\def\spandiagabs#1#2#3#4#5#6#7{%
 \hbox{\begin{diagramme}
  \spandiagabscoord{#1}{#2}{#3}{#4}{#5}{#6}{#7}
 \end{diagramme}}}

\def\spandiagabscoord#1#2#3#4#5#6#7{%
 \tabsc=#1 \divide\tabsc by2
 \let\ss=\scriptstyle
 \latticeUnit=1pt\boxMargin=3pt
 \object(\the\tabsc,#2)={#3}
 \object(0,0)={#4}
 \object(#1,0)={#5}
 \morphism(\the\tabsc,#2)to(0,0)[\ss{#6}][R]
 \morphism(\the\tabsc,#2)to(#1,0)[\ss{#7}]
}

\def\trianglediagabscoord#1#2#3#4#5#6#7#8{%
 \tabsc=#1 \divide\tabsc by2
 \let\ss=\scriptstyle
 \latticeUnit=1pt\boxMargin=3pt
 \object(0,#2)={#3}
 \object(#1,#2)={#4}
 \object(\the\tabsc,0)={#5}
 \morphism(0,#2)to(#1,#2)[\ss{#6}]
 \morphism(0,#2)to(\the\tabsc,0)[\ss #7][R]
 \morphism(#1,#2)to(\the\tabsc,0)[\ss #8]
}

\def\optrianglediagabscoord#1#2#3#4#5#6#7#8{%
 \tabsc=#1 \divide\tabsc by2
 \let\ss=\scriptstyle
 \latticeUnit=1pt\boxMargin=3pt
 \object(\the\tabsc,#2)={#3}
 \object(0,0)={#4}
 \object(#1,0)={#5}
 \morphism(\the\tabsc,#2)to(0,0)[\ss{#6}][R]
 \morphism(\the\tabsc,#2)to(#1,0)[\ss{#7}]
 \morphism(0,0)to(#1,0)[\ss{#8}][R]
}

\def\squarediagabscoord#1#2{%
 \def\hori{#1}\def\vert{#2}%
 \sqda}
\def\sqda#1#2#3#4#5#6#7#8{%
 \let\ss=\scriptstyle
 \latticeUnit=1pt \boxMargin=3pt
 \object(0,\vert)={#1}
 \object(\hori,\vert)={#2}
 \object(0,0)={#3}
 \object(\hori,0)={#4}
 \morphism(0,\vert)to(\hori,\vert)[\ss{#5}]
 \morphism(0,\vert)to(0,0)[\ss{#6}][R]
 \morphism(\hori,\vert)to(\hori,0)[\ss{#7}]
 \morphism(0,0)to(\hori,0)[\ss{#8}][R]
}

\newbox\tcb\setbox\tcb=\hbox{$\scriptstyle \Rightarrow$}
\newdimen\dimH \newdimen\dimI

\def\ptless#1{\expandafter\getfactor\the#1}

\newdimen\dmA \newdimen\dmB \newdimen\dmC \newdimen\dmD \newdimen\dmE

\def\lzda#1#2#3#4#5#6#7#8{%
 \let\ss=\scriptstyle
 \latticeUnit=1pt\boxMargin=3pt
 \object(0,\ptless\dmB)={#1}
 \object(-\ptless\dmA,0)={#2}
 \object(\ptless\dmA,0)={#3}
 \object(0,-\ptless\dmB)={#4}
 \morphism(0,\ptless\dmB)to(-\ptless\dmA,0)[\ss{#5}][R]
 \morphism(0,\ptless\dmB)to(\ptless\dmA,0)[\ss{#6}]
 \morphism(-\ptless\dmA,0)to(0,-\ptless\dmB)[\ss{#7}][R]
 \morphism(\ptless\dmA,0)to(0,-\ptless\dmB)[\ss{#8}]
}

\def\pentagondiagabscoord#1#2#3#4#5#6#7{%
 \def\firstobj{#3}\def\secondobj{#4}\def\thirdobj{#5}%
 \def\fourthobj{#6}\def\fifthobj{#7}%
 \argargpgda{#1}{#2}}

\def\argargpgda#1#2#3#4#5#6#7{%
  \let\ss=\scriptstyle
  \latticeUnit=1pt \boxMargin=3pt
  \pgda{#1}{#2}{#3}{#4}{#5}{#6}{#7}%
}

\def\pgda#1#2#3#4#5#6#7{%
 \dmA=#1pt\dmB=#2pt
 \dmE=\dmB
 \multiply\dmE by95\divide\dmE by100
 \multiply\dmB by59\divide\dmB by100
 \dmC=\dmA
 \multiply\dmC by81\divide\dmC by100
 \dmD=\dmA
 \multiply\dmD by19\divide\dmD by100
 \divide\dmA by2
 \object(\ptless\dmA,\ptless\dmE)={\firstobj}
 \object(0,\ptless\dmB)={\secondobj}
 \object(#1,\ptless\dmB)={\thirdobj}
 \object(\ptless\dmD,0)={\fourthobj}
 \object(\ptless\dmC,0)={\fifthobj}
 \morphism(\ptless\dmA,\ptless\dmE)to(0,\ptless\dmB)[\ss{#3}][R]
 \morphism(\ptless\dmA,\ptless\dmE)to(#1,\ptless\dmB)[\ss{#4}]
 \morphism(0,\ptless\dmB)to(\ptless\dmD,0)[\ss{#5}][R]
 \morphism(#1,\ptless\dmB)to(\ptless\dmC,0)[\ss{#6}]
 \morphism(\ptless\dmD,0)to(\ptless\dmC,0)[\ss{#7}][R]
}

\def\oppentagondiagabscoord#1#2#3#4#5#6#7{%
 \def\firstobj{#3}\def\secondobj{#4}\def\thirdobj{#5}%
 \def\fourthobj{#6}\def\fifthobj{#7}%
 \opargargpgda{#1}{#2}}

\def\opargargpgda#1#2#3#4#5#6#7{%
  \let\ss=\scriptstyle
  \latticeUnit=1pt \boxMargin=3pt
  \oppgda{#1}{#2}{#3}{#4}{#5}{#6}{#7}%
}

\def\oppgda#1#2#3#4#5#6#7{%
 \dmA=#1pt\dmB=#2pt
 \dmE=\dmB
 \multiply\dmE by95\divide\dmE by100
 \multiply\dmB by59\divide\dmB by100
 \dmC=\dmA
 \multiply\dmC by81\divide\dmC by100
 \dmD=\dmA
 \multiply\dmD by19\divide\dmD by100
 \divide\dmA by2
 \object(\ptless\dmA,-\ptless\dmE)={\fifthobj}
 \object(#1,-\ptless\dmB)={\fourthobj}
 \object(0,-\ptless\dmB)={\thirdobj}
 \object(\ptless\dmC,0)={\secondobj}
 \object(\ptless\dmD,0)={\firstobj}
 \morphism(\ptless\dmD,0)to(\ptless\dmC,0)[\ss{#3}]
 \morphism(\ptless\dmD,0)to(0,-\ptless\dmB)[\ss{#4}][R]
 \morphism(\ptless\dmC,0)to(#1,-\ptless\dmB)[\ss{#5}]
 \morphism(0,-\ptless\dmB)to(\ptless\dmA,-\ptless\dmE)[\ss{#6}][R]
 \morphism(#1,-\ptless\dmB)to(\ptless\dmA,-\ptless\dmE)[\ss{#7}]
}

\def\hexagondiagabscoord#1#2#3#4#5#6#7#8{%
 \def\firstobj{#3}\def\secondobj{#4}\def\thirdobj{#5}%
 \def\fourthobj{#6}\def\fifthobj{#7}\def\sixthobj{#8}%
 \argarghgda{#1}{#2}}

\def\argarghgda#1#2#3#4#5#6#7#8{%
  \let\ss=\scriptstyle
  \latticeUnit=1pt \boxMargin=3pt
  \hgda{#1}{#2}{#3}{#4}{#5}{#6}{#7}{#8}%
}

\def\hgda#1#2#3#4#5#6#7#8{%
 \dmA=#1pt\dmB=#2pt
 \multiply\dmB by58\divide\dmB by100
 \dmC=\dmA \divide\dmC by2
 \dmD=\dmB \divide\dmD by2
 \object(\ptless\dmC,\ptless\dmB)={\firstobj}
 \object(0,\ptless\dmD)={\secondobj}
 \object(#1,\ptless\dmD)={\thirdobj}
 \object(0,-\ptless\dmD)={\fourthobj}
 \object(#1,-\ptless\dmD)={\fifthobj}
 \object(\ptless\dmC,-\ptless\dmB)={\sixthobj}
 \morphism(\ptless\dmC,\ptless\dmB)to(0,\ptless\dmD)[\ss{#3}][R]
 \morphism(\ptless\dmC,\ptless\dmB)to(#1,\ptless\dmD)[\ss{#4}]
 \morphism(0,\ptless\dmD)to(0,-\ptless\dmD)[\ss{#5}][R]
 \morphism(#1,\ptless\dmD)to(#1,-\ptless\dmD)[\ss{#6}]
 \morphism(0,-\ptless\dmD)to(\ptless\dmC,-\ptless\dmB)[\ss{#7}][R]
 \morphism(#1,-\ptless\dmD)to(\ptless\dmC,-\ptless\dmB)[{#8}]
}

\def\putFormula[#1](#2,#3)#4{%
 \vspec{gsave 1 1 scale #2 8.3 mul #3 -8.3 mul translate}%
 \if L#1$\vbox to0pt{\vss\llap{$#4$}\vss}$%
 \else\if R#1$\vbox to0pt{\vss\rlap{$#4$}\vss}$%
 \else $\vbox to0pt{\vss\hbox to0pt{\hss $#4$\hss}\vss}$%
 \fi\fi
 \vspec{grestore}%
}

\def\putSpearHead(#1,#2)(#3,#4){%
 \vspec{gsave
  currentpoint currentpoint translate #3 8.3 mul #4 -8.3 mul translate
  #3 #1 sub #4 #2 sub atan rotate
  neg exch neg exch translate}%
 \tip%
 \vspec{grestore}}

\def\red{\special{ps:1 1 1 sethsbcolor}}%
\def\black{
}%

\newcount\cntDwn

\def\pointNum#1{%
 \putFormula[C](x#1,y#1){\scriptscriptstyle \bf \red#1\black}}%

\def\dwnNext{%
 \ifnum\the\cntDwn=0
 \else
  \pointNum{\the\cntDwn}%
  \advance\cntDwn by-1
  \dwnNext
 \fi
}%

\def\putTensor(#1,#2){%
 \begin{scope}[shift={(#1,#2)}]
  \fill[black] (0,0) circle (5);
  \fill[white] (0,0) circle (4.3);
  \draw[line width=.65] (3.5, 3.5) -- (-3.5, -3.5);
  \draw[line width=.65] (-3.5, 3.5) -- (3.5, -3.5);
 \end{scope}
}

\def\putPar(#1,#2){%
 \begin{scope}[shift={(#1,#2)}]
  \fill[black] (0,0) circle (5);
  \fill[white] (0,0) circle (4.3);
  \fill[black] (0,0) circle (2.2);
  \fill[white] (0,0) circle (1.5);
 \end{scope}
}

\def\putRightDiode(#1,#2){%
 \begin{scope}[thin]
  \fill (#1+1.5, #2+0) -- (#1-2.5, #2+2.5) -- (#1-2.5, #2-2.5) -- cycle;
  \fill (#1+1.5, #2+2.5) -- (#1+1.5, #2-2.5) -- (#1+2.5, #2-2.5) -- (#1+2.5, #2+2.5) -- cycle;
 \end{scope}
}

\def\putLeftDiode(#1,#2){%
 \begin{scope}[thin]
  \fill (#1-1.5, #2+0) -- (#1+2.5, #2+2.5) -- (#1+2.5, #2-2.5) -- cycle;
  \fill (#1-1.5, #2+2.5) -- (#1-1.5, #2-2.5) -- (#1-2.5, #2-2.5) -- (#1-2.5, #2+2.5) -- cycle;
 \end{scope}
}

\def\putPositiveTerminal(#1,#2){%
 \begin{scope}
  \draw[line width=.7] (#1-4,#2) -- (#1+4,#2);
 \end{scope}
}

\def\putUpperNegativeTerminal(#1,#2){%
 \begin{scope}
  \draw[line width=.7] (#1,#2) arc (270:340:4);
  \draw[line width=.7] (#1,#2) arc (270:200:4);
 \end{scope}
}

\def\putLowerNegativeTerminal(#1,#2){%
 \begin{scope}
  \draw[line width=.7] (#1,#2) arc (90:20:4);
  \draw[line width=.7] (#1,#2) arc (90:160:4);
 \end{scope}
}

\def\putRubberBand(#1,#2){
 \begin{scope}
  \draw[white,line width=2.5] (#1-2,#2-1.5) -- (#1+2,#2-1.5);
  \draw[black,line width=.7] (#1-1.02,#2+1.41) arc (110:430:3 and 1.5);
 \end{scope}
}

\def\putConvex(#1,#2){
 \begin{scope}[line width=.7]
 \fill[white] (#1,#2) -- ++(-6,0) -- ++(0,-2) arc (233:307:10) -- ++(0,2) -- cycle;
 \draw (#1,#2) -- ++(-6,0) -- ++(0,-2) arc (233:307:10) -- ++(0,2) -- cycle;
 \end{scope}
}

\def\putConcave(#1,#2){
 \begin{scope}[line width=.7]
 \fill[white] (#1,#2) -- ++(6,0) -- ++(0,-4) arc (53:127:10) -- ++(0,4) -- cycle;
 \draw (#1,#2) -- ++(6,0) -- ++(0,-4) arc (53:127:10) -- ++(0,4) -- cycle;
 \end{scope}
}

\def\putEliminator(#1,#2){
 \begin{scope}[line width=.7]
  \fill[white] (#1,#2) -- ++(6,0) -- ++(-6,-6) -- ++(-6,6) -- cycle;
  \draw (#1,#2) -- ++(6,0) -- ++(-6,-6) -- ++(-6,6) -- cycle;
 \end{scope}
}

\def\putDuplicator(#1,#2){
 \begin{scope}[line width=.7]
  \draw (#1,#2) -- ++(6,0);
  \draw (#1,#2) -- ++(-6,0);
  \draw (#1+2.6,#2) -- ++(3.5,-6.0);
  \draw (#1-2.6,#2) -- ++(-3.5,-6.0);
 \end{scope}
}

\def\putUpperSocket(#1,#2){
 \begin{scope}
  \fill[white] (#1,#2) circle (5);
  \draw[line width=.7] (#1,#2) circle (5);
  \fill[line width=.1] (#1,#2+4) .. controls (#1-1.5,#2+1.2) .. ++(-3.0,-4.6)
   -- ++(0.2,-0.2) -- ++(2.8,2.2) -- ++(2.8,-2.2) -- ++(0.2,0.2)
   .. controls (#1+1.5,#2+1.2) .. ++(-3.0,4.6) -- cycle;
  \fill[line width=.1] (#1-0.35,#2-4.0) -- ++(0.0,7.0)
    -- ++(0.7,0.0) -- ++ (0.0,-7.0) -- cycle;
 \end{scope}
}

\def\putLowerSocket(#1,#2){
 \begin{scope}
  \fill[white] (#1,#2) circle (5);
  \draw[line width=.7] (#1,#2) circle (5);
  \fill[line width=.1] (#1,#2-4) .. controls (#1-1.5,#2-1.2) .. ++(-3.0,4.6)
   -- ++(0.2,0.2) -- ++(2.8,-2.2) -- ++(2.8,2.2) -- ++(0.2,-0.2)
   .. controls (#1+1.5,#2-1.2) .. ++(-3.0,-4.6) -- cycle;
  \fill[line width=.1] (#1-0.35,#2+4.0) -- ++(0.0,-7.0)
    -- ++(0.7,0.0) -- ++ (0.0,7.0) -- cycle;
 \end{scope}
}

\def\putDimple(#1,#2){
 \begin{scope}
  \fill[white] (#1,#2) circle (4.6);
  \draw[line width=.9] (#1-5,#2) arc (180:360:5);
 \end{scope}
}

\def\boardBoundary(#1,#2)(#3,#4){
 \begin{scope}[line width=.9]
  \draw (#1,#2) rectangle (#3,#4);
 \end{scope}
}

\long\def\figbox#1{\vbox{\kern3pt\hbox{\vrule \vbox{\hrule
 \hbox{\vbox{\leftskip2em\rightskip2em\advance\hsize by-.8pt
 #1}}\hrule }\vrule}\kern5pt}}

\newbox\tcd\setbox\tcd=\hbox{$\&$}
\def\linpar{\mathbin{\tikz[inner sep=0] \node[rotate=180] at (0,0) {$\&$};}}
\def\slinpar{\mathbin{\tikz[inner sep=0] \node[rotate=180] at (0,0) {$\scriptstyle\&$};}}

\def\scirc{\mathbin{\vcenter{\hbox{\scriptsize $\circ$}}}}

\def\llapem#1#2{\llap{\hbox to#1em{\rm #2\hss}}}
\let\tempar\par \def\par{{\tempar}}%

\def\bang{\mathord!}

\def\di#1(#2,#3)#4{%
 \dimA=#2\latticeUnit
 \dimB=#3\latticeUnit
 \put(\ptless\dimA,\ptless\dimB){\tikz[inner sep=0]
   \node[rotate=-#4] at (\ptless\dimA,\ptless\dimB)
     {\vbox to0pt{\vss\hbox to0pt{\hss #1\hss}\vss}};}}

\makeatother

\begin{document}

\addtolength{\headsep}{10pt}

\maketitle

\begin{abstract}
We verify a confluence result for the rewriting
 calculus of the linear category introduced in our previous paper.
Together with the termination result proved therein,
 the generalized coherence
 theorem for linear category
 is established.
Namely, we obtain a method to determine if two morphisms
 are equal up to a certain equivalence.
\end{abstract}

\section{Introduction}

The link between the type theory and the category theory is
 bidirectional.
The type theory is a theoretical framework
 for the abstract formalization
 of programming languages,
 mainly developed in computer science.
It is known that categories are suitable
 machinery to  elucidate the mathematical structure
 of type systems, providing a large pool
 of the mathematical models of the artificial
 languages designed for programming purposes.
Conversely, the type theory is employed as the internal
 languages of categories, replacing the lengthy diagram chasing with
 intuitive arguments using familiar logical
 constructions \cite{jaco}.
In various cases, categories give the sound and complete
 semantics of type systems.
This means that the equational theory determined by
 commutative diagrams in the categories exactly correspond
 to the one determined by the equality rules between terms
 of the type systems.
Namely, categories and type systems are not
 quite  similar in their appearances, but
 they are equivalent in their essences.
For example, it is well-known that cartesian closed categories
 give the sound and complete semantics of the simply typed lambda
 calculus with products \cite{lasc}.
In this paper, we consider the categorical semantics of
 linear logic.
It is also sound and complete.

Dynamism is lost in the link between two theories.
The type theory is developed as a mathematical formalization of computation
 on programs.
Naturally, thus, most type systems possess the mechanism of calculi.
The simply typed lambda calculus is typical.
The notion of reductions is incorporated, by which we can
 perform mechanical calculations.
This calculus satisfies the most desirable properties demanded
 to formalized calculi:
 termination and confluence.
A sequence of reductions is assured of terminating
 in a normal form that allows reductions no more.
As the calculus admits non-determinism in the order of reductions,
 different routes of computation may exist.
The confluence ensures that the normal form is unique no matter
 which routes are taken.
Unfortunately, however, the categorical semantics so far fails to capture
 the dynamic aspect of the calculus.
Equivalence between categories and type systems holds
 only up to equalities.
We must neglect the orientation of reductions to validate
 the correspondence.
There is no obvious way to transfer the reductions
 in type systems into categories.
For example, although we are not unable to introduce reductions
 on the free cartesian closed category, the obtained
 calculus is far from being as good as the lambda calculus
 with regard to computational properties.

The main finding in our previous paper is that
 we can install a natural, good calculus on the
 categorical model of linear logic,
 the free linear category \cite{hase2}.
Linear logic is a refinement of the lambda calculus.
What we have disclosed is that
 the cartesian closed category corresponding to the lambda
 calculus is too coarse.
The fine granularity of linear logic enables
 us to turn the linear category into a calculus.
We have developed a rewriting system modulo congruence
 by choosing twenty-three from
 the defining commutative diagrams and making
 them the rewriting rules.
One can rewrite morphisms only in the direction specified
 by the rules.
In the previous paper, we verified a termination
 property for the calculus.

The theme of this paper is confluence, the other half
 of the two properties naturally demanded to reasonable computational
 systems.
We show that our rewriting system is almost confluent.
The exact meaning of being almost confluent is explained
 in a later section.
In brief, a normal form is unique, except
 how and where the isomorphisms related to tensor/cotensor
 units are used.
The (classical) linear category is based
 on the $*$-autonomous category.
The manipulation of the units in the $*$-autonomous category
 is notoriously difficult \cite{bcst,hugh,heho}.
Speaking a little audaciously, however, the problem of the units
 is relatively a minor point.
Hence it may be allowable to say that our calculus is almost confluent.

We provide a semantical proof of confluence.
The parallel reduction is known as
 a standard syntactic method to verify confluence.
As the name suggests, it takes
 all possible parallel combinations of reductions
 as single-step rewritings.
As our calculus has twenty-three rules, however,
 the taking of all combinations is daunting.
We should find other means.
In this paper, we provide a semantical method for confluence.
The link between the syntax and the semantics is usually
 directional.
Models are used to abstract the properties of
 the syntax to obtain semantical information.
However, the other direction is occasionally possible.
Syntactic information can be squeezed out from models.
The evaluation-free normalization by Berger and Schwichtenberg
 extracts a normal form of a term of the typed lambda calculus
 from a model \cite{besh}.
The coherence proof by Joyal and Street finds canonical
 morphisms from concrete categories \cite{jost}.
A series of results for differential nets by de~Carvalho
 and Tortora~de~Falco et al.~is based on a similar idea \cite{carv1,cato,gptf}.
We use the normal functor model to obtain the information
 of normal forms in our calculus.

Remark:
After completing the current work, we noticed the work by
 de~Carvalho \cite{carv2}.
It shares similar ideas with ours, and
 it is plausible that his method can be applied to our calculus.
However, we hope that this paper still
 has some values, since (1) the systems are different, (2)
 the models are different, (3) our method relying on
 the enumerative combinatorics and the number theory may have
 novelty.
We use the linear normal functors, the coefficients of which
 are non-negative integers.
Therefore we can apply the number theory to the coefficients.
In our proof, Fermat's little theorem plays a key role.

We can establish a kind of coherence for the linear category.
The classic coherence theorem by Mac Lane ensures
 that the diagram chasing in certain categories, such as
 the monoidal category, is trivial \cite{kell}.
Any morphisms sharing a domain and a codomain
 must equal.
Most categories, though, fail in having coherence.
However, we occasionally have a certain type of
 characterization of morphisms, which may be regarded as
 the generalization of coherence.
The coherence theorem for the braided monoidal category by
 Joyal and Street relates the morphisms
 to standard braid diagrams \cite{jost}.
Another direction to generalize coherence is to seek an effective
 method to determine if given morphisms are equal.
Blute et al.~have shown that the $*$-autonomous category
 satisfies the generalized coherence in this sense \cite{bcst}.
There is a mechanical way to check if two morphisms equal, freeing us
 from the cumbersome task of forming large commutative diagrams.
The confluence result in this paper, accompanied with a
 termination result in our previous paper, entails that
 the linear category fulfills the generalized coherence
 in the latter sense.
We have a systematic way to determine whether two morphisms
 are almost equal.
Here ``almost equal'' means that
 they are equal
 if we ignore where and how the isomorphisms related
 to units are used.
Our algorithm is simple.
Given two morphisms, turn them to normal forms and
 just compare.

\section{Rewriting system for the linear category}

In our previous paper, we have reformulated the free linear category
 into a rewriting system so that the diagram chasing is
 realized by essentially one-way sequence of computations \cite{hase2}.
A (classical) {\it linear category} is a $*$-autonomous
 category $({\bf C},\otimes,\linpar,{\bf 1},\bot,(\hbox{-})^*)$
 having the following additional data \cite{mmpr}:

\begingroup
	\vskip.5ex
        \hangafter0\hangindent2em
        \noindent
\llapem2{(i)}%
The category ${\bf C}$ is equipped with a symmetric
 monoidal endofunctor $(\mathord!,\tilde\varphi,\varphi_0)$.

	\vskip0ex
        \noindent
\llapem2{(ii)}%
The functor $\mathord!$ has the structure of comonad
 $(\mathord!,\delta,\varepsilon)$ where
 $\delta:\mathord!\rightarrow \mathord!\mathord!$ and
 $\varepsilon:\mathord!\rightarrow{\it Id}$.
Moreover, $\delta$ and $\varepsilon$ are monoidal natural
 transformations.

	\vskip0ex
        \noindent
\llapem2{(iii)}%
Each object of the form $\mathord!A$ has the structure
 of commutative comonoid $(\mathord!A,d_A,e_A)$ where
 $d_A:\mathord!A\rightarrow \mathord!A\otimes \mathord!A$
 and $e_A:\mathord!A\rightarrow {\bf 1}$.
Moreover, the family of $d_A$ and the family of $e_A$ are collectively
 monoidal natural transformations.

	\vskip.5ex
        \hangafter0\hangindent0em
        \noindent
These data should fulfill the following constraints:

	\vskip.5ex
        \hangafter0\hangindent2em
        \noindent
\llapem2{(iv)}%
Each $d_A$ and each $e_A$ are coalgebra morphisms.

	\vskip0ex
        \noindent
\llapem2{(v)}%
Each $\delta_A$ is a comonoid morphism.

	\vskip.5ex
        \hangafter0\hangindent0em
\endgroup
	\noindent
As the formalization of the underlying $*$-autonomous category,
 we take the one having two symmetric monoidal structures
 $(\otimes,{\bf 1})$ and $(\linpar,\bot)$ with
 linear distribution $\partial:A\otimes (B\linpar C)
 \rightarrow (A\otimes B)\linpar C$, tautology $\tau:{\bf 1}\rightarrow
 A\linpar A^*$, and contradiction
 $\gamma:A^*\otimes A\rightarrow\bot$ \cite{cose}.

The definition can be written down as a number
 of commutative diagrams.
We reformulate twenty-three diagrams among them
 as rewriting rules.
In the rules (18) to (21) below, $f$ denotes an arbitrary morphism.

	\vskip2ex
\begingroup\offinterlineskip
\halign to\textwidth{\kern50pt\hfil #\hfil \tabskip1000pt plus100pt minus1000pt
 &\kern10pt$\cdots$\kern5pt{\footnotesize #}\kern20pt\tabskip0pt
 &\kern20pt\hfil #\hfil \tabskip1000pt plus100pt minus1000pt
 &\kern10pt$\cdots$\kern5pt{\footnotesize #}\tabskip0pt\cr
$\vcenter{\footnotesize
\hbox{\begin{diagramme}
\squarediagabscoord{50}{40}%
 {\mathord!A}{\mathord!\mathord!A}{\mathord!\mathord!A}%
 {\mathord!\mathord!\mathord!A}%
 {\delta}{\delta}{\delta}{\mathord!\delta}
 \di{\hbox{$\Rightarrow$}}(25,20){140}
\end{diagramme}}}$%
&(1)%
&$\vcenter{\footnotesize
\hbox{\begin{diagramme}
\trianglediagabscoord{50}{30}%
 {\mathord!A}{\mathord!\mathord!A}{\mathord!A}%
 {\delta}{1}{\varepsilon}%
 \di{\hbox{$\Rightarrow$}}(25,18){140}
\end{diagramme}}}$%
&(2)%
\cr
	\noalign{\vskip5ex}
$\vcenter{\footnotesize
\hbox{\begin{diagramme}
\trianglediagabscoord{50}{30}%
 {\mathord!A}{\mathord!\mathord!A}{\mathord!A}%
 {\delta}{1}{\mathord!\varepsilon}%
 \di{\hbox{$\Rightarrow$}}(25,18){140}
\end{diagramme}}}$%
&(3)%
&$\vcenter{\footnotesize
\hbox{\begin{diagramme}
 \squarediagabscoord{50}{40}%
 {\mathord!A}{\mathord!\mathord!A}%
 {\kern-5pt\mathord!A\otimes \mathord!A}%
 {\mathord!\mathord!A\otimes \mathord!\mathord!A\kern-5pt}%
 {\delta}{d}{d}{\delta \otimes \delta}
 \di{\hbox{$\Rightarrow$}}(25,20){140}
\end{diagramme}}}$%
&(4)%
\cr
	\noalign{\vskip5ex}
$\kern-9pt\vcenter{\footnotesize
\hbox{\begin{diagramme}
\pentagondiagabscoord{70}{60}%
 {\mathord!A}{\mathord!A\otimes \mathord!A}{\mathord!\mathord!A}%
 {\kern-20pt\mathord!\mathord!A\otimes \mathord!\mathord!A}%
 {\mathord!(\mathord!A\otimes \mathord!A)\kern-20pt}
 {d}{\delta}{\delta\otimes\delta}{\mathord!d}{\tilde\varphi} 
 \di{\hbox{$\Rightarrow$}}(35,25){140}
\end{diagramme}}}$%
&(5)%
&$\vcenter{\footnotesize
\hbox{\begin{diagramme}
\trianglediagabscoord{50}{30}%
 {\mathord!A}{\mathord!\mathord!A}{{\bf 1}}%
 {\delta}{e}{e}%
 \di{\hbox{$\Rightarrow$}}(25,18){140}
\end{diagramme}}}$%
&(6)%
\cr
	\noalign{\vskip5ex}
$\vcenter{\footnotesize
\hbox{\begin{diagramme}
\squarediagabscoord{50}{40}%
 {\mathord!A}{\mathord!\mathord!A}{\bf 1}{\mathord!{\bf 1}}%
 {\delta}{e}{\mathord!e}{\varphi_0}
 \di{\hbox{$\Rightarrow$}}(25,20){140}
\end{diagramme}}}$%
&(7)%
&$\kern11pt\vcenter{\footnotesize
\hbox{\begin{diagramme}
\trianglediagabscoord{50}{30}%
 {\mathord!A}{\mathord!A\otimes \mathord!A}{{\bf 1}\otimes \mathord!A}%
 {d}{\sim}{e\otimes 1}
 \di{\hbox{$\Rightarrow$}}(25,18){140}
\end{diagramme}}}$%
&(8)%
\cr
	\noalign{\vskip5ex}
$\vcenter{\footnotesize
\hbox{\begin{diagramme}
\pentagondiagabscoord{70}{60}%
 {\mathord!A\otimes\mathord!B}%
 {\mathord!\mathord!A\otimes \mathord!\mathord!B}%
 {\mathord!(A\otimes B)}%
 {\kern-20pt\mathord!(\mathord!A\otimes \mathord!B)}%
 {\mathord!\mathord!(A\otimes B)\kern-20pt}%
  {\delta\otimes\delta}{\tilde\varphi}{\tilde\varphi}
  {\delta}{\mathord!\tilde\varphi}%
 \di{\hbox{$\Rightarrow$}}(35,25){140}
\end{diagramme}}}$%
&(9)%
&$\vcenter{\footnotesize
\hbox{\begin{diagramme}
\trianglediagabscoord{50}{30}%
 {\kern-10pt\mathord!A\otimes \mathord!B}%
 {\mathord!(A\otimes B)\kern-10pt}{A\otimes B}%
 {\tilde\varphi}{\varepsilon\otimes\varepsilon}{\varepsilon}%
 \di{\hbox{$\Rightarrow$}}(25,18){140}
\end{diagramme}}}$%
&(10)%
\cr
	\noalign{\vskip5ex}
$\kern-8pt\vcenter{\footnotesize
\hbox{\begin{diagramme}
\pentagondiagabscoord{70}{60}%
 {\mathord!A\otimes \mathord!B}%
 {\kern-40pt
    (\mathord!A\otimes \mathord!A)\otimes (\mathord!B\otimes \mathord!B)}%
 {\mathord!(A\otimes B)}%
 {\kern-70pt
    (\mathord!A\otimes \mathord!B)\otimes (\mathord!A\otimes \mathord!B)}%
 {\mathord!(A\otimes B)\otimes \mathord!(A\otimes B)\kern-70pt}%
  {d\otimes d}{\tilde\varphi}{\sim}%
  {d}{\tilde\varphi\otimes\tilde\varphi}
 \di{\hbox{$\Rightarrow$}}(35,25){140}
\end{diagramme}}}$%
&(11)%
&$\kern5pt\vcenter{\footnotesize
\hbox{\begin{diagramme}
\squarediagabscoord{50}{40}%
  {\kern-5pt\mathord!A\otimes\mathord!B}%
  {\mathord!(A\otimes B)\kern-5pt}%
  {{\bf 1}\otimes{\bf 1}}{{\bf 1}}%
  {\tilde\varphi}{e\otimes e}{e}{\sim}
 \di{\hbox{$\Rightarrow$}}(25,20){140}
\end{diagramme}}}$%
&(12)%
\cr
	\noalign{\vskip5ex}
$\vcenter{\footnotesize
\hbox{\begin{diagramme}
\squarediagabscoord{50}{40}%
 {{\bf 1}}{\mathord!{\bf 1}}{\mathord!{\bf 1}}{\mathord!\mathord!{\bf 1}}%
 {\varphi_0}{\varphi_0}{\delta}{\mathord!\varphi_0}
 \di{\hbox{$\Rightarrow$}}(25,20){140}
\end{diagramme}}}$%
&(13)%
&$\vcenter{\footnotesize
\hbox{\begin{diagramme}
\trianglediagabscoord{50}{30}%
 {{\bf 1}}{\mathord!{\bf 1}}{\bf 1}%
 {\varphi_0}{1}{\varepsilon}
 \di{\hbox{$\Rightarrow$}}(25,18){140}
\end{diagramme}}}$%
&(14)%
\cr
	\noalign{\vskip5ex}
$\vcenter{\footnotesize
\hbox{\begin{diagramme}
\squarediagabscoord{50}{40}%
 {{\bf 1}}{\mathord!{\bf 1}}{{\bf 1}\otimes{\bf 1}}%
 {\bang{\bf 1}\otimes\bang{\bf 1}}
 {\varphi_0}{\sim}{d}{\varphi_0\otimes\varphi_0}
 \di{\hbox{$\Rightarrow$}}(25,20){140}
\end{diagramme}}}$%
&(15)%
&$\vcenter{\footnotesize
\hbox{\begin{diagramme}
\trianglediagabscoord{50}{30}%
  {{\bf 1}}{\mathord!{\bf 1}}{{\bf 1}}%
  {\varphi_0}{1}{e}
 \di{\hbox{$\Rightarrow$}}(25,18){140}
\end{diagramme}}}$%
&(16)%
\cr
	\noalign{\vskip5ex}
$\kern4pt\vcenter{\footnotesize
\hbox{\begin{diagramme}
\squarediagabscoord{50}{40}%
 {{\bf 1}\otimes\bang A}{\bang{\bf 1}\otimes\bang A}%
 {\bang A}{\bang({\bf 1}\otimes A)}%
 {\varphi_0\otimes 1}{\sim}{\tilde\varphi}{\sim}
 \di{\hbox{$\Rightarrow$}}(25,20){140}
\end{diagramme}}}$%
&(17)%
&$\vcenter{\footnotesize
\hbox{\begin{diagramme}
\squarediagabscoord{50}{40}%
 {\mathord!A}{\mathord!B}{\mathord!\mathord!A}{\mathord!\mathord!B}%
 {\mathord!f}{\delta}{\delta}{\mathord!\mathord!f}
 \di{\hbox{$\Rightarrow$}}(25,20){140}
\end{diagramme}}}$%
&(18)%
\cr
	\noalign{\vskip5ex}
$\vcenter{\footnotesize
\hbox{\begin{diagramme}
\squarediagabscoord{50}{40}%
 {\mathord!A}{\mathord!B}{A}{B}%
 {\mathord!f}{\varepsilon}{\varepsilon}{f}
 \di{\hbox{$\Rightarrow$}}(25,20){140}
\end{diagramme}}}$%
&(19)%
&$\vcenter{\footnotesize
\hbox{\begin{diagramme}
\squarediagabscoord{50}{40}%
 {\mathord!A}{\mathord!B}{\kern-10pt\mathord!A\otimes \mathord!A}%
 {\mathord!B\otimes \mathord!B\kern-10pt}%
 {\mathord!f}{d}{d}{\mathord!f\otimes \mathord!f}
 \di{\hbox{$\Rightarrow$}}(25,20){140}
\end{diagramme}}}$%
&(20)%
\cr
	\noalign{\vskip5ex}
$\vcenter{\footnotesize
\hbox{\begin{diagramme}
\trianglediagabscoord{50}{30}%
 {\mathord!A}{\mathord!B}{{\bf 1}}%
 {\mathord!f}{e}{e}
 \di{\hbox{$\Rightarrow$}}(25,18){140}
\end{diagramme}}}$%
&(21)%
&$\kern30pt\vcenter{\footnotesize
\def\linpar{\mathbin{\tikz[inner sep=0] \node[rotate=180] at (0,0) {\footnotesize $\&$};}}
\def\slinpar{\mathbin{\tikz[inner sep=0] \node[rotate=180] at (0,0) {\tiny $\&$};}}
\hbox{\begin{diagramme}
\oppentagondiagabscoord{70}{60}%
 {\kern-5pt{\bf 1}\otimes A}{(A\linpar A^*)\otimes A\kern-35pt}{A}%
 {A\linpar (A^*\otimes A)\kern-15pt}{A\linpar \bot}%
 {\tau\otimes 1}{\sim}{\partial'}{\sim}{1\slinpar\gamma}
 \di{\hbox{$\Rightarrow$}}(35,-25){140}
\end{diagramme}}}\kern15pt$%
&(22)%
\cr
$\kern15pt\vcenter{\footnotesize
\def\linpar{\mathbin{\tikz[inner sep=0] \node[rotate=180] at (0,0) {\footnotesize $\&$};}}
\def\slinpar{\mathbin{\tikz[inner sep=0] \node[rotate=180] at (0,0) {\tiny $\&$};}}
\hbox{\begin{diagramme}
\oppentagondiagabscoord{70}{60}%
 {\kern-5ptA^*\otimes{\bf 1}}{A^*\otimes(A\linpar A^*)\kern-40pt}{A^*}%
 {(A^*\otimes A)\linpar A^*\kern-15pt}{\bot\linpar A^*}%
 {1\otimes\tau}{\sim}{\partial}{\sim}{\gamma\slinpar 1}
 \di{\hbox{$\Rightarrow$}}(35,-25){140}
\end{diagramme}}}\kern0pt$%
&(23)%
\cr}
\endgroup

	\vskip3ex
        \noindent
The double arrow $\Rightarrow$ means the rewriting relation.
For example, let us consider rule (1), which is one
 of the coherence conditions for comonad.
Usually, it is interpreted as a commutative diagram.
Namely, $\delta_A;\delta_{\mathord!A}$ and $\delta_A;\mathord!\delta_A$
 are regarded to be equal.
We can replace the former with the latter, or the latter with the former.
Under the interpretation as a rewriting relation, we allow only one-way
 modification.
We can replace
 $\delta_A;\delta_{\mathord!A}$ with $\delta_A;\mathord!\delta_A$, not
 the other way round.
In this way, the modification of morphisms is constrained.

\begin{figure}
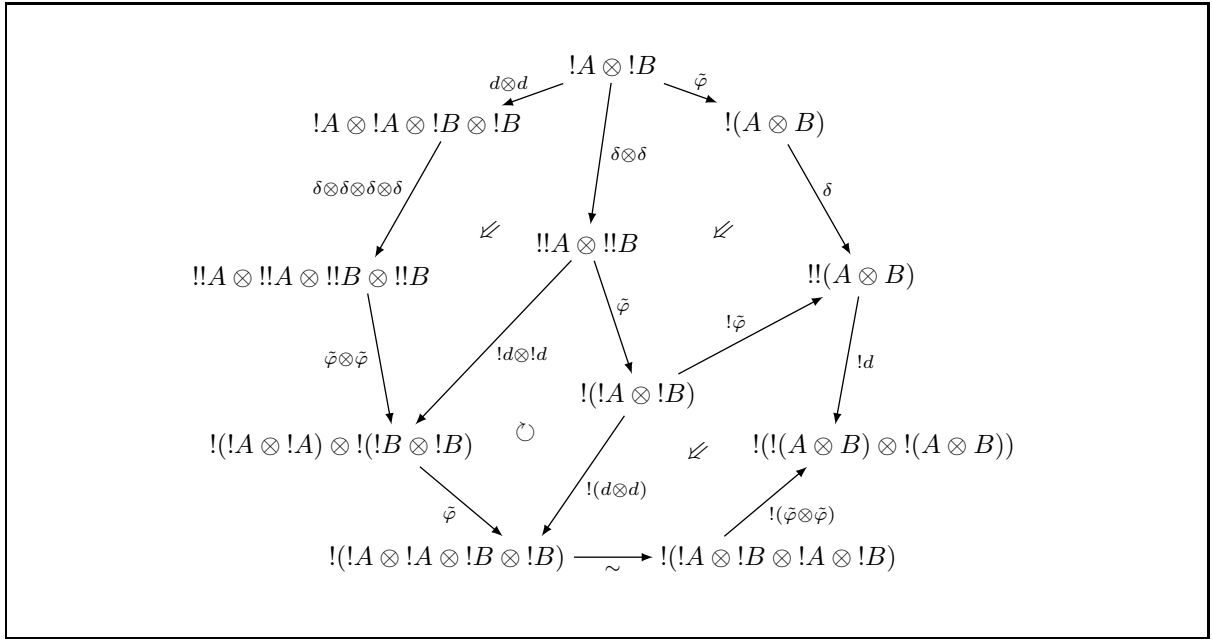

\figbox{%
	\vskip3ex
$$\kern0pt\vcenter{\footnotesize
\hbox{\begin{diagramme}
 \let\ss=\scriptstyle
 \latticeUnit=.96pt \boxMargin=3pt
 \object(0,100)={\bang A\otimes \bang B}
 \object(-64,77)={\kern-25pt\bang A\otimes \bang A
    \otimes\bang B\otimes \bang B}
 \object(64,77)={\bang(A\otimes B)}
 \object(-98,17)={\kern-40pt\bang\bang A\otimes \bang\bang A
    \otimes\bang\bang B\otimes \bang\bang B}
 \object(98,17)={\bang\bang(A\otimes B)}
 \object(-86,-50)={\kern-40pt\bang(\bang A\otimes \bang A)
    \otimes\bang(\bang B\otimes \bang B)}
 \object(86,-50)={\bang(\bang (A\otimes B)
    \otimes \bang (A\otimes B))\kern-40pt}
 \object(-34,-94)={\kern-60pt\bang(\bang A\otimes\bang A
    \otimes\bang B\otimes\bang B)}
 \object(34,-94)={\bang(\bang A\otimes\bang B
    \otimes\bang A\otimes\bang B)\kern-60pt}
 \object(-10,30)={\bang\bang A\otimes\bang\bang B}
 \object(10,-30)={\bang(\bang A\otimes \bang B)}
 \morphism(0,100)to(-64,77)[d\otimes d][R]
 \morphism(0,100)to(64,77)[\tilde\varphi]
 \morphism(-64,77)to(-98,17)[\delta\otimes\delta\otimes\delta\otimes\delta][R]
 \morphism(64,77)to(98,17)[\delta]
 \morphism(-98,17)to(-86,-50)[\tilde\varphi\otimes\tilde\varphi][R]
 \morphism(98,17)to(86,-50)[\bang d]
 \morphism(-86,-50)to(-34,-94)[\tilde\varphi][R]
 \morphism(34,-94)to(86,-50)[\bang(\tilde\varphi\otimes\tilde\varphi)][R]
 \morphism(-34,-94)to(34,-94)[\sim][R]
 \morphism(0,100)to(-10,30)[\delta\otimes\delta]
 \morphism(-10,30)to(10,-30)[\tilde\varphi]
 \morphism(-10,30)to(-86,-50)[\bang d\otimes\bang d]
 \morphism(10,-30)to(98,17)[\bang\tilde\varphi]
 \morphism(10,-30)to(-34,-94)[\bang(d\otimes d)]
 \di{\hbox{$\Rightarrow$}}(42,35){140}
 \di{\hbox{$\Rightarrow$}}(32,-52){140}
 \di{\hbox{$\circlearrowright$}}(-35,-45){0}
 \di{\hbox{$\Rightarrow$}}(-50,35){140}
\end{diagramme}}}$$
	\vskip3ex
}
	\vskip-1ex
\caption{An example of rewriting}\label{pss71}
\end{figure}

There are still a number of diagrams not listed above
in the definition of the linear category.
They are regarded as commutative diagrams in an ordinary sense.
Two legs of a commutative diagram, say $f$ and $g$, are regarded to be equal.
We can replace either $f$ with $g$, or $g$ with $f$.
Hence our calculus is actually a rewriting system
 modulo congruence.
The diagram in Fig.~\ref{pss71} is an example of rewriting in our calculus,
 reproduced from our previous paper.
It contains rewritings by (5), (9), and (11), and one congruence.
All the diagrams defining the linearly distributive category
 is regarded to be congruent.
The diagrams related to $(\hbox{-})^*$ are handled
 as rewriting rules (22) and (23).
For the symmetric monoidal functor $\mathord!$, only one
 diagram (17) gives rise to a rewriting relation.
The others give congruent relations:

	\vskip2ex
\begingroup\offinterlineskip
\halign to\textwidth{\kern120pt\hfil #\hfil \tabskip1000pt plus100pt minus1000pt
 &\kern30pt\hfil #\hfil \kern120pt\tabskip0pt
\cr
$\vcenter{\footnotesize
\hbox to0pt{\hss\begin{diagramme}
\squarediagabscoord{50}{40}%
  {\kern-5pt\mathord!A\otimes \mathord!B}{\mathord!(A\otimes B)\kern-5pt}
  {\kern-10pt\mathord!A'\otimes \mathord!B'}{\mathord!(A'\otimes B')\kern-10pt}
  {\tilde\varphi}{\mathord!f\otimes\mathord!g}%
  {\mathord!(f\otimes g)}{\tilde\varphi}
 \di{\hbox{$\circlearrowright$}}(25,20){0}
\end{diagramme}\hss}}$\kern-4pt
\cr
	\noalign{\vskip4ex}
$\vcenter{\footnotesize
 \hbox to0pt{\hss\begin{diagramme}
  \let\ss=\scriptstyle
  \latticeUnit=1pt \boxMargin=3pt
  \hexagondiagabscoord{80}{70}%
    {(\mathord!A\otimes \mathord!B)\otimes \mathord!C}%
    {\kern-10pt\mathord!A\otimes (\mathord!B\otimes \mathord!C)}%
    {\mathord!(A\otimes B)\otimes \mathord!C\kern-10pt}%
    {\kern-10pt\mathord!A\otimes \mathord!(B\otimes C)}
    {\mathord!((A\otimes B)\otimes C)\kern-10pt}%
    {\mathord!(A\otimes (B\otimes C))}%
     {\global\labelPosition={.9}\alpha}{\tilde\varphi\otimes 1}%
     {\global\labelPosition={.5}1\otimes \tilde\varphi}%
     {\tilde\varphi}%
     {\global\labelPosition={-.2}\tilde\varphi}{\mathord!\alpha}%
 \global\labelPosition={.5}%
 \di{\hbox{$\circlearrowright$}}(40,0){0}%
 \end{diagramme}\kern-4pt\hss}}$
&$\vcenter{\footnotesize
\hbox to0pt{\hss\begin{diagramme}
\squarediagabscoord{50}{40}%
  {\kern-5pt\mathord!A\otimes \mathord!B}{\mathord!(A\otimes B)\kern-5pt}
  {\kern-5pt\mathord!B\otimes \mathord!A}{\mathord!(B\otimes A)\kern-5pt}%
  {\tilde\varphi}{\sigma}{\mathord!\sigma}{\tilde\varphi}
 \di{\hbox{$\circlearrowright$}}(25,20){0}
\end{diagramme}\kern-9pt\hss}}$
\cr}
\endgroup

	\vskip2ex
        \noindent
 where $f$ and $g$ are arbitrary, and where $\alpha$ and $\sigma$ are
 structural isomorphisms.
The defining diagrams of comonad and symmetric comonoid
 comprise twenty rewriting rules, except
 the following two congruent diagrams:

	\vskip2ex
\begingroup\offinterlineskip
\halign to\textwidth{\kern120pt\hfil #\hfil \tabskip1000pt plus100pt minus1000pt
 &\kern30pt\hfil #\hfil \kern120pt\tabskip0pt
\cr
$\vcenter{\footnotesize
\hbox to0pt{\hss\begin{diagramme}
\pentagondiagabscoord{70}{60}%
 {\bang A}{\bang A\otimes \bang A}{\bang A\otimes \bang A}%
 {\kern-40pt(\bang A\otimes \bang A)\otimes \bang A}%
 {\bang A\otimes (\bang A\otimes \bang A)\kern-40pt}%
 dd{d\otimes 1}{1\otimes d}{\alpha}
 \di{\hbox{$\circlearrowright$}}(35,25){0}
\end{diagramme}\hss}}$%
&$\vcenter{\footnotesize
\hbox to0pt{\hss\begin{diagramme}
\trianglediagabscoord{50}{30}%
   {\mathord!A}{\mathord!A\otimes \mathord!A}{\mathord!A\otimes \mathord!A}%
   dd{\sigma}
 \di{\hbox{$\circlearrowright$}}(25,18){0}
\end{diagramme}\kern-16pt\hss}}$
\cr
}
\endgroup

	\vskip2ex
        \noindent
For the complete list, we refer the reader to \cite{hase2}.
Our policy is to reformulate as many diagrams as possible to
 one-way rewriting.
Only those diagrams which lead to an unreasonable system if
 we enforce rewriting are left as congruence.

\section{Linear normal functor model}\label{ypr76}

Multisets are frequently used below.
In this paper, only finite multiset appears.
The multiset union is denoted by $\alpha+
\beta$.

The confluence of the calculus on the
 linear category is proved via
 a specific model $\mathscr{M}$ given in \cite{hase1}.
It is an adaptation of the quantitative model of Girard \cite{gira1}.
This model is based on
 normal functors, which are special cases
 of Joyal's analytic functors \cite{joya}.
The category of normal functors can be regarded
 as a coKleisli category of the exponential comonad given below.
This comonad satisfies all the conditions of the linear category.
In the following definition, ${\bf Set}^A$ denotes
 the category of presheaves where $A$ is a set regarded as a
 discrete category.

	\vskip2ex

\begin{definition}\rm
The objects of $\mathscr{M}$ are sets $A$.
A morphism from $A$ to $B$ is a linear normal functor
 from ${\bf Set}^A$ to ${\bf Set}^B$.
Here linear normal functors are defined as the functors preserving
 equalizers, pullbacks of possibly
 infinite legs, and all colimits \cite{hase1}.
\end{definition}

	\vskip2ex
        \noindent
We comment that, if we replace colimits with
 filtered colimits in the definition,
 we obtain the definition of normal functors \cite{gira1}.
In place of universality, it is
 convenient to use a concrete presentation, which
 we explain shortly.

An object $x\in {\bf Set}^A$ has an analogy to
 a vector in the ordinary
 linear algebra.
To each $a\in A$, a component $x[a]$ corresponds.
Since $A$ is infinite in general, it may be a vector
 of infinite dimension.
In linear algebra, each component is a member of the underlying
 field.
Here, $x[a]$ is an object of ${\bf Set}$, i.e., a set.
If we identify a set with its cardinality,
 $x$ is a vector having cardinal numbers as its components.
In particular, if they are finite sets, the components
 are non-negative integers.
We can perform addition and multiplication of components, but no
 subtraction as they are cardinal numbers.

Linear normal functors have the following characterization.
Suppose that $f:{\bf Set}^A\rightarrow{\bf Set}^B$ is a linear
 normal functor.
If we write $y=f(x)$, there is $M\in {\bf Set}^{A\times B}$
 such that $y[b]=\sum_aM[a;b]x[a]$, where $M[a;b]$ denotes
 the set that is the value of $M$ at $(a,b)\in A\times B$.
The sum is a possibly infinite disjoint sum of sets over $a\in A$,
 and the concatenation signifies the cartesian product.
The equality is not exactly correct, as it usually holds
 only up to isomorphism.
We abuse the equality to emphasize the analogy with
 linear algebra.
In ${\bf Set}$, two sets having the same cardinality
 are isomorphic.
So, if we identify them, the sum and the product obey
 the ordinary cardinal arithmetic.
In particular, if only finite sets are involved,
 they are elementary arithmetic on integers.
In this way, a linear normal functor has an analogy with the
 linear map given by the multiplication of a matrix $M$.
We write $y=Mx$.
Conversely, if a functor is naturally equivalent to
 the functor of the form $y=Mx$ then
 it is a linear normal functor.
All linear normal functors from ${\bf Set}^A$ to
 ${\bf Set}^B$  are determined by
 $M\in{\bf Set}^{A\times B}$ up to natural equivalences.

The coefficient matrix $M\in{\bf Set}^{A\times B}$ is
 also regarded as a span (a profunctor).
If we abuse the same symbol $M$, a diagram
$$\vcenter{\spandiagabs{50}{30}MAB{}{}}$$
 in ${\bf Set}$ corresponds.
Regarding this diagram as a function into $A\times B$,
 the coefficient $M[a;b]$ is the inverse image of
 $(a,b)$ under the function.
The composition of linear normal functors correspond to
 the composition of profunctors, i.e., the pullback.
The coefficient matrices or spans
 give only equivalent presentations of linear normal
 functors.
They do not form a category in an exact sense.
For example, if we consider the category of spans,
 $(L\times_A M)\times_BN$ and $L\times_A
(M\times_BN)$ are isomorphic, but not exactly equal in general.
This type of subtle difference, however, is
 not problematic in this paper.
So let us identify them naively.

Now we install a structure of the linear category on
$\mathscr{M}$.
First of all, the duality is simply ignored.
Namely, $A^*=A$.
Tensor and cotensor are both the cartesian product.
Namely, both of $A\otimes B$ and $A\linpar B$ are given
 as $A\times B$.
The morphism maps are given as the tensor of matrices.
Namely, if $f$ and $g$ correspond to the matrices $M$ and $N$
 respectively, $f\otimes g$ and $f\linpar g$ both
 correspond to $M\otimes N$ given
 as $(M\otimes N)[(a,a');(b,b')]=M[a;b]N[a';b']$.
The tensor unit and cotensor unit are both a singleton ${\bf 1}$.
Its element is denoted by $*$.

We give the endofunctor $\mathord!:\mathscr{M}\rightarrow
\mathscr{M}$.
The object map is $A\mapsto \mathop{\rm exp}A$.
Here, $\mathop{\rm exp}A$
 denotes the set of all finite multisets of members of $A$.
This is an analytic functor in the sense of Joyal
 and this notation is introduced there  \cite{joya}.
It is instructive to identify the multisets in $\mathop{\rm exp}A$
 with monomials in the variables that have a one-to-one correspondence
 to the member of $A$.
With each $a\in A$, we associate a variable $x_a$.
Then, for example, the multiset $\alpha=\{a,a,b\}$ is identified with
 the monomial $x^\alpha=x_a\!{}^2x_b$.
The morphism $\mathord!f:{\bf Set}^{\mathop{\rm exp}A}
\rightarrow {\bf Set}^{\mathop{\rm exp}B}$ is
 naturally defined under this identification.
Suppose that $f$ is given by $M\in {\bf Set}^{A\times B}$.
Substituting for each variable in the monomial $y^\beta$
 by $y_b=\sum_aM[a,b]x_a$,
 we obtain a linear combination of monomials $x^\alpha$.
Namely, 
 we can write $y^\beta=\sum_\alpha\tilde M[\alpha,\beta]x^\alpha$.
This matrix $\tilde M\in {\bf Set}^{\mathop{\rm exp}A\times\mathop{\rm exp}B}$
 gives the definition of $\mathord!f$.
For the reader's convenience, we give a sketch in a simple
 case.
Suppose that $A=\{1,2\}=B$ and $M=\left(\matrix{%
 a & b \cr
 c & d \cr}\right)$.
Namely, $y_1=ax_1+bx_2$ and $y_2=cx_1+dx_2$.
Let us enumerate the members of $\mathop{\rm exp}A$ in the
 order of $1,\,x_1,\,x_2,\,x_1\!{}^2,
\,x_1x_2,\,x_2\!{}^2,\ldots$ and
 those of $\mathop{\rm exp}B$ similarly.
Then the matrix $\tilde M$ is given as
$$\left(\matrix{%
 1\cr
 &a&b\cr &c&d\cr
 &&&a^2&2ab&b^2\cr
 &&&ac&ad+bc&bd\cr
 &&&c^2&2cd&d^2\cr
 &&&&&&\ddots&\cr
}\right).$$
For example, the fourth row comes from
 $y_1\!{}^2=a^2x_1\!{}^2+2ab\,x_1x_2+b^2x_2\!{}^2$.
Since $\mathord!f$ is defined by substitution,
 $\mathord!(f\scirc g)=\mathord!f\scirc\mathord!g$ is
 immediate.
Hence $\mathord!$ turns out to be a functor.

We should observe asymmetry
 of the matrix.
When we view the $3\times 3$ block of degree two monomials
 in the example of $\tilde M$ above,
 the $(1,2)$-component $2ab$ has the coefficient $2$, whereas
 the $(2,1)$-component is $ac$, the coefficient of which is $1$.
In terms of spans, this phenomenon is explained as follows.
Suppose that $f$ corresponds to the span $M$ over $A$ and $B$.
Then $\mathord!f$ does {\it not} correspond to
\bgroup\small $\vcenter{\spandiagabs{40}{24}%
{\mathop{\rm exp}M}{\mathop{\rm exp}A}{\mathop{\rm exp}B}{}{}}$\egroup.
In fact, the analytic functor $\mathop{\rm exp}A$ does not preserve
 pullbacks.
Hence, if we adopted this symmetric span as the definition,
 $\mathord!$ would not preserve
 composition.
Therefore, the definition must intrinsically
 be asymmetric on two legs.
Joyal's analytic functors are introduced in order to
 give a mathematical foundation to the generating
 functions in the enumerative combinatorics.
This asymmetry is the origin that complicates
 and, at the same time, enriches the
 theory of enumeration.
If the symmetry held, the world of enumerative combinatorics
 would be much more languid.

Each function $f:A\rightarrow B$ induces the linear normal
 functor $f^*:{\bf Set}^B\rightarrow {\bf Set}^A$ defined
 by $f^*(y)[a]=x[f(a)]$.
In terms of spans, it corresponds to
$$\vcenter{\spandiagabs{45}{27}%
 ABA{f}{}}$$
In terms of matrices, it corresponds to the matrix $M[b;a]$ that
 equals $1$ whenever $b=f(a)$; otherwise equals $0$.

Using $f^*$, we can give the structures of the linear category
 to $\mathord!$.
The morphism $\delta_A:{\bf Set}^{\mathop{\rm exp}A}
\rightarrow {\bf Set}^{\mathop{\rm exp}\mathop{\rm exp}A}$
 is induced from the function
 $\mathop{\rm exp}\mathop{\rm exp}A
 \rightarrow\mathop{\rm exp}A$
 carrying $\{\alpha_1,\alpha_2,\ldots,\alpha_p\}$ to
 the multiset union $\alpha_1+\alpha_2+\cdots+\alpha_p$.
The morphism $\varepsilon_A:{\bf Set}^{\mathop{\rm exp}A}
 \rightarrow {\bf Set}^A$ is induced from
 the function $A\rightarrow\mathop{\rm exp}A$
 carrying $a$ to the singleton $\{a\}$.
The morphism $d_A:{\bf Set}^{\mathop{\rm exp}A}\rightarrow
{\bf Set}^{\mathop{\rm exp}A\times \mathop{\rm exp}A}$ is induced
 from the function $\mathop{\rm exp}A\times \mathop{\rm exp}A
 \rightarrow \mathop{\rm exp}A$ carrying $(\alpha,\alpha')$ to
 the multiset union $\alpha+\alpha'$.
The morphism $e_A:{\bf Set}^{\mathop{\rm exp}A}\rightarrow
{\bf Set}^{\bf 1}$ is induced from the function ${\bf 1}
 \rightarrow \mathop{\rm exp}A$ carrying $*$ to
 the empty multiset $\emptyset$.
We can directly check that these are natural transformations.
The morphism $\tilde\varphi:{\bf Set}^{\mathop{\rm exp}A\times
\mathop{\rm exp}A}\rightarrow {\bf Set}^{\mathop{\rm exp}(A\times A)}$
 is induced from the function $\mathop{\rm exp}(A\times A)
 \rightarrow \mathop{\rm exp}A\times \mathop{\rm exp}A$ carrying
 $\{(a_1,a'_1),(a_2,a'_2),\ldots,(a_p,a'_p)\}$ to the pair
 $(\{a_1,a_2,\ldots,a_p\},\,\{a'_1,a'_2,\ldots,a'_p\})$.
The morphism $\varphi_0:{\bf Set}^{{\bf 1}}\rightarrow
 {\bf Set}^{\mathop{\rm exp}{\bf 1}}$ is induced from
 the unique function $\mathop{\rm exp}{\bf 1}\rightarrow
 {\bf 1}$.
The naturality of $\tilde\varphi$ is directly checked.
Since all morphisms are of the shape $f^*$ for functions
 $f$, the coherence conditions are easily derived from
 the corresponding equalities between the associated
 functions.

Our calculus is based on the free linear category.
The objects and morphisms of the calculus are thus interpreted
 in the linear category $\mathscr{M}$, once the intepretations
 of the atomic objects $X$ are provided.
In particular, every syntactic morphism $f:A\rightarrow B$
 is interpreted by $[\![f]\!]:[\![A]\!]\rightarrow [\![B]\!]$
 in $\mathscr{M}$.
When we are involved in the matrix representation,
 we write $M_f\in {\bf Set}^{[\![A]\!]\times [\![B]\!]}$.
The sets interpreting the atomic objects are arbitrary.
To our aim, however, we assume infinite sets.

In order to distinguish the objects of the model
 from the objects of the syntactic rewriting system,
 we call the latter objects {\it types}.
Namely $[\![A]\!]$ is the object of $\mathscr{M}$ interpreting
 the type $A$.

The elements of the interpretation of an atomic type are
 called {\it atoms}.
The interpretations of types are either those associated with
 atomic types, ${\bf 1}$, $[\![A]\!]
\times [\![B]\!]$, or
 $\mathop{\rm exp}[\![A]\!]$.
Hence, an element of the interpretation of a
 type is either an atom,
 $*$, a pair $(a,b)$, or a finite multiset $\{a_1,a_2,\ldots,a_n\}$.
Although the category $\mathscr{M}$ does not distinguish between
 $A$ and $A^*$, if we have the information of types.
 we can recover positive and negative occurrences.
To this end, we introduce the notion of signed elements
 using overbars.
Provided that $\alpha$ is a signed element of $[\![A]\!]$, we
 write $\overline \alpha$ if we regard it as a signed element
 of $[\![A^*]\!]$.
Note that $\alpha$ and $\overline \alpha$ are equal as the members of
 $[\![A]\!]=[\![A^*]\!]$.
We say that an occurrence of $\alpha$ is positive (or negative) if
 it occurs under an even (or odd) number of nested overbars.
For example, $\alpha$ is positive and $\beta$ is negative
 in $\overline{(\overline \alpha,\beta)}$.
In the coming argument, signed atoms and positive multiset occurrences
 play important roles.
When we are concerned with $[\![f]\!]:[\![B]\!]\rightarrow
[\![C]\!]$ or $M_f\in {\bf Set}^{[\![B]\!]\times [\![C]\!]}$,
 we give signs as $(\overline \beta;\gamma)$ for $\beta\in [\![B]\!]$
 and $\gamma\in [\![C]\!]$.
Accordingly,
 $\alpha$ occurs positively in $(\beta;\gamma)$
 if it occurs positively
in $\gamma$ or negatively in $\beta$.

\section{Graphics}

We give a graphical presentation of normal forms
 of our calculus.
It is an extension of Blute-Cockett-Seely-Trimble graphs
 for the $*$-autonomous category \cite{bcst}.
The graphical presentation absorbes equivalences
 of the calculus, thus facilitates arguments modulo
 congruence.
We can dispense with diagram chasing with regard to
 the coherence equivalences.
We emphasize that contstructing a graphical calculus
 is not our purpose.
Indeed, we associate graphs with normal forms only.

Graphs are obtained by linking parts by wires (rough analogy with
electric circuits).
The directions up/down and left/right matter.
Each wire has a type.
It is often omitted as it is restored from the shape of the graph.
The tensor parts with the wires attaching to them are
 the following two:
$$
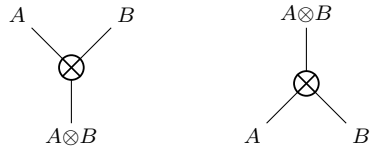
\begin{tikzpicture}[xscale=0.0352778, yscale=0.0352778, thin, inner sep=0]
  \def\p#1#2{%
    \ifcase #1
      \or \ifx#2x  0       \else  0 \fi
    \fi}
  \def\x#1{\p#1x}
  \def\y#1{\p#1y}
  \path[use as bounding box] (\x1-15,\y1-29) rectangle (\x1+15,\y1+23);
  \draw (\x1-15,\y1+15) -- (\x1,\y1);
  \draw (\x1+15,\y1+15) -- (\x1,\y1);
  \draw (\x1,\y1-21) -- (\x1,\y1);
  \putTensor(\x1,\y1)
  \node[above left,inner sep=2] at (\x1-15,\y1+15) {$\scriptstyle A$};
  \node[above right,inner sep=2] at (\x1+15,\y1+15) {$\scriptstyle B$};
  \node[below,inner sep=2] at (\x1,\y1-21) {$\scriptstyle A\otimes B$};
\end{tikzpicture}
	\hskip5em
\begin{tikzpicture}[xscale=0.0352778, yscale=0.0352778, thin, inner sep=0]
  \def\p#1#2{%
    \ifcase #1
      \or \ifx#2x  0       \else  0 \fi
    \fi}
  \def\x#1{\p#1x}
  \def\y#1{\p#1y}
  \path[use as bounding box] (\x1-15,\y1+29) rectangle (\x1+15,\y1-23);
  \draw (\x1-15,\y1-15) -- (\x1,\y1);
  \draw (\x1+15,\y1-15) -- (\x1,\y1);
  \draw (\x1,\y1+21) -- (\x1,\y1);
  \putTensor(\x1,\y1)
  \node[below left,inner sep=2] at (\x1-15,\y1-15) {$\scriptstyle A$};
  \node[below right,inner sep=2] at (\x1+15,\y1-15) {$\scriptstyle B$};
  \node[above,inner sep=2] at (\x1,\y1+21) {$\scriptstyle A\otimes B$};
\end{tikzpicture}
$$
The cotensor parts are the following two:
$$
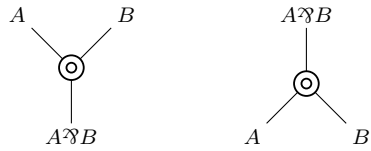
\begin{tikzpicture}[xscale=0.0352778, yscale=0.0352778, thin, inner sep=0]
  \def\p#1#2{%
    \ifcase #1
      \or \ifx#2x  0       \else  0 \fi
    \fi}
  \def\x#1{\p#1x}
  \def\y#1{\p#1y}
  \path[use as bounding box] (\x1-15,\y1-29) rectangle (\x1+15,\y1+23);
  \draw (\x1-15,\y1+15) -- (\x1,\y1);
  \draw (\x1+15,\y1+15) -- (\x1,\y1);
  \draw (\x1,\y1-21) -- (\x1,\y1);
  \putPar(\x1,\y1)
  \node[above left,inner sep=2] at (\x1-15,\y1+15) {$\scriptstyle A$};
  \node[above right,inner sep=2] at (\x1+15,\y1+15) {$\scriptstyle B$};
  \node[below,inner sep=2] at (\x1,\y1-21) {$\scriptstyle A\slinpar B$};
\end{tikzpicture}
	\hskip5em
\begin{tikzpicture}[xscale=0.0352778, yscale=0.0352778, thin, inner sep=0]
  \def\p#1#2{%
    \ifcase #1
      \or \ifx#2x  0       \else  0 \fi
    \fi}
  \def\x#1{\p#1x}
  \def\y#1{\p#1y}
  \path[use as bounding box] (\x1-15,\y1+29) rectangle (\x1+15,\y1-23);
  \draw (\x1-15,\y1-15) -- (\x1,\y1);
  \draw (\x1+15,\y1-15) -- (\x1,\y1);
  \draw (\x1,\y1+21) -- (\x1,\y1);
  \putPar(\x1,\y1)
  \node[below left,inner sep=2] at (\x1-15,\y1-15) {$\scriptstyle A$};
  \node[below right,inner sep=2] at (\x1+15,\y1-15) {$\scriptstyle B$};
  \node[above,inner sep=2] at (\x1,\y1+21) {$\scriptstyle A\slinpar B$};
\end{tikzpicture}
$$
We use the double circle in place of $\linpar$.
The unit parts are
$$
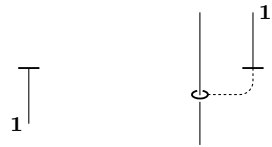
\begin{tikzpicture}[xscale=0.0352778, yscale=0.0352778, thin, inner sep=0]
  \def\p#1#2{%
    \ifcase #1
      \or \ifx#2x  0       \else  0 \fi
    \fi}
  \def\x#1{\p#1x}
  \def\y#1{\p#1y}
  \path[use as bounding box] (\x1-3,\y1-29) rectangle (\x1+3,\y1+29);
  \draw (\x1,\y1-21) -- (\x1,\y1);
  \putPositiveTerminal(\x1,\y1)
  \node[left,inner sep=2] at (\x1,\y1-21) {$\scriptstyle {\bf 1}$};
\end{tikzpicture}
	\hskip5em
\begin{tikzpicture}[xscale=0.0352778, yscale=0.0352778, thin, inner sep=0]
  \def\p#1#2{%
    \ifcase #1
      \or \ifx#2x  0       \else  0 \fi
    \fi}
  \def\x#1{\p#1x}
  \def\y#1{\p#1y}
  \path[use as bounding box] (\x1-23,\y1-29) rectangle (\x1+3,\y1+29);
  \draw (\x1,\y1+21) -- (\x1,\y1);
  \draw (\x1-20,\y1+21) -- (\x1-20,\y1-29);
  \putPositiveTerminal(\x1,\y1)
  \putRubberBand(\x1-20,\y1-10)
  \draw[dash pattern=on 1 off 1] (\x1,\y1)[rounded corners=5] -- (\x1,\y1-10)
    -- (\x1-17,\y1-10);
  \node[right,inner sep=2] at (\x1,\y1+21) {$\scriptstyle {\bf 1}$};
\end{tikzpicture}
$$
The counit parts are
$$
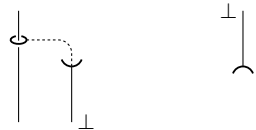
\begin{tikzpicture}[xscale=0.0352778, yscale=0.0352778, thin, inner sep=0]
  \def\p#1#2{%
    \ifcase #1
      \or \ifx#2x  0       \else  0 \fi
    \fi}
  \def\x#1{\p#1x}
  \def\y#1{\p#1y}
  \path[use as bounding box] (\x1-23,\y1-29) rectangle (\x1+3,\y1+29);
  \draw (\x1,\y1-21) -- (\x1,\y1);
  \draw (\x1-20,\y1-21) -- (\x1-20,\y1+21);
  \putUpperNegativeTerminal(\x1,\y1)
  \putRubberBand(\x1-20,\y1+10)
  \draw[dash pattern=on 1 off 1] (\x1,\y1)[rounded corners=5] -- (\x1,\y1+10)
    -- (\x1-17,\y1+10);
  \node[right,inner sep=2] at (\x1,\y1-21) {$\scriptstyle \bot$};
\end{tikzpicture}
	\hskip5em
\begin{tikzpicture}[xscale=0.0352778, yscale=0.0352778, thin, inner sep=0]
  \def\p#1#2{%
    \ifcase #1
      \or \ifx#2x  0       \else  0 \fi
    \fi}
  \def\x#1{\p#1x}
  \def\y#1{\p#1y}
  \path[use as bounding box] (\x1-3,\y1-29) rectangle (\x1+3,\y1+29);
  \draw (\x1,\y1+21) -- (\x1,\y1);
  \putLowerNegativeTerminal(\x1,\y1)
  \node[left,inner sep=2] at (\x1,\y1+21) {$\scriptstyle \bot$};
\end{tikzpicture}
$$
The duality parts are
$$
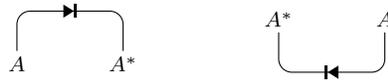
\begin{tikzpicture}[xscale=0.0352778, yscale=0.0352778, thin, inner sep=0]
  \def\p#1#2{%
    \ifcase #1
      \or \ifx#2x  0       \else  0 \fi
    \fi}
  \def\x#1{\p#1x}
  \def\y#1{\p#1y}
  \path[use as bounding box] (\x1-20,\y1-23) rectangle (\x1+20,\y1);
  \draw (\x1+20,\y1-15)[rounded corners=5] -- (\x1+20,\y1)
    -- (\x1-20,\y1) -- (\x1-20,\y1-15);
  \putRightDiode(\x1,\y1)
  \node[below,inner sep=2] at (\x1-20,\y1-15) {$\scriptstyle A$};
  \node[below,inner sep=2] at (\x1+20,\y1-15) {$\scriptstyle A^*$};
\end{tikzpicture}
	\hskip5em
\begin{tikzpicture}[xscale=0.0352778, yscale=0.0352778, thin, inner sep=0]
  \def\p#1#2{%
    \ifcase #1
      \or \ifx#2x  0       \else  0 \fi
    \fi}
  \def\x#1{\p#1x}
  \def\y#1{\p#1y}
  \path[use as bounding box] (\x1-20,\y1) rectangle (\x1+20,\y1+23);
  \draw (\x1+20,\y1+15)[rounded corners=5] -- (\x1+20,\y1)
    -- (\x1-20,\y1) -- (\x1-20,\y1+15);
  \putLeftDiode(\x1,\y1)
  \node[above,inner sep=2] at (\x1-20,\y1+15) {$\scriptstyle A^*$};
  \node[above,inner sep=2] at (\x1+20,\y1+15) {$\scriptstyle A$};
\end{tikzpicture}
$$
In all of these, the left graph
 is called the introduction rule and the right
 the elimination rule.
Wires can cross.
However, the choice of which wire lays over the other does
 not matter as we consider symmetric monoidal categories:
$$
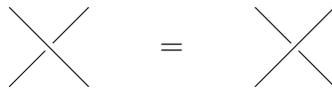
\vcenter{\hbox{%
\begin{tikzpicture}[xscale=0.0352778, yscale=0.0352778, thin, inner sep=0]
  \def\p#1#2{%
    \ifcase #1
      \or \ifx#2x  0       \else  0 \fi
    \fi}
  \def\x#1{\p#1x}
  \def\y#1{\p#1y}
  \path[use as bounding box] (\x1-15,\y1-15) rectangle (\x1+15,\y1+15);
  \draw (\x1+15,\y1+15) -- (\x1-15,\y1-15);
  \draw[white,line width=4.0] (\x1-15,\y1+15) -- (\x1+15,\y1-15);
  \draw (\x1-15,\y1+15) -- (\x1+15,\y1-15);
\end{tikzpicture}%
}}
\qquad=\qquad
\vcenter{\hbox{%
\begin{tikzpicture}[xscale=0.0352778, yscale=0.0352778, thin, inner sep=0]
  \def\p#1#2{%
    \ifcase #1
      \or \ifx#2x  0       \else  0 \fi
    \fi}
  \def\x#1{\p#1x}
  \def\y#1{\p#1y}
  \path[use as bounding box] (\x1-15,\y1-15) rectangle (\x1+15,\y1+15);
  \draw (\x1-15,\y1+15) -- (\x1+15,\y1-15);
  \draw[white,line width=4.0] (\x1+15,\y1+15) -- (\x1-15,\y1-15);
  \draw (\x1+15,\y1+15) -- (\x1-15,\y1-15);
\end{tikzpicture}%
}}
$$
It may be intuitive to regard the dotted wire as
 a fishing line.
It is difficult to see, but certainly exists to connect the
 part to another wire.
A rubberband gives a good analogy for the ring connecting
 the fishing line to a wire.
It is elastic, can move around, and can rotate around a wire:
$$
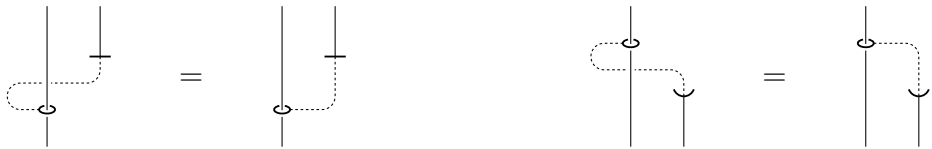
\vcenter{\hbox{%
\begin{tikzpicture}[xscale=0.0352778, yscale=0.0352778, thin, inner sep=0]
  \def\p#1#2{%
    \ifcase #1
      \or \ifx#2x  0       \else  0 \fi
    \fi}
  \def\x#1{\p#1x}
  \def\y#1{\p#1y}
  \path[use as bounding box] (\x1-23,\y1-34) rectangle (\x1+3,\y1+19);
  \draw (\x1,\y1+19) -- (\x1,\y1);
  \draw[dash pattern=on 1 off 1] (\x1,\y1)[rounded corners=5] -- (\x1,\y1-10)
    -- (\x1-35,\y1-10) -- (\x1-35,\y1-20) -- (\x1-23,\y1-20);
  \draw[white,line width=3.0] (\x1-20,\y1+19) -- (\x1-20,\y1-34);
  \draw (\x1-20,\y1+19) -- (\x1-20,\y1-34);
  \putPositiveTerminal(\x1,\y1)
  \putRubberBand(\x1-20,\y1-20)
\end{tikzpicture}%
}}
\qquad=\qquad
\vcenter{\hbox{%
\begin{tikzpicture}[xscale=0.0352778, yscale=0.0352778, thin, inner sep=0]
  \def\p#1#2{%
    \ifcase #1
      \or \ifx#2x  0       \else  0 \fi
    \fi}
  \def\x#1{\p#1x}
  \def\y#1{\p#1y}
  \path[use as bounding box] (\x1-23,\y1-34) rectangle (\x1+3,\y1+19);
  \draw (\x1,\y1+19) -- (\x1,\y1);
  \draw (\x1-20,\y1+19) -- (\x1-20,\y1-34);
  \putPositiveTerminal(\x1,\y1)
  \putRubberBand(\x1-20,\y1-20)
  \draw[dash pattern=on 1 off 1] (\x1,\y1)[rounded corners=5] -- (\x1,\y1-20)
    -- (\x1-17,\y1-20);
\end{tikzpicture}%
}}
	\hskip9em
\vcenter{\hbox{%
\begin{tikzpicture}[xscale=0.0352778, yscale=0.0352778, thin, inner sep=0]
  \def\p#1#2{%
    \ifcase #1
      \or \ifx#2x  0       \else  0 \fi
    \fi}
  \def\x#1{\p#1x}
  \def\y#1{\p#1y}
  \path[use as bounding box] (\x1-23,\y1+34) rectangle (\x1+3,\y1-19);
  \draw (\x1,\y1-19) -- (\x1,\y1);
  \draw[dash pattern=on 1 off 1] (\x1,\y1)[rounded corners=5] -- (\x1,\y1+10)
    -- (\x1-35,\y1+10) -- (\x1-35,\y1+20) -- (\x1-23,\y1+20);
  \draw[white, line width=3.0] (\x1-20,\y1-19) -- (\x1-20,\y1+34);
  \draw (\x1-20,\y1-19) -- (\x1-20,\y1+34);
  \putUpperNegativeTerminal(\x1,\y1)
  \putRubberBand(\x1-20,\y1+20)
\end{tikzpicture}%
}}
\qquad=\qquad
\vcenter{\hbox{%
\begin{tikzpicture}[xscale=0.0352778, yscale=0.0352778, thin, inner sep=0]
  \def\p#1#2{%
    \ifcase #1
      \or \ifx#2x  0       \else  0 \fi
    \fi}
  \def\x#1{\p#1x}
  \def\y#1{\p#1y}
  \path[use as bounding box] (\x1-23,\y1+34) rectangle (\x1+3,\y1-19);
  \draw (\x1,\y1-19) -- (\x1,\y1);
  \draw (\x1-20,\y1-19) -- (\x1-20,\y1+34);
  \putUpperNegativeTerminal(\x1,\y1)
  \putRubberBand(\x1-20,\y1+20)
  \draw[dash pattern=on 1 off 1] (\x1,\y1)[rounded corners=5] -- (\x1,\y1+20)
    -- (\x1-17,\y1+20);
\end{tikzpicture}%
}}
$$
This completes the definition of the graphs of the $*$-autonomous category.
The modification rules of graphs corresponding to
 the coherence conditions of units are given in \cite{bcst}.
Two graphs transferring to one another
 are regarded to be equal.
Not all graphs are legitimate.
We consider only the graphs subject to the well-known
 switching condition~\cite{dare}.

We extend the graph to the classical linear category.
We add the following $\delta$-part, $\varepsilon$-part,
 $d$-part, and $e$-part:
$$
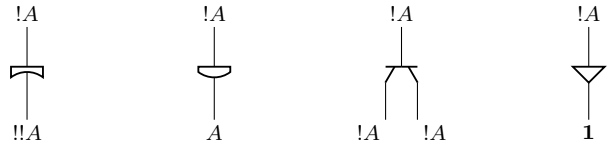
\begin{tikzpicture}[xscale=0.0352778, yscale=0.0352778, thin, inner sep=0]
  \def\p#1#2{%
    \ifcase #1
      \or \ifx#2x  0       \else  0 \fi
    \fi}
  \def\x#1{\p#1x}
  \def\y#1{\p#1y}
  \path[use as bounding box] (\x1-6,\y1+22) rectangle (\x1+6,\y1-27);
  \draw (\x1,\y1+15) -- (\x1,\y1-20);
  \putConcave(\x1,\y1)
  \node[above,inner sep=2] at (\x1,\y1+15) {$\scriptstyle\mathord!A$};
  \node[below,inner sep=2] at (\x1,\y1-20) {$\scriptstyle\mathord!\mathord!A$};
\end{tikzpicture}%
	\hskip5em
\begin{tikzpicture}[xscale=0.0352778, yscale=0.0352778, thin, inner sep=0]
  \def\p#1#2{%
    \ifcase #1
      \or \ifx#2x  0       \else  0 \fi
    \fi}
  \def\x#1{\p#1x}
  \def\y#1{\p#1y}
  \path[use as bounding box] (\x1-6,\y1+22) rectangle (\x1+6,\y1-27);
  \draw (\x1,\y1+15) -- (\x1,\y1-20);
  \putConvex(\x1,\y1)
  \node[above,inner sep=2] at (\x1,\y1+15) {$\scriptstyle\mathord!A$};
  \node[below,inner sep=2] at (\x1,\y1-20) {$\scriptstyle A$};
\end{tikzpicture}%
	\hskip5em
\begin{tikzpicture}[xscale=0.0352778, yscale=0.0352778, thin, inner sep=0]
  \def\p#1#2{%
    \ifcase #1
      \or \ifx#2x  0       \else  0 \fi
    \fi}
  \def\x#1{\p#1x}
  \def\y#1{\p#1y}
  \path[use as bounding box] (\x1-6,\y1+22) rectangle (\x1+6,\y1-27);
  \draw (\x1,\y1+15) -- (\x1,\y1);
  \draw (\x1-6,\y1-6) -- (\x1-6,\y1-20);
  \draw (\x1+6,\y1-6) -- (\x1+6,\y1-20);
  \putDuplicator(\x1,\y1)
  \node[above,inner sep=2] at (\x1,\y1+15) {$\scriptstyle\mathord!A$};
  \node[below left,inner sep=2] at (\x1-6,\y1-20) {$\scriptstyle\mathord!A$};
  \node[below right,inner sep=2] at (\x1+6,\y1-20) {$\scriptstyle\mathord!A$};
\end{tikzpicture}%
	\hskip5em
\begin{tikzpicture}[xscale=0.0352778, yscale=0.0352778, thin, inner sep=0]
  \def\p#1#2{%
    \ifcase #1
      \or \ifx#2x  0       \else  0 \fi
    \fi}
  \def\x#1{\p#1x}
  \def\y#1{\p#1y}
  \path[use as bounding box] (\x1-6,\y1+22) rectangle (\x1+6,\y1-27);
  \draw (\x1,\y1+15) -- (\x1,\y1-20);
  \putEliminator(\x1,\y1)
  \node[above,inner sep=2] at (\x1,\y1+15) {$\scriptstyle\mathord!A$};
  \node[below,inner sep=2] at (\x1,\y1-20) {$\scriptstyle {\bf 1}$};
\end{tikzpicture}%
$$
The optic lense symbols for $\delta$ and $\varepsilon$ are borrowed
 from \cite{aspe,lamp}.
The $d$-part is also called a duplicator, and the $e$-part
 an eliminator.
By one of the coherence conditions, the duplicator satisfy
$$
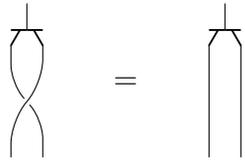
\vcenter{\hbox{%
\begin{tikzpicture}[xscale=0.0352778, yscale=0.0352778, thin, inner sep=0]
  \def\p#1#2{%
    \ifcase #1
      \or \ifx#2x  0       \else  0 \fi
    \fi}
  \def\x#1{\p#1x}
  \def\y#1{\p#1y}
  \path[use as bounding box] (\x1-6,\y1+10) rectangle (\x1+6,\y1-48);
  \draw (\x1,\y1+10) -- (\x1,\y1);
  \draw (\x1-6,\y1-6)[rounded corners=7] -- (\x1-6,\y1-20)
    -- (\x1+6,\y1-34) -- (\x1+6,\y1-48);
  \draw[white,line width=3.0] (\x1+6,\y1-6)[rounded corners=9] -- (\x1+6,\y1-20)
    -- (\x1-6,\y1-34) -- (\x1-6,\y1-48);
  \draw (\x1+6,\y1-6)[rounded corners=9] -- (\x1+6,\y1-20)
    -- (\x1-6,\y1-34) -- (\x1-6,\y1-48);
  \putDuplicator(\x1,\y1)
\end{tikzpicture}%
}}
\qquad=\qquad
\vcenter{\hbox{%
\begin{tikzpicture}[xscale=0.0352778, yscale=0.0352778, thin, inner sep=0]
  \def\p#1#2{%
    \ifcase #1
      \or \ifx#2x  0       \else  0 \fi
    \fi}
  \def\x#1{\p#1x}
  \def\y#1{\p#1y}
  \path[use as bounding box] (\x1-6,\y1+10) rectangle (\x1+6,\y1-48);
  \draw (\x1,\y1+10) -- (\x1,\y1);
  \draw (\x1-6,\y1-6) -- (\x1-6,\y1-48);
  \draw (\x1+6,\y1-6) -- (\x1+6,\y1-48);
  \putDuplicator(\x1,\y1)
\end{tikzpicture}%
}}
$$
Namely, a duplicator can rotate around the axis.
Moreover, the following equality holds:
$$
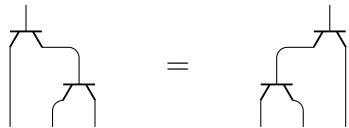
\vcenter{\hbox{%
\begin{tikzpicture}[xscale=0.0352778, yscale=0.0352778, thin, inner sep=0]
  \def\p#1#2{%
    \ifcase #1
      \or \ifx#2x  0       \else  0 \fi
      \or \ifx#2x  \x1+20  \else  \y1-20 \fi
    \fi}
  \def\x#1{\p#1x}
  \def\y#1{\p#1y}
  \path[use as bounding box] (\x1-6,\y1+10) rectangle (\x1+26,\y1-36);
  \draw (\x1,\y1+10) -- (\x1,\y1);
  \draw (\x1+6,\y1-6)[rounded corners=4] -- (\x2,\y1-6) -- (\x2,\y2);
  \draw (\x1-6,\y1-6) -- (\x1-6,\y2-16);
  \draw (\x2-6,\y2-6)[rounded corners=4] -- (\x2-10,\y2-6) -- (\x2-10,\y2-16);
  \draw (\x2+6,\y2-6) -- (\x2+6,\y2-16);
  \putDuplicator(\x1,\y1)
  \putDuplicator(\x2,\y2)
\end{tikzpicture}%
}}
\qquad=\qquad
\vcenter{\hbox{%
\begin{tikzpicture}[xscale=0.0352778, yscale=0.0352778, thin, inner sep=0]
  \def\p#1#2{%
    \ifcase #1
      \or \ifx#2x  0       \else  0 \fi
      \or \ifx#2x  \x1-20  \else  \y1-20 \fi
    \fi}
  \def\x#1{\p#1x}
  \def\y#1{\p#1y}
  \path[use as bounding box] (\x1+6,\y1+10) rectangle (\x1-26,\y1-36);
  \draw (\x1,\y1+10) -- (\x1,\y1);
  \draw (\x1-6,\y1-6)[rounded corners=4] -- (\x2,\y1-6) -- (\x2,\y2);
  \draw (\x1+6,\y1-6) -- (\x1+6,\y2-16);
  \draw (\x2+6,\y2-6)[rounded corners=4] -- (\x2+10,\y2-6) -- (\x2+10,\y2-16);
  \draw (\x2-6,\y2-6) -- (\x2-6,\y2-16);
  \putDuplicator(\x1,\y1)
  \putDuplicator(\x2,\y2)
\end{tikzpicture}%
}}
$$
Namely, the order of duplications does not matter.
Accordingly, it is sometimes convenient to introduce
 the following multi-duplicator:
$$
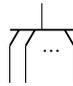
\vcenter{\hbox{%
\begin{tikzpicture}[xscale=0.0352778, yscale=0.0352778, thin, inner sep=0]
  \def\p#1#2{%
    \ifcase #1
      \or \ifx#2x   0         \else  0 \fi
    \fi}
  \def\x#1{\p#1x}
  \def\y#1{\p#1y}
 \begin{scope}[line width=.7]
  \draw (\x1,\y1) -- ++(12,0);
  \draw (\x1,\y1) -- ++(-12,0);
  \draw (\x1+8.6,\y1) -- ++(3.5,-6.0);
  \draw (\x1-8.6,\y1) -- ++(-3.5,-6.0);
  \draw (\x1-2.6,\y1) -- ++(-3.5,-6.0);
 \end{scope}
  \draw (\x1,\y1+10.0) -- (\x1,\y1);
  \draw (\x1-12.1,\y1-6.0) -- (\x1-12.1,\y1-20.0);
  \draw (\x1-6.1,\y1-6.0) -- (\x1-6.1,\y1-20.0);
  \draw (\x1+12.1,\y1-6.0) -- (\x1+12.1,\y1-20.0);
  \node at (\x1+3.5,\y1-8) {$\scriptstyle\cdots$};
\end{tikzpicture}%
}}%
$$
 by integrating successive occurrences of duplicators.
A {\it board} organizes a new part from an already constructed graph.
If we have a graph $G'$ having $m$ wires upward and a single
 wire downward, a board is given as
$$
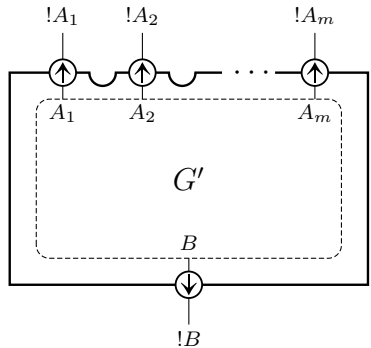
\vcenter{\hbox{%
\begin{tikzpicture}[xscale=0.0352778, yscale=0.0352778, thin, inner sep=0]
  \def\p#1#2{%
    \ifcase #1
      \or \ifx#2x  0       \else  0 \fi
      \or \ifx#2x  \x1+125 \else  \y1+80 \fi
      \or \ifx#2x  \x1+15  \else  \y2 \fi
      \or \ifx#2x  \x1+30  \else  \y2 \fi
      \or \ifx#2x  \x1+45  \else  \y2 \fi
      \or \ifx#2x  \x1+60  \else  \y2 \fi
      \or \ifx#2x  \x2-15  \else  \y2 \fi
    \fi}
  \def\x#1{\p#1x}
  \def\y#1{\p#1y}
 \draw[](\x3,\y3-10) -- (\x3,\y3+15);
 \draw[](\x5,\y5-10) -- (\x5,\y5+15);
 \draw[](\x7,\y7-10) -- (\x7,\y7+15);
 \draw[](\x1+62.5,\y1+10) -- (\x1+62.5,\y1-15);
 \boardBoundary(\x1-5,\y1)(\x2+5,\y2)
 \putLowerSocket(\x1+62.5,\y1)
 \putUpperSocket(\x3,\y3)
 \putDimple(\x4,\y4)
 \putUpperSocket(\x5,\y5)
 \putDimple(\x6,\y6)
 \putUpperSocket(\x7,\y7)
 \draw[white,line width=2.0] (\x1+75,\y2) -- ++(20,0);
 \draw[dash pattern=on 2 off 1,rounded corners=5] (\x1+5,\y1+10) rectangle
   (\x2-5,\y2-10);
 \node at (\x1+86,\y2) {$\cdots$};
 \node at ($ .5*(\x1,\y1)+.5*(\x2,\y2) $) {$G'$};
 \node[above,text depth=4] at (\x3,\y3+15) {$\scriptstyle \mathord!A_1$};
 \node[above,text depth=4] at (\x5,\y5+15) {$\scriptstyle \mathord!A_2$};
 \node[above,text depth=4] at (\x7,\y7+15) {$\scriptstyle \mathord!A_m$};
 \node[below,text height=8] at (\x1+62.5,\y1-15) {$\scriptstyle \mathord!B$};
 \node[below,text depth=4] at (\x3,\y3-12) {$\scriptstyle A_1$};
 \node[below,text depth=4] at (\x5,\y5-12) {$\scriptstyle A_2$};
 \node[below,text depth=4] at (\x7,\y7-12) {$\scriptstyle A_m$};
 \node[above,text height=8,inner sep=1] at (\x1+62.5,\y1+12) {$\scriptstyle B$};
\end{tikzpicture}%
}}%
$$
Here a dotted line signifies the subgraph $G'$ and it is not
 drawn in practice.
The thick arrows with surrounding circles located on the boundary are
 called the {\it gates} of the board.
The single gate on the bottom is a positive gate,
 and the ones on the top are negative gates.
The boards are the two-sided version of the boxes in proof-nets \cite{gira2}.

The switching condition is extended.
The $d$-part is similar to the tensor introduction, thus
 obeys the switching.
The board may be regarded as a sequence of the tensor introduction.
Hence the switching is not involved.

The morphisms $\delta_A,\varepsilon_A,d_A$, and $e_A$ are
 interpreted by the corresponding parts.
The structural morphisms
 $\tilde\varphi_{A,B}$ and $\varphi_0$ are interpreted by
$$
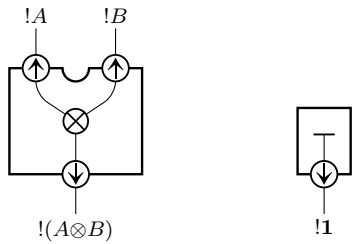
\begin{tikzpicture}[xscale=0.0352778, yscale=0.0352778, thin, inner sep=0]
  \def\p#1#2{%
    \ifcase #1
      \or \ifx#2x  0       \else  0 \fi
      \or \ifx#2x  \x1+25 \else  \y1+20 \fi
      \or \ifx#2x  \x1+50 \else  \y1+40 \fi
    \fi}
  \def\x#1{\p#1x}
  \def\y#1{\p#1y}
  \path[use as bounding box] (0,-26) rectangle (50,63);
  \draw (\x1+10,\y3+15)[rounded corners=5] -- (\x1+10,\y3-10)
    -- (\x2,\y2);
  \draw (\x1+40,\y3+15)[rounded corners=5] -- (\x1+40,\y3-10)
    -- (\x2,\y2);
  \draw (\x2,\y2) -- (\x2,\y1-15);
  \boardBoundary(\x1,\y1)(\x3,\y3)
  \putUpperSocket(\x1+10,\y3)
  \putDimple(\x1+25,\y3)
  \putUpperSocket(\x1+40,\y3)
  \putLowerSocket(\x2,\y1)
  \putTensor(\x2,\y2)
  \node[above,inner sep=2] at (\x1+10,\y3+15) {$\scriptstyle\mathord!A$};
  \node[above,inner sep=2] at (\x1+40,\y3+15) {$\scriptstyle\mathord!B$};
  \node[below,inner sep=2] at (\x2,\y1-15) {$\scriptstyle\mathord!(A\otimes B)$};
\end{tikzpicture}%
	\hskip5em
\begin{tikzpicture}[xscale=0.0352778, yscale=0.0352778, thin, inner sep=0]
  \def\p#1#2{%
    \ifcase #1
      \or \ifx#2x  0      \else  0 \fi
      \or \ifx#2x  \x1+10 \else  \y1+15 \fi
      \or \ifx#2x  \x1+20 \else  \y1+25 \fi
    \fi}
  \def\x#1{\p#1x}
  \def\y#1{\p#1y}
  \path[use as bounding box] (0,-26) rectangle (50,63);
  \draw (\x2,\y2) -- (\x2,\y1-15);
  \boardBoundary(\x1,\y1)(\x3,\y3)
  \putLowerSocket(\x2,\y1)
  \putPositiveTerminal(\x2,\y2)
  \node[below,inner sep=2] at (\x2,\y1-15) {$\scriptstyle\mathord!{\bf 1}$};
\end{tikzpicture}%
$$

The following give the $\beta$-elimination rules:
$$\vcenter
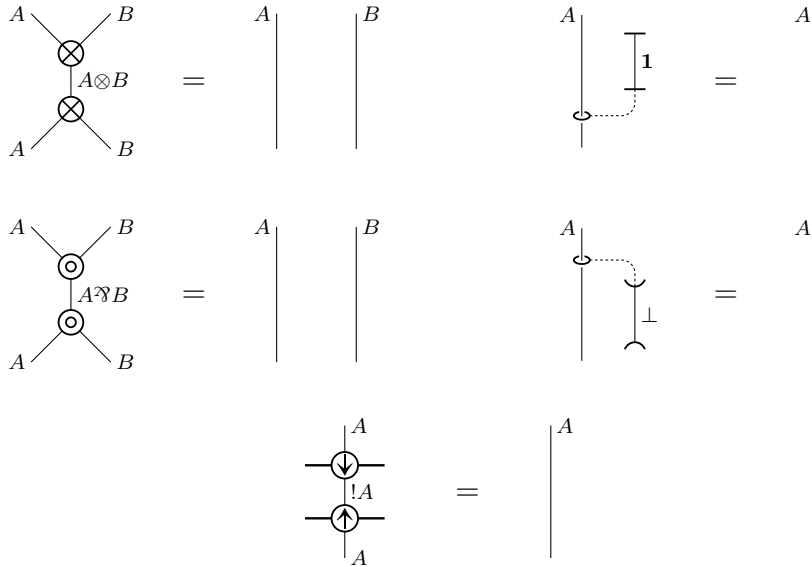
{\hbox{%
\begin{tikzpicture}[xscale=0.0352778, yscale=0.0352778, thin, inner sep=0]
  \def\p#1#2{%
    \ifcase #1
      \or \ifx#2x  0       \else  0 \fi
      \or \ifx#2x  \x1     \else  \y1+21 \fi
    \fi}
  \def\x#1{\p#1x}
  \def\y#1{\p#1y}
  \path[use as bounding box] (\x1-15,\y1-15) rectangle (\x2+15,\y2+15);
  \draw (\x1-15,\y1-15) -- (\x1,\y1);
  \draw (\x1+15,\y1-15) -- (\x1,\y1);
  \draw (\x2,\y2) -- (\x1,\y1);
  \draw (\x2-15,\y2+15) -- (\x2,\y2);
  \draw (\x2+15,\y2+15) -- (\x2,\y2);
  \putTensor(\x1,\y1)
  \putTensor(\x2,\y2)
  \node[right,inner sep=2] at ($ .5*(\x1,\y1)+.5*(\x2,\y2) $)
    {$\scriptstyle A\otimes B$};
  \node[left,inner sep=2] at (\x2-15,\y2+15) {$\scriptstyle A$};
  \node[right,inner sep=2] at (\x2+15,\y2+15) {$\scriptstyle B$};
  \node[left,inner sep=2] at (\x1-15,\y1-15) {$\scriptstyle A$};
  \node[right,inner sep=2] at (\x1+15,\y1-15) {$\scriptstyle B$};
\end{tikzpicture}}}
\qquad =\qquad
\vcenter{\hbox{%
\begin{tikzpicture}[xscale=0.0352778, yscale=0.0352778, thin, inner sep=0]
  \def\p#1#2{%
    \ifcase #1
      \or \ifx#2x  0       \else  0 \fi
      \or \ifx#2x  \x1     \else  \y1+21 \fi
    \fi}
  \def\x#1{\p#1x}
  \def\y#1{\p#1y}
  \path[use as bounding box] (\x1-15,\y1-15) rectangle (\x2+15,\y2+15);
  \draw (\x1-15,\y1-15) -- (\x2-15,\y2+15);
  \draw (\x1+15,\y1-15) -- (\x2+15,\y2+15);
  \node[left,inner sep=2] at (\x2-15,\y2+15) {$\scriptstyle A$};
  \node[right,inner sep=2] at (\x2+15,\y2+15) {$\scriptstyle B$};
\end{tikzpicture}}}
\kern7em
\vcenter{\hbox{%
\begin{tikzpicture}[xscale=0.0352778, yscale=0.0352778, thin, inner sep=0]
  \def\p#1#2{%
    \ifcase #1
      \or \ifx#2x  0       \else  0 \fi
    \fi}
  \def\x#1{\p#1x}
  \def\y#1{\p#1y}
  \path[use as bounding box] (\x1-23,\y1-29) rectangle (\x1+3,\y1+21);
  \draw (\x1,\y1+14) -- (\x1,\y1-7);
  \draw (\x1-20,\y1+21) -- (\x1-20,\y1-29);
  \putPositiveTerminal(\x1,\y1-7)
  \putPositiveTerminal(\x1,\y1+14)
  \putRubberBand(\x1-20,\y1-17)
  \draw[dash pattern=on 1 off 1] (\x1,\y1-7)[rounded corners=5] -- (\x1,\y1-17)
    -- (\x1-17,\y1-17);
  \node[right,inner sep=2] at ($ .5*(\x1,\y1-7)+.5*(\x1,\y1+14) $)
    {$\scriptstyle {\bf 1}$};
  \node[left,inner sep=2] at (\x1-20,\y1+21) {$\scriptstyle A$};
\end{tikzpicture}}}
\qquad=\qquad
\vcenter{\hbox{%
\begin{tikzpicture}[xscale=0.0352778, yscale=0.0352778, thin, inner sep=0]
  \def\p#1#2{%
    \ifcase #1
      \or \ifx#2x  0       \else  0 \fi
    \fi}
  \def\x#1{\p#1x}
  \def\y#1{\p#1y}
  \path[use as bounding box] (\x1-3,\y1-29) rectangle (\x1+3,\y1+21);
  \draw (\x1,\y1+21) -- (\x1,\y1-29);
  \node[left,inner sep=2] at (\x1,\y1+21) {$\scriptstyle A$};
\end{tikzpicture}}}
$$
	\vskip3ex
$$\vcenter{\hbox{%
\begin{tikzpicture}[xscale=0.0352778, yscale=0.0352778, thin, inner sep=0]
  \def\p#1#2{%
    \ifcase #1
      \or \ifx#2x  0       \else  0 \fi
      \or \ifx#2x  \x1     \else  \y1+21 \fi
    \fi}
  \def\x#1{\p#1x}
  \def\y#1{\p#1y}
  \path[use as bounding box] (\x1-15,\y1-15) rectangle (\x2+15,\y2+15);
  \draw (\x1-15,\y1-15) -- (\x1,\y1);
  \draw (\x1+15,\y1-15) -- (\x1,\y1);
  \draw (\x2,\y2) -- (\x1,\y1);
  \draw (\x2-15,\y2+15) -- (\x2,\y2);
  \draw (\x2+15,\y2+15) -- (\x2,\y2);
  \putPar(\x1,\y1)
  \putPar(\x2,\y2)
  \node[right,inner sep=2] at ($ .5*(\x1,\y1)+.5*(\x2,\y2) $)
    {$\scriptstyle A\slinpar B$};
  \node[left,inner sep=2] at (\x2-15,\y2+15) {$\scriptstyle A$};
  \node[right,inner sep=2] at (\x2+15,\y2+15) {$\scriptstyle B$};
  \node[left,inner sep=2] at (\x1-15,\y1-15) {$\scriptstyle A$};
  \node[right,inner sep=2] at (\x1+15,\y1-15) {$\scriptstyle B$};
\end{tikzpicture}}}
\qquad =\qquad
\vcenter{\hbox{%
\begin{tikzpicture}[xscale=0.0352778, yscale=0.0352778, thin, inner sep=0]
  \def\p#1#2{%
    \ifcase #1
      \or \ifx#2x  0       \else  0 \fi
      \or \ifx#2x  \x1     \else  \y1+21 \fi
    \fi}
  \def\x#1{\p#1x}
  \def\y#1{\p#1y}
  \path[use as bounding box] (\x1-15,\y1-15) rectangle (\x2+15,\y2+15);
  \draw (\x1-15,\y1-15) -- (\x2-15,\y2+15);
  \draw (\x1+15,\y1-15) -- (\x2+15,\y2+15);
  \node[left,inner sep=2] at (\x2-15,\y2+15) {$\scriptstyle A$};
  \node[right,inner sep=2] at (\x2+15,\y2+15) {$\scriptstyle B$};
\end{tikzpicture}}}
\kern7em
\vcenter{\hbox{%
\begin{tikzpicture}[xscale=0.0352778, yscale=0.0352778, thin, inner sep=0]
  \def\p#1#2{%
    \ifcase #1
      \or \ifx#2x  0       \else  0 \fi
    \fi}
  \def\x#1{\p#1x}
  \def\y#1{\p#1y}
  \path[use as bounding box] (\x1-23,\y1+29) rectangle (\x1+3,\y1-21);
  \draw (\x1,\y1-14) -- (\x1,\y1+7);
  \draw (\x1-20,\y1-21) -- (\x1-20,\y1+29);
  \putUpperNegativeTerminal(\x1,\y1+7)
  \putLowerNegativeTerminal(\x1,\y1-14)
  \putRubberBand(\x1-20,\y1+17)
  \draw[dash pattern=on 1 off 1] (\x1,\y1+7)[rounded corners=5] -- (\x1,\y1+17)
    -- (\x1-17,\y1+17);
  \node[right,inner sep=2] at ($ .5*(\x1,\y1+7)+.5*(\x1,\y1-14) $)
    {$\scriptstyle \bot$};
  \node[left,inner sep=2] at (\x1-20,\y1+29) {$\scriptstyle A$};
\end{tikzpicture}}}
\qquad=\qquad
\vcenter{\hbox{%
\begin{tikzpicture}[xscale=0.0352778, yscale=0.0352778, thin, inner sep=0]
  \def\p#1#2{%
    \ifcase #1
      \or \ifx#2x  0       \else  0 \fi
    \fi}
  \def\x#1{\p#1x}
  \def\y#1{\p#1y}
  \path[use as bounding box] (\x1-3,\y1+29) rectangle (\x1+3,\y1-21);
  \draw (\x1,\y1-21) -- (\x1,\y1+29);
  \node[left,inner sep=2] at (\x1,\y1+29) {$\scriptstyle A$};
\end{tikzpicture}}}
$$
	\vskip3ex
$$
\vcenter{\hbox{%
\begin{tikzpicture}[xscale=0.0352778, yscale=0.0352778, thin, inner sep=0]
  \def\p#1#2{%
    \ifcase #1
      \or \ifx#2x  0      \else  0 \fi
      \or \ifx#2x  \x1    \else  \y1+20 \fi
    \fi}
  \def\x#1{\p#1x}
  \def\y#1{\p#1y}
  \path[use as bounding box] (-15,-15) rectangle (15,35);
  \draw (\x2,\y2+15) -- (\x1,\y1-15);
  \draw[line width=.9] (\x1-15,\y1) -- (\x1+15,\y1);
  \draw[line width=.9] (\x2-15,\y2) -- (\x2+15,\y2);
  \putUpperSocket(\x1,\y1)
  \putLowerSocket(\x2,\y2)
  \node[right,inner sep=2] at ($ .5*(\x1,\y1)+.5*(\x2,\y2) $)
    {$\scriptstyle\mathord!A$};
  \node[right,inner sep=2] at (\x2,\y2+15) {$\scriptstyle A$};
  \node[right,inner sep=2] at (\x1,\y1-15) {$\scriptstyle A$};
\end{tikzpicture}%
}}
\qquad=\qquad
\vcenter{\hbox{%
\begin{tikzpicture}[xscale=0.0352778, yscale=0.0352778, thin, inner sep=0]
  \def\p#1#2{%
    \ifcase #1
      \or \ifx#2x  0      \else  0 \fi
      \or \ifx#2x  \x1    \else  \y1+20 \fi
    \fi}
  \def\x#1{\p#1x}
  \def\y#1{\p#1y}
  \path[use as bounding box] (0,-15) rectangle (0,35);
  \draw (\x2,\y2+15) -- (\x1,\y1-15);
  \node[right,inner sep=2] at (\x2,\y2+15) {$\scriptstyle A$};
\end{tikzpicture}%
}}
$$
The elimination contracts the left-hand side to the right-hand side.
The contraction is applied
 whenever they are possible.
We omit the $\beta$-rule for the duality parts.
It does not occur
 since we are concerned only with normal forms, in which
 duality $\beta$-redexes have been contracted.
The following are the $\eta$-expansion rules:
$$\vcenter
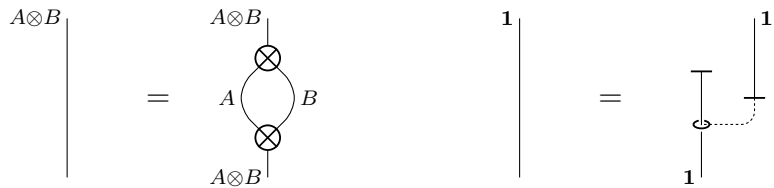
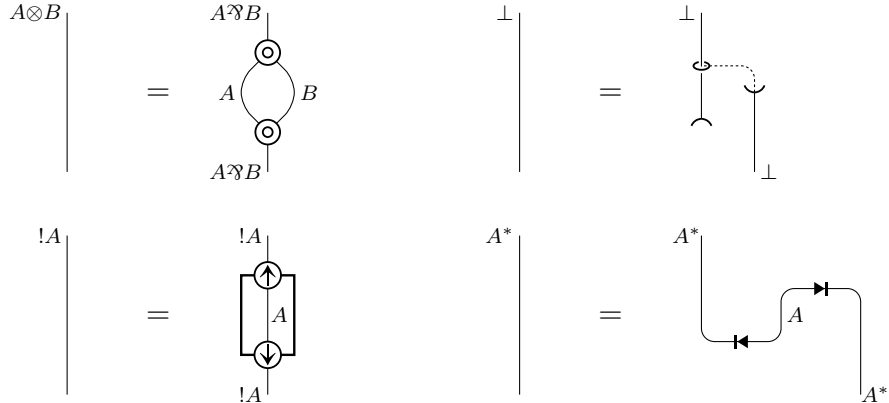
{\hbox{%
\begin{tikzpicture}[xscale=0.0352778, yscale=0.0352778, thin, inner sep=0]
  \def\p#1#2{%
    \ifcase #1
      \or \ifx#2x  0       \else  0 \fi
      \or \ifx#2x  \x1     \else  \y1+30 \fi
    \fi}
  \def\x#1{\p#1x}
  \def\y#1{\p#1y}
  \path[use as bounding box] (\x1-3,\y1-15) rectangle (\x2+3,\y2+15);
  \draw (\x1,\y1-15) -- (\x2,\y2+15);
  \node[left,inner sep=2] at (\x2,\y2+15) {$\scriptstyle A\otimes B$};
\end{tikzpicture}}}
\qquad=\qquad
\vcenter{\hbox{%
\begin{tikzpicture}[xscale=0.0352778, yscale=0.0352778, thin, inner sep=0]
  \def\p#1#2{%
    \ifcase #1
      \or \ifx#2x  0       \else  0 \fi
      \or \ifx#2x  \x1     \else  \y1+30 \fi
    \fi}
  \def\x#1{\p#1x}
  \def\y#1{\p#1y}
  \path[use as bounding box] (\x1-10,\y1-15) rectangle (\x2+10,\y2+15);
  \draw (\x1,\y1)[rounded corners=5] -- (\x1-10,\y1+10) -- (\x1-10,\y2-10) -- (\x2,\y2);
  \draw (\x1,\y1)[rounded corners=5] -- (\x1+10,\y1+10) -- (\x1+10,\y2-10) -- (\x2,\y2);
  \draw (\x1,\y1) -- (\x1,\y1-15);
  \draw (\x2,\y2) -- (\x2,\y2+15);
  \putTensor(\x1,\y1)
  \putTensor(\x2,\y2)
  \node[left,inner sep=2] at (\x1,\y1-15) {$\scriptstyle A\otimes B$};
  \node[left,inner sep=2] at (\x2,\y2+15) {$\scriptstyle A\otimes B$};
  \node[left,inner sep=2] at ($ .5*(\x1-10,\y1)+.5*(\x2-10,\y2) $)
    {$\scriptstyle A$};
  \node[right,inner sep=2] at ($ .5*(\x1+10,\y1)+.5*(\x2+10,\y2) $)
    {$\scriptstyle B$};
\end{tikzpicture}}}
\kern7em
\vcenter{\hbox{%
\begin{tikzpicture}[xscale=0.0352778, yscale=0.0352778, thin, inner sep=0]
  \def\p#1#2{%
    \ifcase #1
      \or \ifx#2x  0       \else  0 \fi
      \or \ifx#2x  \x1     \else  \y1+30 \fi
    \fi}
  \def\x#1{\p#1x}
  \def\y#1{\p#1y}
  \path[use as bounding box] (\x1-3,\y1-15) rectangle (\x2+3,\y2+15);
  \draw (\x1,\y1-15) -- (\x2,\y2+15);
  \node[left,inner sep=2] at (\x2,\y2+15) {$\scriptstyle {\bf 1}$};
\end{tikzpicture}}}
\qquad=\qquad
\vcenter{\hbox{%
\begin{tikzpicture}[xscale=0.0352778, yscale=0.0352778, thin, inner sep=0]
  \def\p#1#2{%
    \ifcase #1
      \or \ifx#2x  0       \else  0 \fi
    \fi}
  \def\x#1{\p#1x}
  \def\y#1{\p#1y}
  \path[use as bounding box] (\x1-3,\y1-30) rectangle (\x1+20,\y1+30);
  \draw (\x1,\y1+10) -- (\x1,\y1-30);
  \draw (\x1+20,\y1) -- (\x1+20,\y1+30);
  \putPositiveTerminal(\x1,\y1+10)
  \putPositiveTerminal(\x1+20,\y1)
  \putRubberBand(\x1,\y1-10)
  \draw[dash pattern=on 1 off 1] (\x1+20,\y1)[rounded corners=5] -- (\x1+20,\y1-10)
    -- (\x1,\y1-10);
  \node[right,inner sep=2] at (\x1+20,\y1+30) {$\scriptstyle {\bf 1}$};
  \node[left,inner sep=2] at (\x1,\y1-30) {$\scriptstyle {\bf 1}$};
\end{tikzpicture}}}
$$
	\vskip3ex
$$\vcenter{\hbox{%
\begin{tikzpicture}[xscale=0.0352778, yscale=0.0352778, thin, inner sep=0]
  \def\p#1#2{%
    \ifcase #1
      \or \ifx#2x  0       \else  0 \fi
      \or \ifx#2x  \x1     \else  \y1+30 \fi
    \fi}
  \def\x#1{\p#1x}
  \def\y#1{\p#1y}
  \path[use as bounding box] (\x1-3,\y1-15) rectangle (\x2+3,\y2+15);
  \draw (\x1,\y1-15) -- (\x2,\y2+15);
  \node[left,inner sep=2] at (\x2,\y2+15) {$\scriptstyle A\otimes B$};
\end{tikzpicture}}}
\qquad=\qquad
\vcenter{\hbox{%
\begin{tikzpicture}[xscale=0.0352778, yscale=0.0352778, thin, inner sep=0]
  \def\p#1#2{%
    \ifcase #1
      \or \ifx#2x  0       \else  0 \fi
      \or \ifx#2x  \x1     \else  \y1+30 \fi
    \fi}
  \def\x#1{\p#1x}
  \def\y#1{\p#1y}
  \path[use as bounding box] (\x1-10,\y1-15) rectangle (\x2+10,\y2+15);
  \draw (\x1,\y1)[rounded corners=5] -- (\x1-10,\y1+10) -- (\x1-10,\y2-10) -- (\x2,\y2);
  \draw (\x1,\y1)[rounded corners=5] -- (\x1+10,\y1+10) -- (\x1+10,\y2-10) -- (\x2,\y2);
  \draw (\x1,\y1) -- (\x1,\y1-15);
  \draw (\x2,\y2) -- (\x2,\y2+15);
  \putPar(\x1,\y1)
  \putPar(\x2,\y2)
  \node[left,inner sep=2] at (\x1,\y1-15) {$\scriptstyle A\slinpar B$};
  \node[left,inner sep=2] at (\x2,\y2+15) {$\scriptstyle A\slinpar B$};
  \node[left,inner sep=2] at ($ .5*(\x1-10,\y1)+.5*(\x2-10,\y2) $)
    {$\scriptstyle A$};
  \node[right,inner sep=2] at ($ .5*(\x1+10,\y1)+.5*(\x2+10,\y2) $)
    {$\scriptstyle B$};
\end{tikzpicture}}}
\kern7em
\vcenter{\hbox{%
\begin{tikzpicture}[xscale=0.0352778, yscale=0.0352778, thin, inner sep=0]
  \def\p#1#2{%
    \ifcase #1
      \or \ifx#2x  0       \else  0 \fi
      \or \ifx#2x  \x1     \else  \y1+30 \fi
    \fi}
  \def\x#1{\p#1x}
  \def\y#1{\p#1y}
  \path[use as bounding box] (\x1-3,\y1-15) rectangle (\x2+3,\y2+15);
  \draw (\x1,\y1-15) -- (\x2,\y2+15);
  \node[left,inner sep=2] at (\x2,\y2+15) {$\scriptstyle \bot$};
\end{tikzpicture}}}
\qquad=\qquad
\vcenter{\hbox{%
\begin{tikzpicture}[xscale=0.0352778, yscale=0.0352778, thin, inner sep=0]
  \def\p#1#2{%
    \ifcase #1
      \or \ifx#2x  0       \else  0 \fi
    \fi}
  \def\x#1{\p#1x}
  \def\y#1{\p#1y}
  \path[use as bounding box] (\x1-3,\y1-30) rectangle (\x1+20,\y1+30);
  \draw (\x1,\y1-10) -- (\x1,\y1+30);
  \draw (\x1+20,\y1) -- (\x1+20,\y1-30);
  \putLowerNegativeTerminal(\x1,\y1-10)
  \putUpperNegativeTerminal(\x1+20,\y1)
  \putRubberBand(\x1,\y1+10)
  \draw[dash pattern=on 1 off 1] (\x1+20,\y1)[rounded corners=5] -- (\x1+20,\y1+10)
    -- (\x1,\y1+10);
  \node[right,inner sep=2] at (\x1+20,\y1-30) {$\scriptstyle \bot$};
  \node[left,inner sep=2] at (\x1,\y1+30) {$\scriptstyle \bot$};
\end{tikzpicture}}}
$$
	\vskip3ex
$$\vcenter{\hbox{%
\begin{tikzpicture}[xscale=0.0352778, yscale=0.0352778, thin, inner sep=0]
  \def\p#1#2{%
    \ifcase #1
      \or \ifx#2x  0       \else  0 \fi
      \or \ifx#2x  \x1     \else  \y1+30 \fi
    \fi}
  \def\x#1{\p#1x}
  \def\y#1{\p#1y}
  \path[use as bounding box] (\x1-3,\y1-15) rectangle (\x2+3,\y2+15);
  \draw (\x1,\y1-15) -- (\x2,\y2+15);
  \node[left,inner sep=2] at (\x2,\y2+15) {$\scriptstyle \mathord!A$};
\end{tikzpicture}}}
\qquad=\qquad
\vcenter{\hbox{%
\begin{tikzpicture}[xscale=0.0352778, yscale=0.0352778, thin, inner sep=0]
  \def\p#1#2{%
    \ifcase #1
      \or \ifx#2x  0       \else  0 \fi
      \or \ifx#2x  \x1     \else  \y1+30 \fi
    \fi}
  \def\x#1{\p#1x}
  \def\y#1{\p#1y}
  \path[use as bounding box] (\x1-10,\y1-15) rectangle (\x2+10,\y2+15);
  \boardBoundary(\x1-10,\y1)(\x2+10,\y2)
  \draw (\x1,\y1-15) -- (\x2,\y2+15);
  \putLowerSocket(\x1,\y1)
  \putUpperSocket(\x2,\y2)
  \node[left,inner sep=2] at (\x2,\y2+15) {$\scriptstyle \mathord!A$};
  \node[left,inner sep=2] at (\x1,\y1-15) {$\scriptstyle \mathord!A$};
  \node[right,inner sep=1] at ($ .5*(\x1,\y1)+.5*(\x2,\y2) $)
    {$\scriptstyle A$};
\end{tikzpicture}}}
\kern7em
\vcenter{\hbox{%
\begin{tikzpicture}[xscale=0.0352778, yscale=0.0352778, thin, inner sep=0]
  \def\p#1#2{%
    \ifcase #1
      \or \ifx#2x  0       \else  0 \fi
      \or \ifx#2x  \x1     \else  \y1+30 \fi
    \fi}
  \def\x#1{\p#1x}
  \def\y#1{\p#1y}
  \path[use as bounding box] (\x1-3,\y1-15) rectangle (\x2+3,\y2+15);
  \draw (\x1,\y1-15) -- (\x2,\y2+15);
  \node[left,inner sep=2] at (\x2,\y2+15) {$\scriptstyle A^*$};
\end{tikzpicture}}}
\qquad=\qquad
\vcenter{\hbox{%
\begin{tikzpicture}[xscale=0.0352778, yscale=0.0352778, thin, inner sep=0]
  \def\p#1#2{%
    \ifcase #1
      \or \ifx#2x  0       \else  0 \fi
      \or \ifx#2x  \x1+15  \else  \y1-10 \fi
      \or \ifx#2x  \x2+15  \else  \y1 \fi
      \or \ifx#2x  \x3+15  \else  \y3+10 \fi
      \or \ifx#2x  \x4+15  \else  \y1 \fi
    \fi}
  \def\x#1{\p#1x}
  \def\y#1{\p#1y}
  \path[use as bounding box] (\x1-3,\y1-30) rectangle (\x1+20,\y1+30);
  \draw (\x1,\y1+30)[rounded corners=5] -- (\x1,\y1-10) -- (\x3,\y3-10)
    -- (\x3,\y3+10) -- (\x5,\y5+10) -- (\x5,\y5-30);
  \putLeftDiode(\x2,\y2)
  \putRightDiode(\x4,\y4)
  \node[right,inner sep=2] at (\x3,\y3) {$\scriptstyle A$};
  \node[right,inner sep=0] at (\x5,\y5-30) {$\scriptstyle A^*$};
  \node[left,inner sep=0] at (\x1,\y1+30) {$\scriptstyle A^*$};
\end{tikzpicture}}}
$$
The expansion rules transform the left graphs to the
 right graphs.
The rules are applied as far as they do not create new
 redexes.
First, if the $\eta$-expansion creates a
 new $\beta$-redex (including those of the duality),
 the expansion is not applied.
Second, the rule for $\mathord!A$
 is not applied immdiately above
 (i.e., the flat side of) a $\delta$-part,
 an $\varepsilon$-part, a $d$-part, or an $e$-part,
 since it creates a new naturality redex, i.e.,
 a redex of rules (18) to (21).

The equivalence between morphisms is absorbed
 by graphical equivalence under the $\beta$-rules.
For example, the naturality of $\tilde\varphi$ turns
 out to be the following sequence of equivalences:
$$
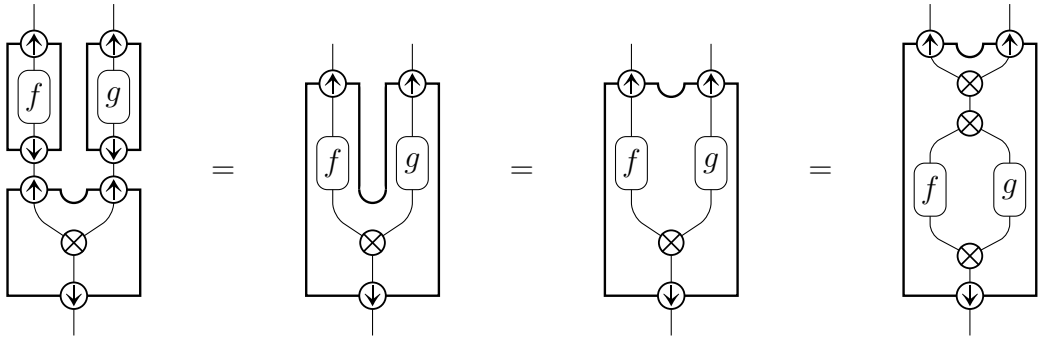
\vcenter{\hbox{%
\begin{tikzpicture}[xscale=0.0352778, yscale=0.0352778, thin, inner sep=0]
  \def\p#1#2{%
    \ifcase #1
      \or \ifx#2x  0      \else  0 \fi
      \or \ifx#2x  \x1    \else  \y1+20 \fi
      \or \ifx#2x  \x1-15 \else  \y1+40 \fi
      \or \ifx#2x  \x1+15 \else  \y3 \fi
      \or \ifx#2x  \x1-15 \else  \y3+15 \fi
      \or \ifx#2x  \x1+15 \else  \y5 \fi
      \or \ifx#2x  \x1-15 \else  \y5+40 \fi
      \or \ifx#2x  \x1+15 \else  \y7 \fi
    \fi}
  \def\x#1{\p#1x}
  \def\y#1{\p#1y}
  \path[use as bounding box] (-25,-15) rectangle (25,110);
  \draw (\x2,\y2) -- (\x2,\y1-15);
  \draw (\x5,\y5+10)[rounded corners=5] -- (\x3,\y3-10) -- (\x2,\y2);
  \draw (\x6,\y6+10)[rounded corners=5] -- (\x4,\y4-10) -- (\x2,\y2);
  \draw (\x7,\y7+15) -- (\x7,\y7-10);
  \draw (\x8,\y8+15) -- (\x8,\y8-10);
  \draw[rounded corners=4] (\x5-6,\y5+10) rectangle
    (\x7+6,\y7-10);
  \draw[rounded corners=4] (\x6-6,\y6+10) rectangle
    (\x8+6,\y8-10);
  \boardBoundary(\x3-10,\y1)(\x4+10,\y4)
  \boardBoundary(\x5-10,\y5)(\x7+10,\y7)
  \boardBoundary(\x6-10,\y6)(\x8+10,\y8)
  \putUpperSocket(\x3,\y3)
  \putDimple(\x3+15,\y3)
  \putUpperSocket(\x4,\y4)
  \putLowerSocket(\x1,\y1)
  \putLowerSocket(\x5,\y5)
  \putLowerSocket(\x6,\y6)
  \putUpperSocket(\x7,\y7)
  \putUpperSocket(\x8,\y8)
  \putTensor(\x2,\y2)
  \node at (\x5,\y5+20) {$f$};
  \node at (\x6,\y6+20) {$g$};
\end{tikzpicture}%
}}
\qquad=\qquad
\vcenter{\hbox{%
\begin{tikzpicture}[xscale=0.0352778, yscale=0.0352778, thin, inner sep=0]
  \def\p#1#2{%
    \ifcase #1
      \or \ifx#2x  0      \else  0 \fi
      \or \ifx#2x  \x1    \else  \y1+20 \fi
      \or \ifx#2x  \x1-15 \else  \y1+80 \fi
      \or \ifx#2x  \x1+15 \else  \y3 \fi
      \or \ifx#2x  \x1-15 \else  \y1+40 \fi
      \or \ifx#2x  \x1+15 \else  \y5 \fi
    \fi}
  \def\x#1{\p#1x}
  \def\y#1{\p#1y}
  \path[use as bounding box] (-25,-15) rectangle (25,110);
  \draw (\x2,\y2) -- (\x2,\y1-15);
  \draw (\x5,\y5)[rounded corners=5] -- (\x5,\y5-10) -- (\x2,\y2);
  \draw (\x6,\y6)[rounded corners=5] -- (\x6,\y6-10) -- (\x2,\y2);
  \draw (\x3,\y3+15) -- (\x3,\y3-20);
  \draw (\x4,\y4+15) -- (\x4,\y4-20);
  \draw[rounded corners=4] (\x5-6,\y5) rectangle
    (\x3+6,\y3-20);
  \draw[rounded corners=4] (\x6-6,\y6) rectangle
    (\x4+6,\y4-20);
  \boardBoundary(\x3-10,\y1)(\x4+10,\y4)
  \putLowerSocket(\x1,\y1)
  \putUpperSocket(\x3,\y3)
  \putUpperSocket(\x4,\y4)
  \draw[line width=.9] (\x3+10,\y3) -- (\x5+10,\y5);
  \draw[line width=.9] (\x4-10,\y4) -- (\x6-10,\y6);
  \fill[white] (\x3+15,\y3) circle (4.6);
  \putDimple(\x5+15,\y5)
  \putTensor(\x2,\y2)
  \node at (\x5,\y5+10) {$f$};
  \node at (\x6,\y6+10) {$g$};
\end{tikzpicture}%
}}
\qquad=\qquad
\vcenter{\hbox{%
\begin{tikzpicture}[xscale=0.0352778, yscale=0.0352778, thin, inner sep=0]
  \def\p#1#2{%
    \ifcase #1
      \or \ifx#2x  0      \else  0 \fi
      \or \ifx#2x  \x1    \else  \y1+20 \fi
      \or \ifx#2x  \x1-15 \else  \y1+80 \fi
      \or \ifx#2x  \x1+15 \else  \y3 \fi
      \or \ifx#2x  \x1-15 \else  \y1+40 \fi
      \or \ifx#2x  \x1+15 \else  \y5 \fi
    \fi}
  \def\x#1{\p#1x}
  \def\y#1{\p#1y}
  \path[use as bounding box] (-25,-15) rectangle (25,110);
  \draw (\x2,\y2) -- (\x2,\y1-15);
  \draw (\x5,\y5)[rounded corners=5] -- (\x5,\y5-10) -- (\x2,\y2);
  \draw (\x6,\y6)[rounded corners=5] -- (\x6,\y6-10) -- (\x2,\y2);
  \draw (\x3,\y3+15) -- (\x3,\y3-20);
  \draw (\x4,\y4+15) -- (\x4,\y4-20);
  \draw[rounded corners=4] (\x5-6,\y5) rectangle
    (\x3+6,\y3-20);
  \draw[rounded corners=4] (\x6-6,\y6) rectangle
    (\x4+6,\y4-20);
  \boardBoundary(\x3-10,\y1)(\x4+10,\y4)
  \putLowerSocket(\x1,\y1)
  \putUpperSocket(\x3,\y3)
  \putUpperSocket(\x4,\y4)
  \putDimple(\x3+15,\y3)
  \putTensor(\x2,\y2)
  \node at (\x5,\y5+10) {$f$};
  \node at (\x6,\y6+10) {$g$};
\end{tikzpicture}%
}}
\qquad=\qquad
\vcenter{\hbox{%
\begin{tikzpicture}[xscale=0.0352778, yscale=0.0352778, thin, inner sep=0]
  \def\p#1#2{%
    \ifcase #1
      \or \ifx#2x  0      \else  0 \fi
      \or \ifx#2x  \x1    \else  \y1+15 \fi
      \or \ifx#2x  \x1    \else  \y2+50 \fi
      \or \ifx#2x  \x1    \else  \y3+15 \fi
      \or \ifx#2x  \x1-15 \else  \y4+15 \fi
      \or \ifx#2x  \x1+15 \else  \y5 \fi
    \fi}
  \def\x#1{\p#1x}
  \def\y#1{\p#1y}
  \path[use as bounding box] (-25,-15) rectangle (25,110);
  \draw (\x2,\y2) -- (\x2,\y1-15);
  \draw (\x2-15,\y2+15)[rounded corners=4] -- (\x2-15,\y2+8) -- (\x2,\y2);
  \draw (\x2+15,\y2+15)[rounded corners=4] -- (\x2+15,\y2+8) -- (\x2,\y2);
  \draw (\x3-15,\y3-15)[rounded corners=4] -- (\x3-15,\y3-8) -- (\x3,\y3);
  \draw (\x3+15,\y3-15)[rounded corners=4] -- (\x3+15,\y3-8) -- (\x3,\y3);
  \draw (\x4,\y4) -- (\x3,\y3);
  \draw (\x5,\y5+15)[rounded corners=4] -- (\x5,\y4+8) -- (\x4,\y4);
  \draw (\x6,\y6+15)[rounded corners=4] -- (\x6,\y4+8) -- (\x4,\y4);
  \draw[rounded corners=4] (\x2-21,\y2+15) rectangle
    (\x2-9,\y3-15);
  \draw[rounded corners=4] (\x2+21,\y2+15) rectangle
    (\x2+9,\y3-15);
  \boardBoundary(\x1-25,\y1)(\x1+25,\y6)
  \putLowerSocket(\x1,\y1)
  \putTensor(\x2,\y2)
  \putTensor(\x3,\y3)
  \putTensor(\x4,\y4)
  \putUpperSocket(\x5,\y5)
  \putUpperSocket(\x6,\y6)
  \putDimple(\x5+15,\y5)
  \node at (\x2-15,\y2+25) {$f$};
  \node at (\x2+15,\y2+25) {$g$};
\end{tikzpicture}%
}}
$$
Here the first and the last equalities are $\beta$-rules
 and the middle simply changes the depth of the slit
(we implicitly assume that the depths of slits
 do not matter).
For the coherences involved in ${\bf 1}$ and $\bot$,
 the equivalence rules for the $*$-autonomous case
 should be appropriately extended.
We leave these to the reader since the results in
 this paper are not concerned with them.

We transform only normal forms into graphs.
Hence not all combinations are legitimate.
For example, an $e$-part is not allowed to
 be linked to a leg of a $d$-part.
In particular, by the rewriting rules (1) to (7)
 and the convention of the $\eta$-expansions,
 $\delta$-parts always occur in
 the shape
$$
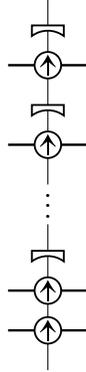
\vcenter{\hbox{%
\begin{tikzpicture}[xscale=0.0352778, yscale=0.0352778, thin, inner sep=0]
  \def\p#1#2{%
    \ifcase #1
      \or \ifx#2x  0      \else  0 \fi
      \or \ifx#2x  \x1    \else  \y1+15 \fi
      \or \ifx#2x  \x2    \else  \y2+15 \fi
      \or \ifx#2x  \x3    \else  \y3+40 \fi
      \or \ifx#2x  \x4    \else  \y4+15 \fi
      \or \ifx#2x  \x5    \else  \y5+15 \fi
      \or \ifx#2x  \x6    \else  \y6+15 \fi
    \fi}
  \def\x#1{\p#1x}
  \def\y#1{\p#1y}
  \path[use as bounding box] (-15,-15) rectangle (15,125);
  \draw (\x1,\y1-15) -- (\x3,\y3+10);
  \draw (\x4,\y4-15) -- (\x7,\y7+10);
  \draw[line width=.9] (\x1-15,\y1) -- (\x1+15,\y1);
  \draw[line width=.9] (\x2-15,\y2) -- (\x2+15,\y2);
  \draw[line width=.9] (\x4-15,\y4) -- (\x4+15,\y4);
  \draw[line width=.9] (\x6-15,\y6) -- (\x6+15,\y6);
  \putUpperSocket(\x1,\y1)
  \putUpperSocket(\x2,\y2)
  \putConcave(\x3,\y3)
  \putUpperSocket(\x4,\y4)
  \putConcave(\x5,\y5)
  \putUpperSocket(\x6,\y6)
  \putConcave(\x7,\y7)
  \node at (\x3,\y3+20) {$\vdots$};
\end{tikzpicture}%
}}
$$
Namely, the leg of each $\delta$-part must be linked to
 a negative gate of a board.
Moreover, the lowest $\delta$-part must have
 two negative gates in a row.
This corresponds to
$\delta;\mathord!\delta;\cdots;\mathord!^{n-1}\delta$.

Remark:
We explain why dimples are put on the upper edge of
 a board.
We recall that we associate graphs only with normal forms.
For a moment, however, suppose that we are trying to develop
 a graphical rewriting system, and let us see
 how the $(\varphi_0;\tilde\varphi)$-type rule (17)
 looks.
It contracts ${\bf 1}\otimes \mathord!A\mor{\varphi_0\cdot}
\mathord!{\bf 1}\otimes \mathord!A\mor{\tilde\varphi}
\mathord!({\bf 1}\otimes A)$.
This sequence of morphisms translates into
$$
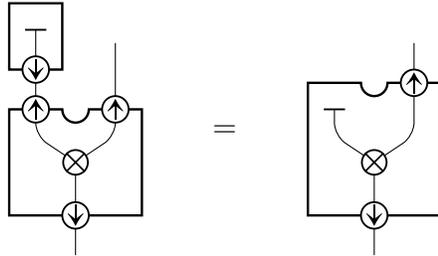
\vcenter{\hbox{%
\begin{tikzpicture}[xscale=0.0352778, yscale=0.0352778, thin, inner sep=0]
  \def\p#1#2{%
    \ifcase #1
      \or \ifx#2x  0      \else  0 \fi
      \or \ifx#2x  \x1    \else  \y1+20 \fi
      \or \ifx#2x  \x1-15 \else  \y1+40 \fi
      \or \ifx#2x  \x1+15 \else  \y3 \fi
      \or \ifx#2x  \x3    \else  \y3+15 \fi
      \or \ifx#2x  \x5    \else  \y5+15 \fi
      \or \ifx#2x  \x6    \else  \y6+10 \fi
    \fi}
  \def\x#1{\p#1x}
  \def\y#1{\p#1y}
  \path[use as bounding box] (-25,-15) rectangle (25,80);
  \draw (\x2,\y2) -- (\x2,\y1-15);
  \draw (\x5,\y5+15)[rounded corners=5] -- (\x3,\y3-10) -- (\x2,\y2);
  \draw (\x4,\y4+25)[rounded corners=5] -- (\x4,\y4-10) -- (\x2,\y2);
  \boardBoundary(\x3-10,\y1)(\x4+10,\y4)
  \boardBoundary(\x5-10,\y5)(\x7+10,\y7)
  \putUpperSocket(\x3,\y3)
  \putDimple(\x3+15,\y3)
  \putLowerSocket(\x1,\y1)
  \putUpperSocket(\x4,\y4)
  \putLowerSocket(\x5,\y5)
  \putPositiveTerminal(\x6,\y6)
  \putTensor(\x2,\y2)
\end{tikzpicture}%
}}
\qquad=\qquad
\vcenter{\hbox{%
\begin{tikzpicture}[xscale=0.0352778, yscale=0.0352778, thin, inner sep=0]
  \def\p#1#2{%
    \ifcase #1
      \or \ifx#2x  0      \else  0 \fi
      \or \ifx#2x  \x1    \else  \y1+20 \fi
      \or \ifx#2x  \x1-15 \else  \y1+40 \fi
      \or \ifx#2x  \x1+15 \else  \y3+10 \fi
    \fi}
  \def\x#1{\p#1x}
  \def\y#1{\p#1y}
  \path[use as bounding box] (-25,-15) rectangle (25,80);
  \draw (\x2,\y2) -- (\x2,\y1-15);
  \draw (\x3,\y3)[rounded corners=5] -- (\x3,\y2+10) -- (\x2,\y2);
  \draw (\x4,\y4+15)[rounded corners=5] -- (\x4,\y2+10) -- (\x2,\y2);
  \boardBoundary(\x3-10,\y1)(\x4+10,\y4)
  \putUpperSocket(\x4,\y4)
  \putDimple(\x4-15,\y4)
  \putLowerSocket(\x1,\y1)
  \putPositiveTerminal(\x3,\y3)
  \putTensor(\x2,\y2)
\end{tikzpicture}%
}}
$$
Observe that the dimple is not sandwiched between
 two gates.
In contrast, every dimple has a gate on each side
 if the graph comes from a normal form.
Namely, the rule (17) would be regarded as a rule
 eliminating an isolated dimple that does not lay between
 two gates.

We can restore a morphism from a legitimate graph by the process
 of sequentialization.
It is naturally defined by extending the definition
 in \cite{bcst}.
The process is not unique but yields equivalent
 morphisms regardless of the choice.
Although we do not give details here,
 let us address a point not stressed in \cite{bcst}.
A crossing of wires does not make sense in its own right.
We have two symmetries $A\otimes B\cong B\otimes A$ and
 $A\linpar B\cong B\linpar A$, both represented by crossing.
Indeed, some crossing is unable to be interpreted by
 any of these.
Suppose that we have
$$
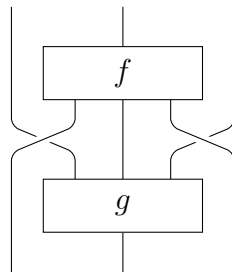
\begin{tikzpicture}[xscale=0.0352778, yscale=0.0352778, thin, inner sep=0]
  \def\p#1#2{%
    \ifcase #1
      \or \ifx#2x  0       \else  0 \fi
      \or \ifx#2x  \x1     \else  \y1-50 \fi
      \or \ifx#2x  \x1-42  \else  0 \fi
      \or \ifx#2x  \x2-18  \else  0 \fi
      \or \ifx#2x  \x1+42  \else  0 \fi
      \or \ifx#2x  \x2+18  \else  0 \fi
    \fi}
  \def\x#1{\p#1x}
  \def\y#1{\p#1y}
 \path[use as bounding box] (-42,-75) rectangle (42,25);
 \draw (\x1,\y1+25) -- (\x2,\y2-25);
 \draw (\x3,\y1+25)[rounded corners=3] -- (\x3,\y1-20) -- (\x4,\y2+20) -- (\x4,\y2);
 \draw (\x3,\y2-25)[rounded corners=3,white,line width=4.0] -- (\x3,\y2+20) -- (\x4,\y1-20) -- (\x4,\y1);
 \draw (\x3,\y2-25)[rounded corners=3] -- (\x3,\y2+20) -- (\x4,\y1-20) -- (\x4,\y1);
 \draw (\x5,\y1+25)[rounded corners=3] -- (\x5,\y1-20) -- (\x6,\y2+20) -- (\x6,\y2);
 \draw (\x5,\y2-25)[rounded corners=3,white,line width=4.0] -- (\x5,\y2+20) -- (\x6,\y1-20) -- (\x6,\y1);
 \draw (\x5,\y2-25)[rounded corners=3] -- (\x5,\y2+20) -- (\x6,\y1-20) -- (\x6,\y1);
 \filldraw[fill=white] (\x1-30,\y1-10) rectangle (\x1+30,\y1+10);
 \filldraw[fill=white] (\x2-30,\y2-10) rectangle (\x2+30,\y2+10);
 \node at (\x1,\y1) {$f$};
 \node at (\x2,\y2) {$g$};
\end{tikzpicture}
$$
We have trouble if we want to compose $f$ and $g$ as a step
 of sequentialization, since
 wires of $f$ cross over wires of $g$.
Two crossings are neither that of $\otimes$ nor $\linpar$.
We should redraw the graph as follows:
$$
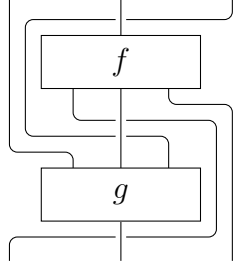
\begin{tikzpicture}[xscale=0.0352778, yscale=0.0352778, thin, inner sep=0]
  \def\p#1#2{%
    \ifcase #1
      \or \ifx#2x  0       \else  0 \fi
      \or \ifx#2x  \x1     \else  \y1-50 \fi
      \or \ifx#2x  \x1-42  \else  0 \fi
      \or \ifx#2x  \x2-18  \else  0 \fi
      \or \ifx#2x  \x1+42  \else  0 \fi
      \or \ifx#2x  \x2+18  \else  0 \fi
      \or \ifx#2x  \x1-36  \else  0 \fi
      \or \ifx#2x  \x1+36  \else  0 \fi
    \fi}
  \def\x#1{\p#1x}
  \def\y#1{\p#1y}
 \path[use as bounding box] (-42,-75) rectangle (42,25);
 \draw (\x3,\y1+25)[rounded corners=3] -- (\x3,\y2+16) -- (\x4,\y2+16) -- (\x4,\y2);
 \draw (\x3,\y2-25)[rounded corners=3] -- (\x3,\y2-16) -- (\x8,\y2-16) -- (\x8,\y1-22)
  -- (\x4,\y1-22) -- (\x4,\y1);
 \draw (\x5,\y1+25)[rounded corners=3] -- (\x5,\y1+16) -- (\x7,\y1+16) -- (\x7,\y2+22)
  -- (\x6,\y2+22) -- (\x6,\y2);
 \draw (\x5,\y2-25)[rounded corners=3] -- (\x5,\y1-16) -- (\x6,\y1-16) -- (\x6,\y1);
 \draw[white,line width=4.0] (\x1,\y1+25) -- (\x2,\y2-25);
 \draw (\x1,\y1+25) -- (\x2,\y2-25);
 \filldraw[fill=white] (\x1-30,\y1-10) rectangle (\x1+30,\y1+10);
 \filldraw[fill=white] (\x2-30,\y2-10) rectangle (\x2+30,\y2+10);
 \node at (\x1,\y1) {$f$};
 \node at (\x2,\y2) {$g$};
\end{tikzpicture}
$$
Now the crossing wires between two boxes are both connected to $f$ or
 both connected to $g$.
The former is interpreted
 as the symmetry of cotensor, and the latter as that of tensor.
The outer crossings are handled after composing $f$ and $g$.

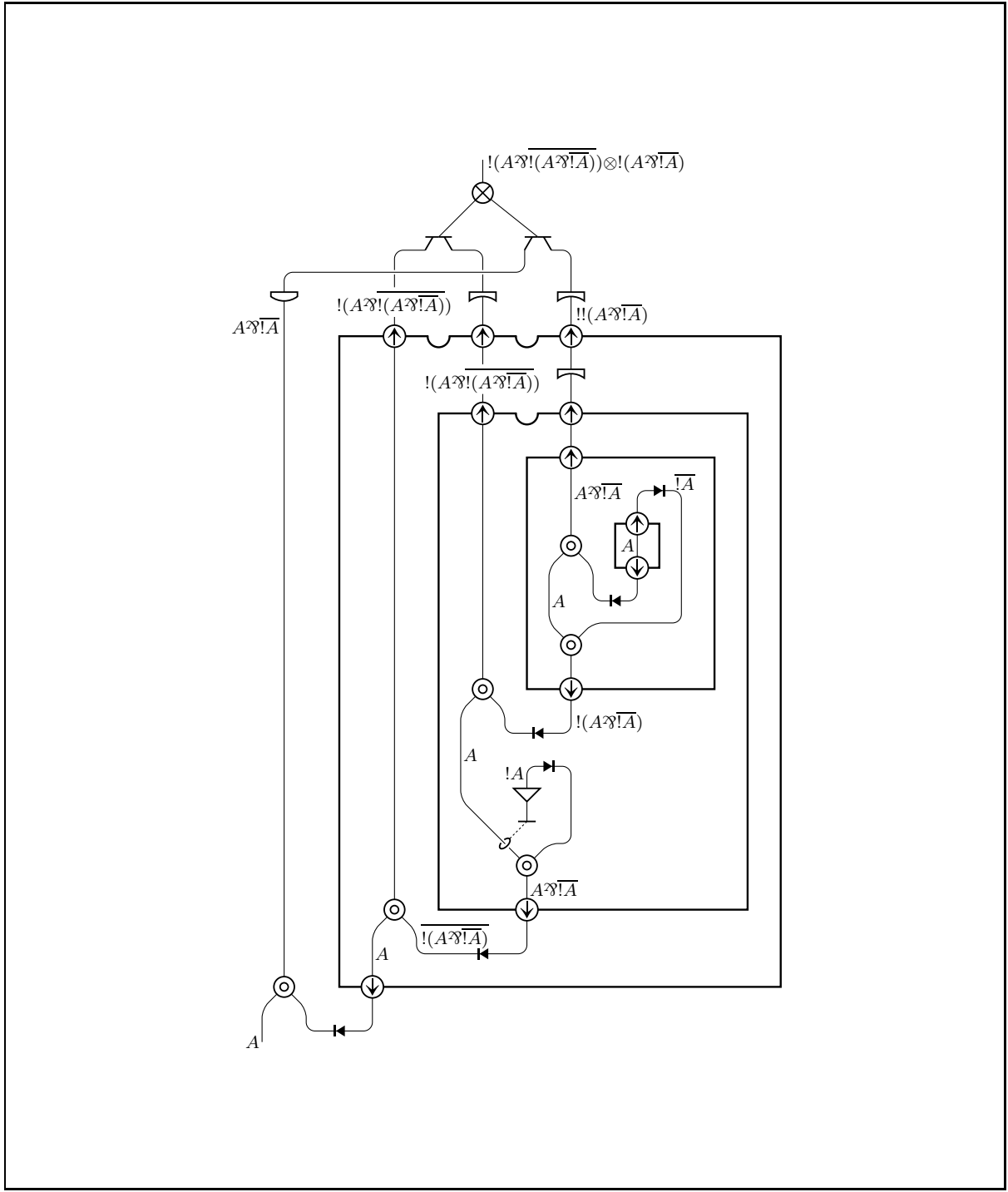
\begin{figure}[t]
\figbox{%
$$
\begin{tikzpicture}[xscale=0.0352778, yscale=0.0352778, thin, inner sep=0]
  \def\pa#1#2{%
    \ifcase #1
      \or \ifx#2x  0           \else  0 \fi
      \or \ifx#2x  \xa1+25     \else  \ya1-20 \fi
      \or \ifx#2x  \xa2+15     \else  \ya2+20 \fi
      \or \ifx#2x  \xa3+10     \else  \ya3+35 \fi
      \or \ifx#2x  \xa4+40     \else  \ya4-20 \fi
      \or \ifx#2x  \xa5+20     \else  \ya5+20 \fi
      \or \ifx#2x  \xa6        \else  \ya6+20 \fi
      \or \ifx#2x  \xa7+10     \else  \ya7+45 \fi
      \or \ifx#2x  \xa8-10     \else  \ya8-10 \fi
    \fi}
  \def\xa#1{\pa#1x}
  \def\ya#1{\pa#1y}
  \def\pb#1#2{%
    \ifcase #1
          \ifx#2x  \xa9        \else  \ya9-15 \fi
      \or \ifx#2x  \xa7-20     \else  \ya7+80 \fi
      \or \ifx#2x  \xb1+25     \else  \yb1-20 \fi
      \or \ifx#2x  \xb2+15     \else  \yb2+20 \fi
      \or \ifx#2x  \xb3        \else  \yb3+20 \fi
      \or \ifx#2x  \xb4+40     \else  \yb4+70 \fi
      \or \ifx#2x  \xb5-10     \else  \yb5-15 \fi
      \or \ifx#2x  \xb6        \else  \yb6-20 \fi
      \or \ifx#2x  \xb7-10     \else  \yb7-15 \fi
      \or \ifx#2x  \xb8-20     \else  \yb8+25 \fi
    \fi}
  \def\xb#1{\pb#1x}
  \def\yb#1{\pb#1y}
  \def\pc#1#2{%
    \ifcase #1
          \ifx#2x  \xb9        \else  \yb9+40 \fi
      \or \ifx#2x  \xc0        \else  \yc0+20 \fi
      \or \ifx#2x  \xc1        \else  \yc1+20 \fi
      \or \ifx#2x  \xc2        \else  \yc2+15 \fi
      \or \ifx#2x  \xc3        \else  \yc3+20 \fi
      \or \ifx#2x  \xc4-15     \else  \yc4+25 \fi
      \or \ifx#2x  \xc5-25     \else  \yc5+20 \fi
      \or \ifx#2x  \xb1        \else  \yc1 \fi
      \or \ifx#2x  \xc7        \else  \yc3 \fi
      \or \ifx#2x  \xc8-20     \else  \yc5 \fi
    \fi}
  \def\xc#1{\pc#1x}
  \def\yc#1{\pc#1y}
  \def\pd#1#2{%
    \ifcase #1
          \ifx#2x  \xa4        \else  \yc3 \fi
      \or \ifx#2x  \xa1        \else  \yc4 \fi
      \or \ifx#2x  \xc7        \else  \yc4 \fi
    \fi}
  \def\xd#1{\pd#1x}
  \def\yd#1{\pd#1y}
%
 \path[use as bounding box] (-25,-85) rectangle (225,445);
 \draw (\xa1,\ya1)[rounded corners=4] -- (\xa1-10,\ya1-10) -- (\xa1-10,\ya1-25);
 \draw (\xa1,\ya1)[rounded corners=4] -- (\xa1+10,\ya1-10) -- (\xa1+10,\ya1-20) -- (\xa2,\ya2);
 \draw (\xa2,\ya2)[rounded corners=4] -- (\xa3,\ya3-20) -- (\xa3,\ya3);
 \draw (\xa3,\ya3)[rounded corners=4] -- (\xa4-10,\ya4-10) -- (\xa4,\ya4);
 \draw (\xa4,\ya4)[rounded corners=4] -- (\xa4+10,\ya4-10) -- (\xa4+10,\ya5) -- (\xa5,\ya5);
 \draw (\xa5,\ya5)[rounded corners=4] -- (\xa6,\ya5) -- (\xa6,\ya6);
 \draw (\xa6,\ya6) -- (\xa7,\ya7);
 \draw (\xa7,\ya7)[rounded corners=4] -- (\xa7+10,\ya7+10) -- (\xa8+10,\ya7+10)
   -- (\xa8+10,\ya8) -- (\xa8,\ya8);
 \draw (\xa8,\ya8)[rounded corners=4] -- (\xa9,\ya8) -- (\xa9,\ya9);
 \draw (\xa9,\ya9) -- (\xb0,\yb0);
 \draw (\xa7,\ya7)[rounded corners=4] -- (\xa7-30,\ya7+30) -- (\xb1-10,\yb1-10)
   -- (\xb1,\yb1);
 \draw (\xb1,\yb1)[rounded corners=4] -- (\xb1+10,\yb1-10) -- (\xb1+10,\yb2) -- (\xb2,\yb2);
 \draw (\xb2,\yb2)[rounded corners=4] -- (\xb3,\yb2) -- (\xb3,\yb3);
 \draw (\xb3,\yb3) -- (\xb4,\yb4);
 \draw (\xb4,\yb4)[rounded corners=4] -- (\xb4+10,\yb4+10) -- (\xb5+10,\yb4+10)
   -- (\xb5+10,\yb5) -- (\xb5,\yb5);
 \draw (\xb5,\yb5)[rounded corners=4] -- (\xb6,\yb5) -- (\xb6,\yb6);
 \draw (\xb6,\yb6) -- (\xb7,\yb7);
 \draw (\xb7,\yb7)[rounded corners=4] -- (\xb7,\yb8) -- (\xb8,\yb8);
 \draw (\xb8,\yb8)[rounded corners=4] -- (\xb9+10,\yb8) -- (\xb9+10,\yb9-10) -- (\xb9,\yb9);
 \draw (\xb4,\yb4)[rounded corners=4] -- (\xb4-10,\yb4+10) -- (\xb9-10,\yb9-10) -- (\xb9,\yb9);
 \draw (\xb9,\yb9) -- (\xc0,\yc0);
 \draw (\xc0,\yc0) -- (\xc2,\yc2);
 \draw (\xc2,\yc2) -- (\xc3,\yc3);
 \draw (\xc3,\yc3) -- (\xc4,\yc4);
 \draw (\xc4,\yc4)[rounded corners=4] -- (\xc4,\yc5-6) -- (\xc5+6,\yc5-6);
 \draw (\xc5,\yc5) -- (\xc6,\yc6);
 \draw (\xb1,\yb1) -- (\xc7,\yc7);
 \draw (\xc7,\yc7) -- (\xd2,\yd2);
 \draw (\xd2,\yd2)[rounded corners=4] -- (\xc8,\yc9-6) -- (\xc9+6,\yc9-6);
 \draw (\xc6,\yc6) -- (\xc9,\yc9);
 \draw (\xa4,\ya4) -- (\xd0,\yd0);
 \draw (\xd0,\yd0)[rounded corners=4] -- (\xd0,\yc9-6) -- (\xc9-6,\yc9-6);
 \draw (\xa1,\ya1) -- (\xd1,\yd1);
 \draw[white,line width=2.0] (\xd1,\yc5-16) -- (\xc5-6,\yc5-16);
 \draw (\xd1,\yd1)[rounded corners=4] -- (\xd1,\yc5-16) -- (\xc5-6,\yc5-16)
   -- (\xc5-6,\yc5-6);
 \draw (\xc6,\yc6) -- (\xc6,\yc6+15);
 \draw[dash pattern=on 1 off 1] (\xb0,\yb0) -- (\xa7-8,\ya7+12);
 \begin{scope}[shift={(\xa7-10,\ya7+10)},rotate=45]
  \putRubberBand(0,0)
 \end{scope}
 \boardBoundary(\xc3+95,\yc3)(\xa3-15,\ya3)
 \boardBoundary(\xc7-20,\yc7)(\xa6+100,\ya6)
 \boardBoundary(\xb3-20,\yb3)(\xc0+65,\yc0)
 \boardBoundary(\xb6-10,\yb6)(\xb7+10,\yb7)
 \putPar(\xa1,\ya1)
 \putLeftDiode(\xa2,\ya2)
 \putLowerSocket(\xa3,\ya3)
 \putPar(\xa4,\ya4)
 \putLeftDiode(\xa5,\ya5)
 \putLowerSocket(\xa6,\ya6)
 \putPar(\xa7,\ya7)
 \putRightDiode(\xa8,\ya8)
 \putEliminator(\xa9,\ya9)
 \putPositiveTerminal(\xb0,\yb0)
 \putPar(\xb1,\yb1)
 \putLeftDiode(\xb2,\yb2)
 \putLowerSocket(\xb3,\yb3)
 \putPar(\xb4,\yb4)
 \putRightDiode(\xb5,\yb5)
 \putUpperSocket(\xb6,\yb6)
 \putLowerSocket(\xb7,\yb7)
 \putLeftDiode(\xb8,\yb8)
 \putPar(\xb9,\yb9)
 \putUpperSocket(\xc0,\yc0)
 \putUpperSocket(\xc1,\yc1)
 \putConcave(\xc2,\yc2)
 \putUpperSocket(\xc3,\yc3)
 \putConcave(\xc4,\yc4)
 \putDuplicator(\xc5,\yc5)
 \putTensor(\xc6,\yc6)
 \putUpperSocket(\xc7,\yc7)
 \putUpperSocket(\xc8,\yc8)
 \putDuplicator(\xc9,\yc9)
 \putUpperSocket(\xd0,\yd0)
 \putConvex(\xd1,\yd1)
 \putConcave(\xd2,\yd2)
 \putDimple(\xc1-20,\yc1)
 \putDimple(\xc3-20,\yc3)
 \putDimple(\xc8-20,\yc8)
 \node[left] at (\xa1-11,\ya1-25) {$\scriptstyle A$};
 \node[right] at (\xa3+1,\ya3+15) {$\scriptstyle A$};
 \node[right] at (\xa4+12,\ya4-12) {$\scriptstyle
   \overline{\mathord!(A\slinpar\overline{\mathord!A})}$};
 \node[right] at (\xa6+1,\ya6+10) {$\scriptstyle
   A\slinpar\overline{\mathord!A}$};
 \node[left] at (\xa9-1,\ya9+7) {$\scriptstyle \mathord!A$};
 \node[right] at (\xb1-9,\yb1-30) {$\scriptstyle A$};
 \node[right] at (\xb3+2,\yb3-15) {$\scriptstyle
   \mathord!(A\slinpar\overline{\mathord!A})$};
 \node[right] at (\xb4-9,\yb4+20) {$\scriptstyle A$};
 \node[left] at (\xb7-1,\yb7+10) {$\scriptstyle A$};
 \node[right] at (\xb9+1,\yb9+25) {$\scriptstyle
   A\slinpar\overline{\mathord!A}$};
 \node[right] at (\xb5+7,\yb5+4) {$\scriptstyle
   \overline{\mathord!A}$};
 \filldraw[white] (\xc7-1,\yc7+9) rectangle (\xc7+1,\yc7+22);
 \node at (\xc7,\yc7+15) {$\scriptstyle
   \mathord!(A\slinpar
   \overline{\mathord!(A\slinpar\overline{\mathord!A})})$};
 \node[right] at (\xc3+2,\yc3+10) {$\scriptstyle
   \mathord!\mathord!(A\slinpar\overline{\mathord!A})$};
 \filldraw[white] (\xd0-1,\yd0+9) rectangle (\xd0+1,\yd0+22);
 \node at (\xd0,\yd0+15) {$\scriptstyle
   \mathord!(A\slinpar
   \overline{\mathord!(A\slinpar\overline{\mathord!A})})$};
 \node[left] at (\xd1-2,\yd1-15) {$\scriptstyle
   A\slinpar\overline{\mathord!A}$};
 \node[right] at (\xc6+2,\yc6+15) {$\scriptstyle
   \mathord!(A\slinpar
   \overline{\mathord!(A\slinpar\overline{\mathord!A})})\otimes
   \mathord!(A\slinpar\overline{\mathord!A})$};
\end{tikzpicture}
$$
}%
\caption{An example of a graph}\label{yes66}
\end{figure}

Figure~\ref{yes66} gives an example of a graph.
\afterpage{\clearpage}%
It is the translation of the following morphisms.
To facilitate the parsing, we denote the dual by $\overline A$
 in place of $A^*$.
Let $f$ be\footnote{%
We use dots in place of appropriate identity morphisms in order to save space
 and promote readability.}%

	\vskip1ex
        \noindent\kern3em
\begingroup\small
\def\linpar{\mathbin{\tikz[inner sep=0] \node[rotate=180] at (0,0) {\small $\&$};}}%
\vbox{\halign{$#$\hfil\cr
  \mathord!(A\linpar
    \overline{\mathord!(A\linpar\overline{\mathord!A})})\otimes
    \mathord!\mathord!(A\linpar\overline{\mathord!A})
  \mor{\tilde\varphi}
  \mathord!((A\linpar \overline{\mathord!(A\linpar\overline{\mathord!A})})\otimes
    \mathord!(A\linpar\overline{\mathord!A}))
  \mor{\mathord!\partial'}\cr
  \qquad\mathord!(A\linpar (\overline{\mathord!(A\linpar\overline{\mathord!A})}\otimes
    \mathord!(A\linpar\overline{\mathord!A})))
  \mor{\mathord!(\cdot\gamma)}
  \mathord!(A\linpar\bot)
  \mor{\sim}\mathord!A\mor{\sim}\mathord!(A\otimes{\bf 1})
  \mor{\mathord!(\cdot\tau)}\cr
  \qquad\mathord!(A\otimes(\mathord!A\linpar\overline{\mathord!A}))
  \mor{\mathord!(\cdot(e\cdot))}
  \mathord!(A\otimes({\bf 1}\linpar\overline{\mathord!A}))
  \mor{\mathord!\partial}
  \mathord!((A\otimes{\bf 1})\linpar\overline{\mathord!A})
  \mor{\sim}
  \mathord!(A\linpar\overline{\mathord!A})\cr
}}
\endgroup

	\vskip1ex
        \noindent
 and let $g$ be the following, which uses $f$:
	\vskip1ex
        \noindent\kern3em
\begingroup\small
\def\linpar{\mathbin{\tikz[inner sep=0] \node[rotate=180] at (0,0) {\small $\&$};}}%
\vbox{\halign{$#$\hfil\cr
  \mathord!(A\linpar
    \overline{\mathord!(A\linpar\overline{\mathord!A})})\otimes
    \mathord!\mathord!(A\linpar
    \overline{\mathord!(A\linpar\overline{\mathord!A})})\otimes
    \mathord!\mathord!(A\linpar\overline{\mathord!A})
  \mor{\tilde\varphi^2}\cr
  \qquad\mathord!((A\linpar
    \overline{\mathord!(A\linpar\overline{\mathord!A})})\otimes
    \mathord!(A\linpar
    \overline{\mathord!(A\linpar\overline{\mathord!A})})\otimes
    \mathord!(A\linpar\overline{\mathord!A}))
  \mor{\mathord!(\cdot\cdot\delta)}\cr
  \qquad\mathord!((A\linpar
    \overline{\mathord!(A\linpar\overline{\mathord!A})})\otimes
    \mathord!(A\linpar
    \overline{\mathord!(A\linpar\overline{\mathord!A})})\otimes
    \mathord!\mathord!(A\linpar\overline{\mathord!A}))
  \mor{\mathord!(\cdot f)}\cr
  \qquad\mathord!((A\linpar
    \overline{\mathord!(A\linpar\overline{\mathord!A})})\otimes
    \mathord!(A\linpar\overline{\mathord!A}))
  \mor{\mathord!\partial'}
  \mathord!(A\linpar
    (\overline{\mathord!(A\linpar\overline{\mathord!A})}\otimes
    \mathord!(A\linpar\overline{\mathord!A})))
  \mor{\mathord!(\cdot\gamma)}\cr
  \qquad\mathord!(A\linpar\bot)\mor\sim\mathord!A.\cr
}}
\endgroup

	\vskip1ex
        \noindent
The figure is the graph of the following morphism containing
 $g$:
	\vskip1ex
        \noindent\kern3em
\begingroup\small
\def\linpar{\mathbin{\tikz[inner sep=0] \node[rotate=180] at (0,0) {\small $\&$};}}%
\vbox{\halign{$#$\hfil\cr
  \mathord!(A\linpar
    \overline{\mathord!(A\linpar\overline{\mathord!A})})\otimes
    \mathord!(A\linpar\overline{\mathord!A})
  \mor{dd}\cr
  \qquad\mathord!(A\linpar
    \overline{\mathord!(A\linpar\overline{\mathord!A})})\otimes
    \mathord!(A\linpar
    \overline{\mathord!(A\linpar\overline{\mathord!A})})\otimes
    \mathord!(A\linpar\overline{\mathord!A})\otimes
    \mathord!(A\linpar\overline{\mathord!A})
  \mor\sim\cr
  \qquad\mathord!(A\linpar\overline{\mathord!A})\otimes
    \mathord!(A\linpar
    \overline{\mathord!(A\linpar\overline{\mathord!A})})\otimes
    \mathord!(A\linpar
    \overline{\mathord!(A\linpar\overline{\mathord!A})})\otimes
    \mathord!(A\linpar\overline{\mathord!A})
  \mor{\varepsilon\cdot\delta\delta}\cr
  \qquad(A\linpar\overline{\mathord!A})\otimes
    \mathord!(A\linpar
    \overline{\mathord!(A\linpar\overline{\mathord!A})})\otimes
    \mathord!\mathord!(A\linpar
    \overline{\mathord!(A\linpar\overline{\mathord!A})})\otimes
        \mathord!\mathord!(A\linpar\overline{\mathord!A})
  \mor{\cdot g}\cr
  \qquad(A\linpar\overline{\mathord!A})\otimes\mathord!A
  \mor{\partial'}
  A\linpar(\overline{\mathord!A}\otimes\mathord!A)
  \mor{\cdot\gamma}A\linpar\bot\mor\sim A\cr
}}
\endgroup

	\vskip1ex
        \noindent
The morphism $f$ corresponds to the second board from the outermost,
 while $g$ to the outermost.
The inner two boards are created by $\eta$-expansions.
The types of wires are omitted in part by space restriction.
Although the morphisms are long, they are derived from
 the typing judgment $\Gamma\vdash y(x(\lambda z.\,xy))$ of the
 lambda calculus.
Here $\Gamma$ consists of $x\mathbin:(A\Rightarrow A)\Rightarrow A$
 and $y:A\Rightarrow A$.

\begin{definition}\label{zef08}\rm
The graphs $G$ and $G'$ are {\it almost equal} if they are equal
 when all the dotted lines and the attached rings of units/counits
 are erased.
We write $G\sim G'$.
Accordingly, if $f$ and $f'$ are morphisms in normal forms and
 their associated graphs are almost equal, we write $f\sim f'$.
\end{definition}

	\vskip1ex

Namely, $f\sim f'$ holds if they are equal except how the
 isomorphisms $A\otimes{\bf 1}\cong A$ and $A\linpar \bot \cong A$
 are used.
The goal of this paper is to show $f\sim f'$ whenever
 $f$ and $f'$ are the normal forms obtained by contracting
 a common morphism.

\section{Enumeration}\label{yhv89}

We return to the analysis of the model $\mathscr{M}$.
Hereafter, we often omit the brackets in $[\![A]\!]$,
 simply writing $A$.
This is not harmful since we are not concerned with
 equalities between objects.
In contrast,
 brackets enclosing morphisms should not be omitted,
 since we are motivated to relate the equality of morphisms
 of the free linear category with that of
 the model $\mathscr{M}$.

In general, it is difficult to compute
 the matrix component of $M_f$.
In fact, $M_f$ may be involved in
 the matrix multiplication $(NM)[a,c]=
 \sum_b N[b,c]M[a,b]$ in which
 $b$ ranges over
 an infinite set, thus not computable in finite time.
However, we have an effective method to compute the
 matrix components.
We first contract $f$ to a normal form,
 then translate it to a graph, whence the following
 process works.

The effective method is given by an enumeration process on the graph.
Suppose that the graph
$$
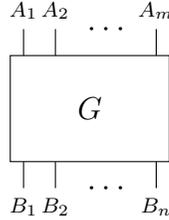
\begin{tikzpicture}[xscale=0.0352778, yscale=0.0352778, thin, inner sep=0]
  \def\p#1#2{%
    \ifcase #1
      \or \ifx#2x  5       \else  0 \fi
      \or \ifx#2x  \x1+12     \else  \y1 \fi
      \or \ifx#2x  \x1+50     \else  \y1 \fi
    \fi}
  \def\x#1{\p#1x}
  \def\y#1{\p#1y}
 \path[use as bounding box] (0,-20) rectangle (60,60);
 \draw (\x1,\y1-10) -- (\x1,\y1+50);
 \draw (\x2,\y2-10) -- (\x2,\y2+50);
 \draw (\x3,\y3-10) -- (\x3,\y3+50);
 \filldraw[fill=white] (0,0) rectangle (60,40);
 \node at (30,20) {$G$};
 \node at ($ .5*(\x3,\y3-10)+.5*(\x2,\y2-10) $) {$\cdots$};
 \node at ($ .5*(\x3,\y3+50)+.5*(\x2,\y2+50) $) {$\cdots$};
 \node at (\x1,\y1+57) {$\scriptstyle A_1$};
 \node at (\x2,\y2+57) {$\scriptstyle A_2$};
 \node at (\x3,\y3+57) {$\scriptstyle A_m$};
 \node at (\x1,\y1-17) {$\scriptstyle B_1$};
 \node at (\x2,\y2-17) {$\scriptstyle B_2$};
 \node at (\x3,\y3-17) {$\scriptstyle B_n$};
\end{tikzpicture}
$$
 is obtained from a normal form $f:
 A_1\otimes A_2\otimes\cdots\otimes A_m
 \rightarrow B_1\linpar B_2\linpar\cdots\linpar B_n$.
We want to compute the component
 $M_f[(\alpha_1,\alpha_2,\ldots,\alpha_m);(\beta_1,\beta_2,
 \ldots,\beta_n)]$ for $\alpha_i\in A_i$
 and $\beta_j\in B_j$.
We define a process
 $\pi(G)^{\alpha_1\alpha_2\cdots\alpha_m}_{\beta_1\beta_2\cdots\beta_n}$.
The return value of the process is a non-negative integer.
We have
 $M_f[(\alpha_1,\alpha_2,\ldots,\alpha_m);(\beta_1,\beta_2,
 \ldots,\beta_n)]
=\pi(G)^{\alpha_1\alpha_2\cdots\alpha_m}_{\beta_1\beta_2\cdots\beta_n}$.
The process may fork meantime, spawning several
 subprocesses to run concurrently.
Moreover, speculative executions are carried out, thus
 some processes may fail at some point and are aborted.
The aborted processes are enforced to return $0$.
In other words, they do not contribute to the enumeration.
Successful processes return positive integers.

We start with annotating the outermost wires by
 the elements as
$$
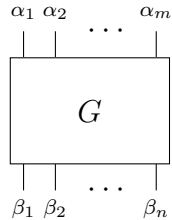
\begin{tikzpicture}[xscale=0.0352778, yscale=0.0352778, thin, inner sep=0]
  \def\p#1#2{%
    \ifcase #1
      \or \ifx#2x  5       \else  0 \fi
      \or \ifx#2x  \x1+12     \else  \y1 \fi
      \or \ifx#2x  \x1+50     \else  \y1 \fi
    \fi}
  \def\x#1{\p#1x}
  \def\y#1{\p#1y}
 \path[use as bounding box] (0,-20) rectangle (60,60);
 \draw (\x1,\y1-10) -- (\x1,\y1+50);
 \draw (\x2,\y2-10) -- (\x2,\y2+50);
 \draw (\x3,\y3-10) -- (\x3,\y3+50);
 \filldraw[fill=white] (0,0) rectangle (60,40);
 \node at (30,20) {$G$};
 \node at ($ .5*(\x3,\y3-10)+.5*(\x2,\y2-10) $) {$\cdots$};
 \node at ($ .5*(\x3,\y3+50)+.5*(\x2,\y2+50) $) {$\cdots$};
 \node at (\x1,\y1+57) {$\scriptstyle \alpha_1$};
 \node at (\x2,\y2+57) {$\scriptstyle \alpha_2$};
 \node at (\x3,\y3+57) {$\scriptstyle \alpha_m$};
 \node at (\x1,\y1-17) {$\scriptstyle \beta_1$};
 \node at (\x2,\y2-17) {$\scriptstyle \beta_2$};
 \node at (\x3,\y3-17) {$\scriptstyle \beta_n$};
\end{tikzpicture}
$$
The process $\pi(G)=\pi(G)^{\alpha_1\alpha_2\cdots\alpha_m}_{\beta_1\beta_2\cdots\beta_n}$
 decomposes the graph in the reverse order of sequentialization.
The decomposition is occasionally non-deterministic, though
 it is deterministic in many cases.
If the annotated graph $G$ is the left of the following
$$
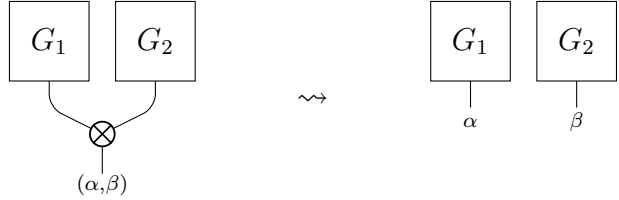
\vcenter{\hbox{%
\begin{tikzpicture}[xscale=0.0352778, yscale=0.0352778, thin, inner sep=0]
  \def\p#1#2{%
    \ifcase #1
      \or \ifx#2x  0       \else  0 \fi
      \or \ifx#2x  \x1+40     \else  \y1 \fi
      \or \ifx#2x  \x1+35     \else  \y1-20 \fi
    \fi}
  \def\x#1{\p#1x}
  \def\y#1{\p#1y}
 \path[use as bounding box] (0,-43) rectangle (70,30);
 \draw (\x1,\y1) rectangle (\x1+30,\y1+30);
 \draw (\x2,\y2) rectangle (\x2+30,\y2+30);
 \draw (\x1+15,\y1)[rounded corners=5] -- (\x1+15,\y1-10) -- (\x3,\y3);
 \draw (\x2+15,\y2)[rounded corners=5] -- (\x2+15,\y2-10) -- (\x3,\y3);
 \draw (\x3,\y3) -- (\x3, \y3-15);
 \node at (\x1+15,\y1+15) {$G_1$};
 \node at (\x2+15,\y2+15) {$G_2$};
 \node[below] at (\x3,\y3-15) {$\scriptstyle(\alpha,\beta)$};
 \putTensor(\x3,\y3)
\end{tikzpicture}%
}}%
\kern3em\rightsquigarrow\kern3em
\vcenter{\hbox{%
\begin{tikzpicture}[xscale=0.0352778, yscale=0.0352778, thin, inner sep=0]
  \def\p#1#2{%
    \ifcase #1
      \or \ifx#2x  0       \else  0 \fi
      \or \ifx#2x  \x1+40     \else  \y1 \fi
    \fi}
  \def\x#1{\p#1x}
  \def\y#1{\p#1y}
 \path[use as bounding box] (0,-43) rectangle (70,30);
 \draw (\x1,\y1) rectangle (\x1+30,\y1+30);
 \draw (\x2,\y2) rectangle (\x2+30,\y2+30);
 \draw (\x1+15,\y1) -- (\x1+15,\y1-10);
 \draw (\x2+15,\y2) -- (\x2+15,\y2-10);
 \node at (\x1+15,\y1+15) {$G_1$};
 \node at (\x2+15,\y2+15) {$G_2$};
 \node[below,inner sep=2,text height=5] at (\x1+15,\y1-10) {$\scriptstyle \alpha$};
 \node[below,inner sep=2,text height=5] at (\x2+15,\y2-10) {$\scriptstyle \beta$};
\end{tikzpicture}%
}}%
$$
 then it is deterministically decomposed into
 the right annotated graph.
Here the graphs may have more outgoing wires, which are suppressed
 for simplicity.
The same comment is applied to all of the following cases, and
 is repeated no more.
Given the return values of the processes $\pi(G_1)_\alpha$ and
 $\pi(G_2)_\beta$, we set $\pi(G)_{(\alpha,\beta)}=
\pi(G_1)_\alpha \cdot \pi(G_2)_\beta$.
If the graph $G$ is the left of the following,
$$\vcenter
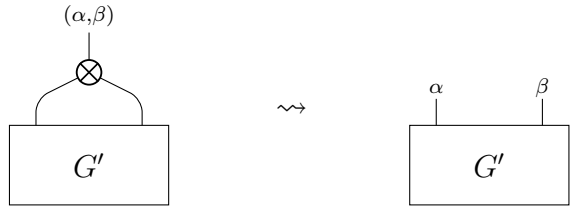
{\hbox{%
\begin{tikzpicture}[xscale=0.0352778, yscale=0.0352778, thin, inner sep=0]
  \def\p#1#2{%
    \ifcase #1
      \or \ifx#2x  0       \else  0 \fi
      \or \ifx#2x  \x1+10     \else  \y1+30 \fi
      \or \ifx#2x  \x1+50     \else  \y1+30 \fi
      \or \ifx#2x  \x1+60     \else  \y1+30 \fi
      \or \ifx#2x  \x1+30     \else  \y2+20 \fi
    \fi}
  \def\x#1{\p#1x}
  \def\y#1{\p#1y}
 \path[use as bounding box] (0,0) rectangle (62,73);
 \draw (\x1,\y1) rectangle (\x4,\y4);
 \draw (\x2,\y2)[rounded corners=5] -- (\x2,\y2+10) -- (\x5,\y5);
 \draw (\x3,\y3)[rounded corners=5] -- (\x3,\y3+10) -- (\x5,\y5);
 \draw (\x5,\y5) -- (\x5,\y5+15);
 \node[above,inner sep=2] at (\x5,\y5+15) {$\scriptstyle (\alpha,\beta)$};
 \node at (\x5,\y1+15) {$G'$};
 \putTensor(\x5,\y5)
\end{tikzpicture}%
}}%
\kern3em\rightsquigarrow\kern3em
\vcenter{\hbox{%
\begin{tikzpicture}[xscale=0.0352778, yscale=0.0352778, thin, inner sep=0]
  \def\p#1#2{%
    \ifcase #1
      \or \ifx#2x  0       \else  0 \fi
      \or \ifx#2x  \x1+10     \else  \y1+30 \fi
      \or \ifx#2x  \x1+50     \else  \y1+30 \fi
      \or \ifx#2x  \x1+60     \else  \y1+30 \fi
      \or \ifx#2x  \x1+30     \else  \y2+20 \fi
    \fi}
  \def\x#1{\p#1x}
  \def\y#1{\p#1y}
 \path[use as bounding box] (0,0) rectangle (62,73);
 \draw (\x1,\y1) rectangle (\x4,\y4);
 \draw (\x2,\y2)[rounded corners=5] -- (\x2,\y2+10);
 \draw (\x3,\y3)[rounded corners=5] -- (\x3,\y3+10);
 \node[above,text depth=2] at (\x2,\y2+10) {$\scriptstyle \alpha$};
 \node[above,text depth=2] at (\x3,\y3+10) {$\scriptstyle \beta$};
 \node at (\x5,\y1+15) {$G'$};
\end{tikzpicture}%
}}$$
 it is decomposed into the right, and $\pi(G)^{(\alpha,\beta)}=
\pi(G')^{\alpha\beta}$.
In the case of cotensor, the decomposition is similar
 except the graphs are turned upside down.
The unit introduction is manipulated as follows:
$$
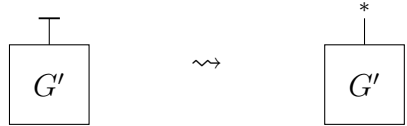
\vcenter{\hbox{%
\begin{tikzpicture}[xscale=0.0352778, yscale=0.0352778, thin, inner sep=0]
  \def\p#1#2{%
    \ifcase #1
      \or \ifx#2x  0       \else  0 \fi
    \fi}
  \def\x#1{\p#1x}
  \def\y#1{\p#1y}
 \path[use as bounding box] (0,0) rectangle (30,46);
 \draw (\x1,\y1) rectangle (\x1+30,\y1+30);
 \draw (\x1+15,\y1+30) -- (\x1+15,\y1+40);
 \node at (\x1+15,\y1+15) {$G'$};
 \putPositiveTerminal(\x1+15,\y1+40)
\end{tikzpicture}%
}}%
\kern3em\rightsquigarrow\kern3em
\vcenter{\hbox{%
\begin{tikzpicture}[xscale=0.0352778, yscale=0.0352778, thin, inner sep=0]
  \def\p#1#2{%
    \ifcase #1
      \or \ifx#2x  0       \else  0 \fi
    \fi}
  \def\x#1{\p#1x}
  \def\y#1{\p#1y}
 \path[use as bounding box] (0,0) rectangle (30,46);
 \draw (\x1,\y1) rectangle (\x1+30,\y1+30);
 \draw (\x1+15,\y1+30) -- (\x1+15,\y1+40);
 \node at (\x1+15,\y1+15) {$G'$};
 \node[above,inner sep=2] at (\x1+15,\y1+40) {$\scriptstyle*$};
\end{tikzpicture}%
}}%
$$
We set $\pi(G)=\pi(G')^*$.
The unit elimination rule is manipulated as follows:
$$
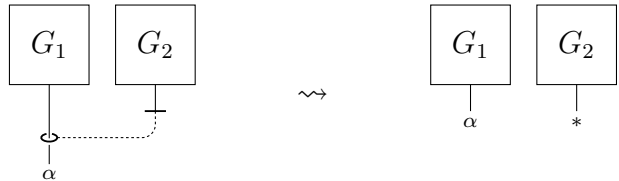
\vcenter{\hbox{%
\begin{tikzpicture}[xscale=0.0352778, yscale=0.0352778, thin, inner sep=0]
  \def\p#1#2{%
    \ifcase #1
      \or \ifx#2x  0       \else  0 \fi
      \or \ifx#2x  \x1+40     \else  \y1 \fi
    \fi}
  \def\x#1{\p#1x}
  \def\y#1{\p#1y}
 \path[use as bounding box] (0,-37) rectangle (70,30);
 \draw (\x1,\y1) rectangle (\x1+30,\y1+30);
 \draw (\x2,\y2) rectangle (\x2+30,\y2+30);
 \draw (\x1+15,\y1) -- (\x1+15,\y1-30);
 \draw (\x2+15,\y2) -- (\x2+15,\y2-10);
 \node at (\x1+15,\y1+15) {$G_1$};
 \node at (\x2+15,\y2+15) {$G_2$};
 \putPositiveTerminal(\x2+15,\y2-10)
 \draw[dash pattern=on 1 off 1] (\x2+15,\y2-10)[rounded corners=5]
  -- (\x2+15,\y2-20) -- (\x1+18,\y2-20);
 \putRubberBand(\x1+15,\y2-20)
 \node[below,inner sep=2] at (\x1+15,\y1-30) {$\scriptstyle\alpha$};
\end{tikzpicture}%
}}%
\kern3em\rightsquigarrow\kern3em
\vcenter{\hbox{%
\begin{tikzpicture}[xscale=0.0352778, yscale=0.0352778, thin, inner sep=0]
  \def\p#1#2{%
    \ifcase #1
      \or \ifx#2x  0       \else  0 \fi
      \or \ifx#2x  \x1+40     \else  \y1 \fi
    \fi}
  \def\x#1{\p#1x}
  \def\y#1{\p#1y}
 \path[use as bounding box] (0,-37) rectangle (70,30);
 \draw (\x1,\y1) rectangle (\x1+30,\y1+30);
 \draw (\x2,\y2) rectangle (\x2+30,\y2+30);
 \draw (\x1+15,\y1) -- (\x1+15,\y1-10);
 \draw (\x2+15,\y2) -- (\x2+15,\y2-10);
 \node at (\x1+15,\y1+15) {$G_1$};
 \node at (\x2+15,\y2+15) {$G_2$};
 \node[below,text height=6] at (\x1+15,\y1-10) {$\scriptstyle\alpha$};
 \node[below,text height=6] at (\x2+15,\y2-10) {$\scriptstyle*$};
\end{tikzpicture}%
}}%
$$
We set $\pi(G)_\alpha=\pi(G_1)_\alpha\cdot\pi(G_2)_*$.
The counit rules are symmetric.
The duality introduction rule is handled as follows:
$$
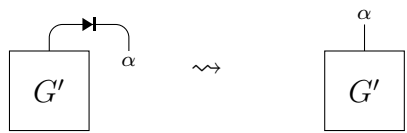
\vcenter{\hbox{%
\begin{tikzpicture}[xscale=0.0352778, yscale=0.0352778, thin, inner sep=0]
  \def\p#1#2{%
    \ifcase #1
      \or \ifx#2x  0       \else  0 \fi
      \or \ifx#2x  \x1+30  \else  \y1+40 \fi
    \fi}
  \def\x#1{\p#1x}
  \def\y#1{\p#1y}
 \path[use as bounding box] (0,0) rectangle (30,47);
 \draw (\x1,\y1) rectangle (\x1+30,\y1+30);
 \draw (\x1+15,\y1+30)[rounded corners=5] -- (\x1+15,\y2)
   -- (\x1+45,\y2) -- (\x1+45,\y1+30);
 \node at (\x1+15,\y1+15) {$G'$};
 \putRightDiode(\x2,\y2)
 \node[below,inner sep=2] at (\x1+45,\y1+30) {$\scriptstyle \alpha$};
\end{tikzpicture}%
}}%
\kern3em\rightsquigarrow\kern3em
\vcenter{\hbox{%
\begin{tikzpicture}[xscale=0.0352778, yscale=0.0352778, thin, inner sep=0]
  \def\p#1#2{%
    \ifcase #1
      \or \ifx#2x  0       \else  0 \fi
      \or \ifx#2x  \x1+30  \else  \y1+40 \fi
    \fi}
  \def\x#1{\p#1x}
  \def\y#1{\p#1y}
 \path[use as bounding box] (0,0) rectangle (30,47);
 \draw (\x1,\y1) rectangle (\x1+30,\y1+30);
 \draw (\x1+15,\y1+30) -- (\x1+15,\y2);
 \node at (\x1+15,\y1+15) {$G'$};
 \node[above,inner sep=2] at (\x1+15,\y2) {$\scriptstyle \alpha$};
\end{tikzpicture}%
}}%
$$
We set $\pi(G)_\alpha=\pi(G')^\alpha$.
The elimination rule is symmetric.

Next, we consider the parts related to $\mathord!$.
Here the speculative executions matter.
In the case of an $e$-part, if $G$ is the
 left-hand side of
$$
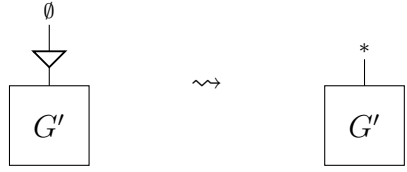
\vcenter{\hbox{%
\begin{tikzpicture}[xscale=0.0352778, yscale=0.0352778, thin, inner sep=0]
  \def\p#1#2{%
    \ifcase #1
      \or \ifx#2x  0       \else  0 \fi
    \fi}
  \def\x#1{\p#1x}
  \def\y#1{\p#1y}
 \path[use as bounding box] (-15,-30) rectangle (15,32);
 \draw (\x1-15,\y1) rectangle (\x1+15,\y1-30);
 \draw (\x1,\y1) -- (\x1,\y1+23);
 \node at (\x1,\y1-15) {$G'$};
 \putEliminator(\x1,\y1+13)
 \node[above,inner sep=2] at (\x1,\y1+23) {$\scriptstyle\emptyset$};
\end{tikzpicture}%
}}%
\kern3em\rightsquigarrow\kern3em
\vcenter{\hbox{%
\begin{tikzpicture}[xscale=0.0352778, yscale=0.0352778, thin, inner sep=0]
  \def\p#1#2{%
    \ifcase #1
      \or \ifx#2x  0       \else  0 \fi
    \fi}
  \def\x#1{\p#1x}
  \def\y#1{\p#1y}
 \path[use as bounding box] (-15,-30) rectangle (15,32);
 \draw (\x1-15,\y1) rectangle (\x1+15,\y1-30);
 \draw (\x1,\y1) -- (\x1,\y1+10);
 \node at (\x1,\y1-15) {$G'$};
 \node[above,inner sep=2] at (\x1,\y1+10) {$\scriptstyle*$};
\end{tikzpicture}%
}}%
$$
 then it is modified to the right graph.
We set $\pi(G)^{\emptyset}=\pi(G')^*$.
If the wire in the left
 graph is annotated by $\alpha\neq\emptyset$,
 the process fails and is aborted immediately.
We impose $\pi(G)^{\alpha}=0$, not proceeding to
 $G'$.
The case of an $\varepsilon$-part is handled as follows:
$$
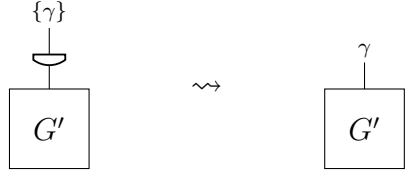
\vcenter{\hbox{%
\begin{tikzpicture}[xscale=0.0352778, yscale=0.0352778, thin, inner sep=0]
  \def\p#1#2{%
    \ifcase #1
      \or \ifx#2x  0       \else  0 \fi
    \fi}
  \def\x#1{\p#1x}
  \def\y#1{\p#1y}
 \path[use as bounding box] (-15,-30) rectangle (15,32);
 \draw (\x1-15,\y1) rectangle (\x1+15,\y1-30);
 \draw (\x1,\y1) -- (\x1,\y1+23);
 \node at (\x1,\y1-15) {$G'$};
 \putConvex(\x1,\y1+13)
 \node[above,inner sep=2] at (\x1,\y1+23) {$\scriptstyle\{\gamma\}$};
\end{tikzpicture}%
}}%
\kern3em\rightsquigarrow\kern3em
\vcenter{\hbox{%
\begin{tikzpicture}[xscale=0.0352778, yscale=0.0352778, thin, inner sep=0]
  \def\p#1#2{%
    \ifcase #1
      \or \ifx#2x  0       \else  0 \fi
    \fi}
  \def\x#1{\p#1x}
  \def\y#1{\p#1y}
 \path[use as bounding box] (-15,-30) rectangle (15,32);
 \draw (\x1-15,\y1) rectangle (\x1+15,\y1-30);
 \draw (\x1,\y1) -- (\x1,\y1+10);
 \node at (\x1,\y1-15) {$G'$};
 \node[above,inner sep=2] at (\x1,\y1+10) {$\scriptstyle\gamma$};
\end{tikzpicture}%
}}%
$$
We set $\pi(G)^{\{\gamma\}}=\pi(G')^\gamma$.
If the wire in the left
 graph is annotated by $\alpha$ that is not a singleton,
 the process is aborted immediately and returns
 $\pi(G)^{\alpha}=0$.
The remaining three cases are involved in non-determinacy.
The decomposition for a $d$-part is given as
$$\vcenter
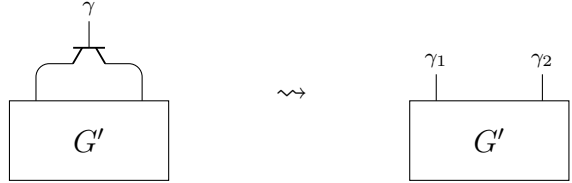
{\hbox{%
\begin{tikzpicture}[xscale=0.0352778, yscale=0.0352778, thin, inner sep=0]
  \def\p#1#2{%
    \ifcase #1
      \or \ifx#2x  0       \else  0 \fi
      \or \ifx#2x  \x1+10     \else  \y1+30 \fi
      \or \ifx#2x  \x1+50     \else  \y1+30 \fi
      \or \ifx#2x  \x1+60     \else  \y1+30 \fi
      \or \ifx#2x  \x1+30     \else  \y2+20 \fi
    \fi}
  \def\x#1{\p#1x}
  \def\y#1{\p#1y}
 \path[use as bounding box] (0,0) rectangle (62,66);
 \draw (\x1,\y1) rectangle (\x4,\y4);
 \draw (\x2,\y2)[rounded corners=5] -- (\x2,\y5-6) -- (\x5-6,\y5-6);
 \draw (\x3,\y3)[rounded corners=5] -- (\x3,\y5-6) -- (\x5+6,\y5-6);
 \draw (\x5,\y5) -- (\x5,\y5+10);
 \node[above,inner sep=2] at (\x5,\y5+10) {$\scriptstyle \gamma$};
 \node at (\x5,\y1+15) {$G'$};
 \putDuplicator(\x5,\y5)
\end{tikzpicture}%
}}%
\kern3em\rightsquigarrow\kern3em
\vcenter{\hbox{%
\begin{tikzpicture}[xscale=0.0352778, yscale=0.0352778, thin, inner sep=0]
  \def\p#1#2{%
    \ifcase #1
      \or \ifx#2x  0       \else  0 \fi
      \or \ifx#2x  \x1+10     \else  \y1+30 \fi
      \or \ifx#2x  \x1+50     \else  \y1+30 \fi
      \or \ifx#2x  \x1+60     \else  \y1+30 \fi
      \or \ifx#2x  \x1+30     \else  \y2+20 \fi
    \fi}
  \def\x#1{\p#1x}
  \def\y#1{\p#1y}
 \path[use as bounding box] (0,0) rectangle (62,66);
 \draw (\x1,\y1) rectangle (\x4,\y4);
 \draw (\x2,\y2)[rounded corners=5] -- (\x2,\y2+10);
 \draw (\x3,\y3)[rounded corners=5] -- (\x3,\y3+10);
 \node[above,text depth=4] at (\x2,\y2+10) {$\scriptstyle \gamma_1$};
 \node[above,text depth=4] at (\x3,\y3+10) {$\scriptstyle \gamma_2$};
 \node at (\x5,\y1+15) {$G'$};
\end{tikzpicture}%
}}$$
 where the sub-mutlisets $\gamma_1$ and $\gamma_2$ are chosen
 so that $\gamma=\gamma_1+\gamma_2$ holds.
For each of such pairs, a subprocess $\pi(G')^{\gamma_1\gamma_2}$ is
 invoked.
The number of created subprocesses is equal to the number of
 different ordered pairs $(\gamma_1,\gamma_2)$.
We remark that $\gamma_1+\gamma_2$ and $\gamma_2+\gamma_1$ are
 separately handled unless $\gamma_1=\gamma_2$.
Given the returns of the subprocesses, we
 set $\pi(G)^\gamma=\sum\pi(G')^{\gamma_1,\gamma_2}$.
The sum ranges over all successful subprocesses.
The process $\pi(G)$ succeeds if one of the subprocesses succeeds.
As the failing subprocess returns $0$, we
 may equivalently take the sum to
 range over all different ordered pairs $(\gamma_1,\gamma_2)$.
The decomposition of a $\delta$-part is given as follows:
$$
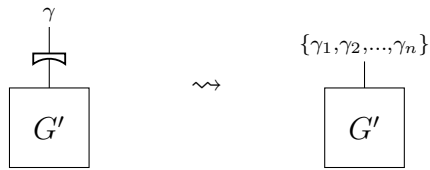
\vcenter{\hbox{%
\begin{tikzpicture}[xscale=0.0352778, yscale=0.0352778, thin, inner sep=0]
  \def\p#1#2{%
    \ifcase #1
      \or \ifx#2x  0       \else  0 \fi
    \fi}
  \def\x#1{\p#1x}
  \def\y#1{\p#1y}
 \path[use as bounding box] (-15,-30) rectangle (15,32);
 \draw (\x1-15,\y1) rectangle (\x1+15,\y1-30);
 \draw (\x1,\y1) -- (\x1,\y1+23);
 \node at (\x1,\y1-15) {$G'$};
 \putConcave(\x1,\y1+13)
 \node[above,inner sep=2] at (\x1,\y1+23) {$\scriptstyle\gamma$};
\end{tikzpicture}%
}}%
\kern3em\rightsquigarrow\kern3em
\vcenter{\hbox{%
\begin{tikzpicture}[xscale=0.0352778, yscale=0.0352778, thin, inner sep=0]
  \def\p#1#2{%
    \ifcase #1
      \or \ifx#2x  0       \else  0 \fi
    \fi}
  \def\x#1{\p#1x}
  \def\y#1{\p#1y}
 \path[use as bounding box] (-15,-30) rectangle (15,32);
 \draw (\x1-15,\y1) rectangle (\x1+15,\y1-30);
 \draw (\x1,\y1) -- (\x1,\y1+10);
 \node at (\x1,\y1-15) {$G'$};
 \node[above,inner sep=2] at (\x1,\y1+10) {$\scriptstyle\{\gamma_1,\gamma_2,
  \ldots,\gamma_n\}$};
\end{tikzpicture}%
}}%
$$
We suppose $\gamma=\gamma_1+\gamma_2+\cdots+\gamma_n$ where $n\geq 0$.
The order of $\gamma_i$ is irrelevant.
We invokde $\pi(G')^{\{\gamma_1,\gamma_2,\ldots,\gamma_n\}}$ for
 each of such decompositions.
We set $\pi(G)^\gamma=\sum \pi(G')^{\{\gamma_1,\gamma_2,\ldots,\gamma_n\}}$.

Next, We consider the decomposition of a board.
Suppose that the outgoing wires of a board are annotated
 as follows:
$$
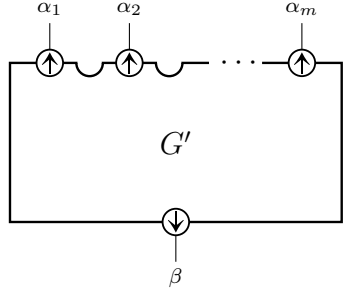
\vcenter{\hbox{%
\begin{tikzpicture}[xscale=0.0352778, yscale=0.0352778, thin, inner sep=0]
  \def\p#1#2{%
    \ifcase #1
      \or \ifx#2x  0       \else  0 \fi
      \or \ifx#2x  \x1+125 \else  \y1+60 \fi
      \or \ifx#2x  \x1+15  \else  \y2 \fi
      \or \ifx#2x  \x1+30  \else  \y2 \fi
      \or \ifx#2x  \x1+45  \else  \y2 \fi
      \or \ifx#2x  \x1+60  \else  \y2 \fi
      \or \ifx#2x  \x2-15  \else  \y2 \fi
    \fi}
  \def\x#1{\p#1x}
  \def\y#1{\p#1y}
 \draw (\x3,\y3) -- (\x3,\y3+15);
 \draw (\x5,\y5) -- (\x5,\y5+15);
 \draw (\x7,\y7) -- (\x7,\y7+15);
 \draw (\x1+62.5,\y1) -- (\x1+62.5,\y1-15);
 \boardBoundary(\x1,\y1)(\x2,\y2)
 \putLowerSocket(\x1+62.5,\y1)
 \putUpperSocket(\x3,\y3)
 \putDimple(\x4,\y4)
 \putUpperSocket(\x5,\y5)
 \putDimple(\x6,\y6)
 \putUpperSocket(\x7,\y7)
 \draw[white,line width=2.0] (\x1+75,\y2) -- ++(20,0);
 \node at (\x1+86,\y2) {$\cdots$};
 \node at ($ .5*(\x1,\y1)+.5*(\x2,\y2) $) {$G'$};
 \node[above,text depth=4] at (\x3,\y3+15) {$\scriptstyle\alpha_1$};
 \node[above,text depth=4] at (\x5,\y5+15) {$\scriptstyle\alpha_2$};
 \node[above,text depth=4] at (\x7,\y7+15) {$\scriptstyle\alpha_m$};
 \node[below,text height=8] at (\x1+62.5,\y1-15) {$\scriptstyle\beta$};
\end{tikzpicture}%
}}%
$$
First, we check whether all of $\alpha_1,\alpha_2,
\ldots,\alpha_m$ and $\beta$ have
 the same number of elements.
Otherwise the process fails and returns
 $\pi(G)^{\alpha_1\alpha_2\cdots\alpha_m}_\beta=0$.
Provided that the multisets have the same cardinality $n\geq 0$,
 let us put $\alpha_i=\{a_{i1},a_{i2},\ldots,a_{in}\}$ and
 $\beta=\{b_1,b_2,\ldots,b_n\}$.
We choose and fix a linear ordering
 $b_1,b_2,\ldots,b_n$ once and for all.
A {\it linear disposition} of $\alpha_i$ is a linear order obtained by
 shuffling the $n$ element of $\alpha_i$.
Namely, if $\sigma$ is a permutation over $n$ letters,
 a linear disposition of $\alpha_i$ is $a_{i\sigma(1)},a_{i\sigma(2)},
 \ldots,a_{i\sigma(n)}$.
Let $l_i$ denote the number of all linear dispositions of $\alpha_i$.
For instance, if the elements of $\alpha_i$ are all different then
$l_i=n\mathord!$; if all are identical then $l_i=1$.
We form an $m\times n$ matrix
$$\matrix{%
 a'_{11} & a'_{12} & \cdots & a'_{1n} \cr
 a'_{21} & a'_{22} & \cdots & a'_{2n} \cr
 \vdots & \vdots &         & \vdots \cr
 a'_{m1} & a'_{m2} & \cdots & a'_{mn} \cr
}$$
 where the $i$-th row is a linear disposition of $\alpha_i$.
The number of different matrices is equal to $l_1l_2\cdots l_m$.
Given a matrix, we spawn a subprocesses per column.
Namely, we invoke $n$ processes
 $\pi(G')^{a'_{11}a'_{21}\cdots a'_{m1}}_{b_1},\>
\pi(G')^{a'_{21}a'_{22}\cdots a'_{m2}}_{b_2},
\>\ldots,
\pi(G')^{a'_{1n}a'_{2n}\cdots a'_{mn}}_{b_n}$.
Then we set

	\vskip2ex
        \noindent\kern5em
$\displaystyle
\pi(G)^{\alpha_1\alpha_2\cdots\alpha_m}_\beta\ =
\ \sum \pi(G')^{a'_{11}a'_{21}\cdots a'_{m1}}_{b_1}
\pi(G')^{a'_{12}a'_{22}\cdots a'_{m2}}_{b_2}\>\cdots\>
\pi(G')^{a'_{1n}a'_{2n}\cdots a'_{mn}}_{b_n}$.
	\vskip2ex
        \noindent
The sum ranges over the different matrices.
We note that the order of the elements of $\beta$ is fixed,
 whereas the elements of $\alpha_i$ are shuffled.
This corresponds to the asymmetry of $\mathord!f$
 in the model $\mathscr{M}$.

Finally, suppose that the process reaches a wire of an atomic type
$$
\vcenter{\hbox{%
\begin{tikzpicture}[xscale=0.0352778, yscale=0.0352778, thin, inner sep=0]
 \draw (0,0) -- (0,18);
 \node[right,inner sep=2] at (0,18) {$\scriptstyle a$};
 \node[right,inner sep=2] at (0,0) {$\scriptstyle b$};
\end{tikzpicture}%
}}$$
If $a=b$, we set $\pi(G)^a_a=1$;
otherwise $\pi(G)^a_b$ fails and returns $0$.

Two types of concurrency appear.
One is non-determinacy, which is invoked at $d$-parts,
 $\delta$-parts, and boards by the choices of decomposition.
In the terminology of alternating computation models,
 it is the $\exists$-type non-determinism.
Namely, a process succeeds if one of the subprocesses succeeds.
In the decomposition of boards, we also have
 the $\forall$-type concurrency that
 succeeds only when all of the subprocesses succeed.
At a board $G$, the process create subprocesses
 by $\exists$-type non-determinism,
 which in turn create further subproceeses by
 the $\forall$-type concurrency.
Hence the value of $\pi(G)$ is defined by a sum of multiplications.

	\vskip2ex

\begin{lemma}
Let $G$ denote the graph associated with $f:A\rightarrow B$ in normal form.
Then $\pi(G)^{\alpha}_{\beta}=
M_f[\alpha;\beta]$ holds.
\end{lemma}

\proof
The definition of $\pi(G)$ simply rephrases the structure of $\mathscr{M}$
 in the terminology of processes.
\endproof

	\vskip2ex

\begin{lemma}
Every components $M_f[\alpha;\beta]$ is finite for every $f$
 of the calculus.
\end{lemma}

\proof
Contract $f$ to a normal form and turn it into a graph.
The process $\pi(G)$ only decomposes the graph.
The number of subprocesses is bounded by the tally
 of multiset decompositions.
Hence the numbers are finite.
\endproof

\section{Generic forms}

We introduce the notion of the generic form $\varphi(G)$
 of a graph $G$.
It is a kind of regular expressions.
In the theory of automata, we can restore an automaton from
 a regular expression \cite{sips}.
Moreover, the latter can be used as machinery to
 enumerate the language recognized by the automaton \cite{sefl}.
The generic form has similar functionality.
It can be used to restore the graph $G$ to some extent, as well as
 it gives information
 of $M[\alpha;\beta]$ for the corresponding matrix.
Our goal is to explore the relation between the syntax of the
 linear category and the model $\mathscr{M}$.
The generic form is a hinge between syntax and semantics.

We define the orientation of wires of a graph, then
 we associate the generic forms on the wires along
 the flows determined by the orientation.
To save space, however, we simultaneously give 
 the orientation and the association of generic forms.

First of all, we assign different identifiers $i,j,\ldots$
 to the boards occurring in the graph.
To each wire of an atomic type, we put arrowheads in
 both ends.
It is called a bioriented wire.
The generic form is given as
$$
\begin{tikzpicture}[xscale=0.0352778, yscale=0.0352778, thin, inner sep=0]
  \def\p#1#2{%
    \ifcase #1
      \or \ifx#2x  0       \else  0 \fi
    \fi}
  \def\x#1{\p#1x}
  \def\y#1{\p#1y}
 \draw[use as bounding box] (0,0) rectangle (0,20);
 \draw[>=latex,<->] (\x1,\y1) -- (\x1,\y1+20);
 \node[right,inner sep=2] at (\x1,\y1+10) {$\scriptstyle
   x_{i_1i_2\cdots i_n}$};
\end{tikzpicture}
$$
 where $x$ is a variable $i_1,i_2,\ldots,i_n$ are
 the list of the boards that contain this wire, ordered from the innermost.
We assume that different bioriented wires have different variables.
We often call the whole $x_{i_1i_2\ldots i_n}$ a variable.
The orientation and generic form of tensor parts are
$$
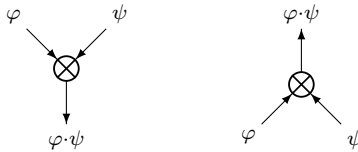
\begin{tikzpicture}[xscale=0.0352778, yscale=0.0352778, thin, inner sep=0]
  \def\p#1#2{%
    \ifcase #1
      \or \ifx#2x  0       \else  0 \fi
    \fi}
  \def\x#1{\p#1x}
  \def\y#1{\p#1y}
  \path[use as bounding box] (\x1-15,\y1-29) rectangle (\x1+15,\y1+23);
  \draw (\x1-15,\y1+15) -- (\x1,\y1);
  \draw (\x1+15,\y1+15) -- (\x1,\y1);
  \draw (\x1,\y1-21) -- (\x1,\y1);
  \draw[>=latex,->] (\x1-4,\y1+4) -- (\x1-3.5,\y1+3.5);
  \draw[>=latex,->] (\x1+4,\y1+4) -- (\x1+3.5,\y1+3.5);
  \draw[>=latex,->] (\x1,\y1-20) -- (\x1,\y1-21);
  \putTensor(\x1,\y1)
  \node[above left,inner sep=2] at (\x1-15,\y1+15) {$\scriptstyle \varphi$};
  \node[above right,inner sep=2] at (\x1+15,\y1+15) {$\scriptstyle \psi$};
  \node[below,inner sep=2] at (\x1,\y1-21) {$\scriptstyle \varphi\cdot\psi$};
\end{tikzpicture}
	\hskip5em
\begin{tikzpicture}[xscale=0.0352778, yscale=0.0352778, thin, inner sep=0]
  \def\p#1#2{%
    \ifcase #1
      \or \ifx#2x  0       \else  0 \fi
    \fi}
  \def\x#1{\p#1x}
  \def\y#1{\p#1y}
  \path[use as bounding box] (\x1-15,\y1+29) rectangle (\x1+15,\y1-23);
  \draw (\x1-15,\y1-15) -- (\x1,\y1);
  \draw (\x1+15,\y1-15) -- (\x1,\y1);
  \draw (\x1,\y1+21) -- (\x1,\y1);
  \draw[>=latex,->] (\x1-4,\y1-4) -- (\x1-3.5,\y1-3.5);
  \draw[>=latex,->] (\x1+4,\y1-4) -- (\x1+3.5,\y1-3.5);
  \draw[>=latex,->] (\x1,\y1+20) -- (\x1,\y1+21);
  \putTensor(\x1,\y1)
  \node[below left,inner sep=2] at (\x1-15,\y1-15) {$\scriptstyle \varphi$};
  \node[below right,inner sep=2] at (\x1+15,\y1-15) {$\scriptstyle \psi$};
  \node[above,inner sep=2] at (\x1,\y1+21) {$\scriptstyle \varphi\cdot\psi$};
\end{tikzpicture}
$$
For the cotensor parts, we similarly set as follows:
$$
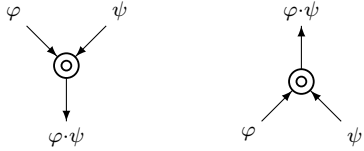
\begin{tikzpicture}[xscale=0.0352778, yscale=0.0352778, thin, inner sep=0]
  \def\p#1#2{%
    \ifcase #1
      \or \ifx#2x  0       \else  0 \fi
    \fi}
  \def\x#1{\p#1x}
  \def\y#1{\p#1y}
  \path[use as bounding box] (\x1-15,\y1-29) rectangle (\x1+15,\y1+23);
  \draw (\x1-15,\y1+15) -- (\x1,\y1);
  \draw (\x1+15,\y1+15) -- (\x1,\y1);
  \draw (\x1,\y1-21) -- (\x1,\y1);
  \draw[>=latex,->] (\x1-4,\y1+4) -- (\x1-3.5,\y1+3.5);
  \draw[>=latex,->] (\x1+4,\y1+4) -- (\x1+3.5,\y1+3.5);
  \draw[>=latex,->] (\x1,\y1-20) -- (\x1,\y1-21);
  \putPar(\x1,\y1)
  \node[above left,inner sep=2] at (\x1-15,\y1+15) {$\scriptstyle \varphi$};
  \node[above right,inner sep=2] at (\x1+15,\y1+15) {$\scriptstyle \psi$};
  \node[below,inner sep=2] at (\x1,\y1-21) {$\scriptstyle \varphi\cdot\psi$};
\end{tikzpicture}
	\hskip5em
\begin{tikzpicture}[xscale=0.0352778, yscale=0.0352778, thin, inner sep=0]
  \def\p#1#2{%
    \ifcase #1
      \or \ifx#2x  0       \else  0 \fi
    \fi}
  \def\x#1{\p#1x}
  \def\y#1{\p#1y}
  \path[use as bounding box] (\x1-15,\y1+29) rectangle (\x1+15,\y1-23);
  \draw (\x1-15,\y1-15) -- (\x1,\y1);
  \draw (\x1+15,\y1-15) -- (\x1,\y1);
  \draw (\x1,\y1+21) -- (\x1,\y1);
  \draw[>=latex,->] (\x1-4,\y1-4) -- (\x1-3.5,\y1-3.5);
  \draw[>=latex,->] (\x1+4,\y1-4) -- (\x1+3.5,\y1-3.5);
  \draw[>=latex,->] (\x1,\y1+20) -- (\x1,\y1+21);
  \putPar(\x1,\y1)
  \node[below left,inner sep=2] at (\x1-15,\y1-15) {$\scriptstyle \varphi$};
  \node[below right,inner sep=2] at (\x1+15,\y1-15) {$\scriptstyle \psi$};
  \node[above,inner sep=2] at (\x1,\y1+21) {$\scriptstyle \varphi\cdot\psi$};
\end{tikzpicture}
$$
For the units and counits, we set
$$
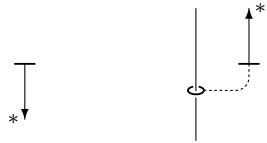
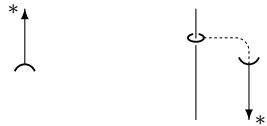
\begin{tikzpicture}[xscale=0.0352778, yscale=0.0352778, thin, inner sep=0]
  \def\p#1#2{%
    \ifcase #1
      \or \ifx#2x  0       \else  0 \fi
    \fi}
  \def\x#1{\p#1x}
  \def\y#1{\p#1y}
  \path[use as bounding box] (\x1-3,\y1-29) rectangle (\x1+3,\y1+29);
  \draw (\x1,\y1-21) -- (\x1,\y1);
  \draw[>=latex,->] (\x1,\y1-20) -- (\x1,\y1-21);
  \putPositiveTerminal(\x1,\y1)
  \node[left,inner sep=2] at (\x1,\y1-21) {$\scriptstyle *$};
\end{tikzpicture}
	\hskip5em
\begin{tikzpicture}[xscale=0.0352778, yscale=0.0352778, thin, inner sep=0]
  \def\p#1#2{%
    \ifcase #1
      \or \ifx#2x  0       \else  0 \fi
    \fi}
  \def\x#1{\p#1x}
  \def\y#1{\p#1y}
  \path[use as bounding box] (\x1-23,\y1-29) rectangle (\x1+3,\y1+29);
  \draw (\x1,\y1+21) -- (\x1,\y1);
  \draw (\x1-20,\y1+21) -- (\x1-20,\y1-29);
  \draw[>=latex,->] (\x1,\y1+20) -- (\x1,\y1+21);
  \putPositiveTerminal(\x1,\y1)
  \putRubberBand(\x1-20,\y1-10)
  \draw[dash pattern=on 1 off 1] (\x1,\y1)[rounded corners=5] -- (\x1,\y1-10)
    -- (\x1-17,\y1-10);
  \node[right,inner sep=2] at (\x1,\y1+21) {$\scriptstyle *$};
\end{tikzpicture}
$$
 and
$$
\begin{tikzpicture}[xscale=0.0352778, yscale=0.0352778, thin, inner sep=0]
  \def\p#1#2{%
    \ifcase #1
      \or \ifx#2x  0       \else  0 \fi
    \fi}
  \def\x#1{\p#1x}
  \def\y#1{\p#1y}
  \path[use as bounding box] (\x1-3,\y1-29) rectangle (\x1+3,\y1+29);
  \draw (\x1,\y1+21) -- (\x1,\y1);
  \draw[>=latex,->] (\x1,\y1+20) -- (\x1,\y1+21);
  \putLowerNegativeTerminal(\x1,\y1)
  \node[left,inner sep=2] at (\x1,\y1+21) {$\scriptstyle *$};
\end{tikzpicture}
	\hskip5em
\begin{tikzpicture}[xscale=0.0352778, yscale=0.0352778, thin, inner sep=0]
  \def\p#1#2{%
    \ifcase #1
      \or \ifx#2x  0       \else  0 \fi
    \fi}
  \def\x#1{\p#1x}
  \def\y#1{\p#1y}
  \path[use as bounding box] (\x1-23,\y1-29) rectangle (\x1+3,\y1+29);
  \draw (\x1,\y1-21) -- (\x1,\y1);
  \draw (\x1-20,\y1-21) -- (\x1-20,\y1+21);
  \draw[>=latex,->] (\x1,\y1-20) -- (\x1,\y1-21);
  \putUpperNegativeTerminal(\x1,\y1)
  \putRubberBand(\x1-20,\y1+10)
  \draw[dash pattern=on 1 off 1] (\x1,\y1)[rounded corners=5] -- (\x1,\y1+10)
    -- (\x1-17,\y1+10);
  \node[right,inner sep=2] at (\x1,\y1-21) {$\scriptstyle *$};
\end{tikzpicture}
$$
For the duality parts, we set
$$
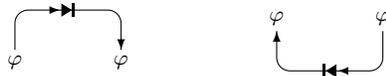
\begin{tikzpicture}[xscale=0.0352778, yscale=0.0352778, thin, inner sep=0]
  \def\p#1#2{%
    \ifcase #1
      \or \ifx#2x  0       \else  0 \fi
    \fi}
  \def\x#1{\p#1x}
  \def\y#1{\p#1y}
  \path[use as bounding box] (\x1-20,\y1-23) rectangle (\x1+20,\y1);
  \draw (\x1+20,\y1-15)[rounded corners=5] -- (\x1+20,\y1) -- (\x1,\y1);
  \draw (\x1,\y1)[rounded corners=5] -- (\x1-20,\y1) -- (\x1-20,\y1-15);
  \draw[>=latex,->] (\x1-4,\y1) -- (\x1-3,\y1);
  \draw[>=latex,->] (\x1+20,\y1-14) -- (\x1+20,\y1-15);
  \putRightDiode(\x1,\y1)
  \node[below,inner sep=2] at (\x1-20,\y1-15) {$\scriptstyle \varphi$};
  \node[below,inner sep=2] at (\x1+20,\y1-15) {$\scriptstyle \varphi$};
\end{tikzpicture}
	\hskip5em
\begin{tikzpicture}[xscale=0.0352778, yscale=0.0352778, thin, inner sep=0]
  \def\p#1#2{%
    \ifcase #1
      \or \ifx#2x  0       \else  0 \fi
    \fi}
  \def\x#1{\p#1x}
  \def\y#1{\p#1y}
  \path[use as bounding box] (\x1-20,\y1) rectangle (\x1+20,\y1+23);
  \draw (\x1+20,\y1+15)[rounded corners=5] -- (\x1+20,\y1) -- (\x1,\y1);
  \draw (\x1,\y1)[rounded corners=5] -- (\x1-20,\y1) -- (\x1-20,\y1+15);
  \draw[>=latex,->] (\x1+4,\y1) -- (\x1+3,\y1);
  \draw[>=latex,->] (\x1-20,\y1+14) -- (\x1-20,\y1+15);
  \putLeftDiode(\x1,\y1)
  \node[above,inner sep=2] at (\x1-20,\y1+15) {$\scriptstyle \varphi$};
  \node[above,inner sep=2] at (\x1+20,\y1+15) {$\scriptstyle \varphi$};
\end{tikzpicture}
$$
The generic form does not alter, while the orientation
 of the flow rotates $180^\circ$.
The orientations and the generic forms for
 a $\delta$-part, an $\varepsilon$-part, a $d$-part, and an
 $e$-part is
$$
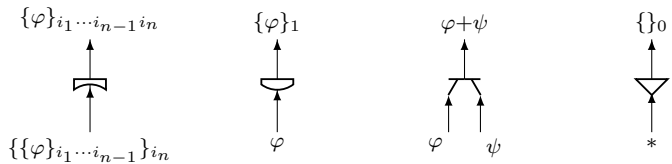
\begin{tikzpicture}[xscale=0.0352778, yscale=0.0352778, thin, inner sep=0]
  \def\p#1#2{%
    \ifcase #1
      \or \ifx#2x  0       \else  0 \fi
    \fi}
  \def\x#1{\p#1x}
  \def\y#1{\p#1y}
  \path[use as bounding box] (\x1-6,\y1+22) rectangle (\x1+6,\y1-27);
  \draw (\x1,\y1+15) -- (\x1,\y1-20);
  \draw[>=latex,->] (\x1,\y1-4) -- (\x1,\y1-3);
  \draw[>=latex,->] (\x1,\y1+14) -- (\x1,\y1+15);
  \putConcave(\x1,\y1)
  \node[above,inner sep=2] at (\x1,\y1+15) {$\scriptstyle\{\varphi\}_{i_1\cdots i_{n-1}i_n}$};
  \node[below,inner sep=2] at (\x1,\y1-20) {$\scriptstyle\{\{\varphi\}_{i_1\cdots i_{n-1}}\}_{i_n}$};
\end{tikzpicture}%
	\hskip5em
\begin{tikzpicture}[xscale=0.0352778, yscale=0.0352778, thin, inner sep=0]
  \def\p#1#2{%
    \ifcase #1
      \or \ifx#2x  0       \else  0 \fi
    \fi}
  \def\x#1{\p#1x}
  \def\y#1{\p#1y}
  \path[use as bounding box] (\x1-6,\y1+22) rectangle (\x1+6,\y1-27);
  \draw (\x1,\y1+15) -- (\x1,\y1-20);
  \draw[>=latex,->] (\x1,\y1-5) -- (\x1,\y1-4);
  \draw[>=latex,->] (\x1,\y1+14) -- (\x1,\y1+15);
  \putConvex(\x1,\y1)
  \node[above,inner sep=2] at (\x1,\y1+15) {$\scriptstyle\{\varphi\}_1$};
  \node[below,inner sep=2] at (\x1,\y1-20) {$\scriptstyle \varphi$};
\end{tikzpicture}%
	\hskip5em
\begin{tikzpicture}[xscale=0.0352778, yscale=0.0352778, thin, inner sep=0]
  \def\p#1#2{%
    \ifcase #1
      \or \ifx#2x  0       \else  0 \fi
    \fi}
  \def\x#1{\p#1x}
  \def\y#1{\p#1y}
  \path[use as bounding box] (\x1-6,\y1+22) rectangle (\x1+6,\y1-27);
  \draw (\x1,\y1+15) -- (\x1,\y1);
  \draw[>=latex,->] (\x1-6,\y1-20) -- (\x1-6,\y1-6);
  \draw[>=latex,->] (\x1+6,\y1-20) -- (\x1+6,\y1-6);
  \draw[>=latex,->] (\x1,\y1+14) -- (\x1,\y1+15);
  \putDuplicator(\x1,\y1)
  \node[above,inner sep=2] at (\x1,\y1+15) {$\scriptstyle\mathord\varphi+\psi$};
  \node[below left,inner sep=2] at (\x1-6,\y1-20) {$\scriptstyle\varphi$};
  \node[below right,inner sep=2] at (\x1+6,\y1-20) {$\scriptstyle\psi$};
\end{tikzpicture}%
	\hskip5em
\begin{tikzpicture}[xscale=0.0352778, yscale=0.0352778, thin, inner sep=0]
  \def\p#1#2{%
    \ifcase #1
      \or \ifx#2x  0       \else  0 \fi
    \fi}
  \def\x#1{\p#1x}
  \def\y#1{\p#1y}
  \path[use as bounding box] (\x1-6,\y1+22) rectangle (\x1+6,\y1-27);
  \draw (\x1,\y1+15) -- (\x1,\y1-20);
  \draw[>=latex,->] (\x1,\y1-6) -- (\x1,\y1-5);
  \draw[>=latex,->] (\x1,\y1+14) -- (\x1,\y1+15);
  \putEliminator(\x1,\y1)
  \node[above,inner sep=2] at (\x1,\y1+15) {$\scriptstyle\{\}_0$};
  \node[below,inner sep=2] at (\x1,\y1-20) {$\scriptstyle *$};
\end{tikzpicture}%
$$
We recall that a $\delta$-part must have a gate of
 a board immediately below.
Hence the generic form beneath the $\delta$-part is justified
 from the case of a board below.
For the $d$-part, we do not distinguish between $\varphi+\psi$
 and $\psi+\varphi$.
For a board $i$, we set
$$
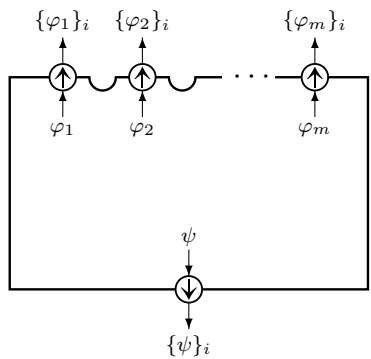
\vcenter{\hbox{%
\begin{tikzpicture}[xscale=0.0352778, yscale=0.0352778, thin, inner sep=0]
  \def\p#1#2{%
    \ifcase #1
      \or \ifx#2x  0       \else  0 \fi
      \or \ifx#2x  \x1+125 \else  \y1+80 \fi
      \or \ifx#2x  \x1+15  \else  \y2 \fi
      \or \ifx#2x  \x1+30  \else  \y2 \fi
      \or \ifx#2x  \x1+45  \else  \y2 \fi
      \or \ifx#2x  \x1+60  \else  \y2 \fi
      \or \ifx#2x  \x2-15  \else  \y2 \fi
    \fi}
  \def\x#1{\p#1x}
  \def\y#1{\p#1y}
 \draw[](\x3,\y3-15) -- (\x3,\y3+15);
 \draw[](\x5,\y5-15) -- (\x5,\y5+15);
 \draw[](\x7,\y7-15) -- (\x7,\y7+15);
 \draw[](\x1+62.5,\y1+15) -- (\x1+62.5,\y1-15);
 \draw[>=latex,->] (\x3,\y3+14) -- (\x3,\y3+15);
 \draw[>=latex,->] (\x5,\y5+14) -- (\x5,\y5+15);
 \draw[>=latex,->] (\x7,\y7+14) -- (\x7,\y7+15);
 \draw[>=latex,->] (\x1+62.5,\y1-14) -- (\x1+62.5,\y1-15);
 \draw[>=latex,->] (\x3,\y3-6) -- (\x3,\y3-5);
 \draw[>=latex,->] (\x5,\y5-6) -- (\x5,\y5-5);
 \draw[>=latex,->] (\x7,\y7-6) -- (\x7,\y7-5);
 \draw[>=latex,->] (\x1+62.5,\y1+6) -- (\x1+62.5,\y1+5);
 \boardBoundary(\x1-5,\y1)(\x2+5,\y2)
 \putLowerSocket(\x1+62.5,\y1)
 \putUpperSocket(\x3,\y3)
 \putDimple(\x4,\y4)
 \putUpperSocket(\x5,\y5)
 \putDimple(\x6,\y6)
 \putUpperSocket(\x7,\y7)
 \draw[white,line width=2.0] (\x1+75,\y2) -- ++(20,0);
 \node at (\x1+86,\y2) {$\cdots$};
 \node[above,text depth=4] at (\x3,\y3+15) {$\scriptstyle \{\varphi_1\}_i$};
 \node[above,text depth=4] at (\x5,\y5+15) {$\scriptstyle \{\varphi_2\}_i$};
 \node[above,text depth=4] at (\x7,\y7+15) {$\scriptstyle \{\varphi_m\}_i$};
 \node[below,text height=8] at (\x1+62.5,\y1-15) {$\scriptstyle \{\psi\}_i$};
 \node[below,text depth=4,inner sep=2] at (\x3,\y3-15) {$\scriptstyle \varphi_1$};
 \node[below,text depth=4,inner sep=2] at (\x5,\y5-15) {$\scriptstyle \varphi_2$};
 \node[below,text depth=4,inner sep=2] at (\x7,\y7-15) {$\scriptstyle \varphi_m$};
 \node[above,text height=8,inner sep=1] at (\x1+62.5,\y1+15) {$\scriptstyle \psi$};
\end{tikzpicture}%
}}%
$$

	\vskip3ex

\begin{lemma}
Flows never collide and never loop.
\end{lemma}

\proof
Collision occurs at a $\beta$-redex or at a board located immediately
 above a $\delta$-part, an $\varepsilon$-part, a $d$-part, or an
 $e$-part.
By our convention on graphs, there are no $\beta$-redexes
 and no $\eta$-expansion is applied above these four parts.
So there is no collision.
The changes of directions up/down are caused only by
 duality parts, which changes the type $A$ to more complex $A^*$.
Hence, if a flow returned to the same place, the type of the
 wire could not be equal to the original.
So there is no loop.
\endproof

	\vskip2ex

\begin{figure}[t]
\figbox{%
$$
\begin{tikzpicture}[xscale=0.0352778, yscale=0.0352778, thin, inner sep=0]
  \def\pa#1#2{%
    \ifcase #1
      \or \ifx#2x  0           \else  0 \fi
      \or \ifx#2x  \xa1+25     \else  \ya1-20 \fi
      \or \ifx#2x  \xa2+15     \else  \ya2+20 \fi
      \or \ifx#2x  \xa3+10     \else  \ya3+35 \fi
      \or \ifx#2x  \xa4+40     \else  \ya4-20 \fi
      \or \ifx#2x  \xa5+20     \else  \ya5+20 \fi
      \or \ifx#2x  \xa6        \else  \ya6+20 \fi
      \or \ifx#2x  \xa7+10     \else  \ya7+45 \fi
      \or \ifx#2x  \xa8-10     \else  \ya8-10 \fi
    \fi}
  \def\xa#1{\pa#1x}
  \def\ya#1{\pa#1y}
  \def\pb#1#2{%
    \ifcase #1
          \ifx#2x  \xa9        \else  \ya9-15 \fi
      \or \ifx#2x  \xa7-20     \else  \ya7+80 \fi
      \or \ifx#2x  \xb1+25     \else  \yb1-20 \fi
      \or \ifx#2x  \xb2+15     \else  \yb2+20 \fi
      \or \ifx#2x  \xb3        \else  \yb3+20 \fi
      \or \ifx#2x  \xb4+40     \else  \yb4+70 \fi
      \or \ifx#2x  \xb5-10     \else  \yb5-15 \fi
      \or \ifx#2x  \xb6        \else  \yb6-20 \fi
      \or \ifx#2x  \xb7-10     \else  \yb7-15 \fi
      \or \ifx#2x  \xb8-20     \else  \yb8+25 \fi
    \fi}
  \def\xb#1{\pb#1x}
  \def\yb#1{\pb#1y}
  \def\pc#1#2{%
    \ifcase #1
          \ifx#2x  \xb9        \else  \yb9+40 \fi
      \or \ifx#2x  \xc0        \else  \yc0+20 \fi
      \or \ifx#2x  \xc1        \else  \yc1+20 \fi
      \or \ifx#2x  \xc2        \else  \yc2+15 \fi
      \or \ifx#2x  \xc3        \else  \yc3+20 \fi
      \or \ifx#2x  \xc4-15     \else  \yc4+25 \fi
      \or \ifx#2x  \xc5-25     \else  \yc5+20 \fi
      \or \ifx#2x  \xb1        \else  \yc1 \fi
      \or \ifx#2x  \xc7        \else  \yc3 \fi
      \or \ifx#2x  \xc8-20     \else  \yc5 \fi
    \fi}
  \def\xc#1{\pc#1x}
  \def\yc#1{\pc#1y}
  \def\pd#1#2{%
    \ifcase #1
          \ifx#2x  \xa4        \else  \yc3 \fi
      \or \ifx#2x  \xa1        \else  \yc4 \fi
      \or \ifx#2x  \xc7        \else  \yc4 \fi
    \fi}
  \def\xd#1{\pd#1x}
  \def\yd#1{\pd#1y}
%
 \path[use as bounding box] (-25,-85) rectangle (225,445);
 \draw (\xa1,\ya1)[rounded corners=4] -- (\xa1-10,\ya1-10) -- (\xa1-10,\ya1-25);
 \draw (\xa1,\ya1)[rounded corners=4] -- (\xa1+10,\ya1-10) -- (\xa1+10,\ya1-20) -- (\xa2,\ya2);
 \draw (\xa2,\ya2)[rounded corners=4] -- (\xa3,\ya3-20) -- (\xa3,\ya3);
 \draw (\xa3,\ya3)[rounded corners=4] -- (\xa4-10,\ya4-10) -- (\xa4,\ya4);
 \draw (\xa4,\ya4)[rounded corners=4] -- (\xa4+10,\ya4-10) -- (\xa4+10,\ya5) -- (\xa5,\ya5);
 \draw (\xa5,\ya5)[rounded corners=4] -- (\xa6,\ya5) -- (\xa6,\ya6);
 \draw (\xa6,\ya6) -- (\xa7,\ya7);
 \draw (\xa7,\ya7)[rounded corners=4] -- (\xa7+10,\ya7+10) -- (\xa8+10,\ya7+10)
   -- (\xa8+10,\ya8) -- (\xa8,\ya8);
 \draw (\xa8,\ya8)[rounded corners=4] -- (\xa9,\ya8) -- (\xa9,\ya9);
 \draw (\xa9,\ya9) -- (\xb0,\yb0);
 \draw (\xa7,\ya7)[rounded corners=4] -- (\xa7-30,\ya7+30) -- (\xb1-10,\yb1-10)
   -- (\xb1,\yb1);
 \draw (\xb1,\yb1)[rounded corners=4] -- (\xb1+10,\yb1-10) -- (\xb1+10,\yb2) -- (\xb2,\yb2);
 \draw (\xb2,\yb2)[rounded corners=4] -- (\xb3,\yb2) -- (\xb3,\yb3);
 \draw (\xb3,\yb3) -- (\xb4,\yb4);
 \draw (\xb4,\yb4)[rounded corners=4] -- (\xb4+10,\yb4+10) -- (\xb5+10,\yb4+10)
   -- (\xb5+10,\yb5) -- (\xb5,\yb5);
 \draw (\xb5,\yb5)[rounded corners=4] -- (\xb6,\yb5) -- (\xb6,\yb6);
 \draw (\xb6,\yb6) -- (\xb7,\yb7);
 \draw (\xb7,\yb7)[rounded corners=4] -- (\xb7,\yb8) -- (\xb8,\yb8);
 \draw (\xb8,\yb8)[rounded corners=4] -- (\xb9+10,\yb8) -- (\xb9+10,\yb9-10) -- (\xb9,\yb9);
 \draw (\xb4,\yb4)[rounded corners=4] -- (\xb4-10,\yb4+10) -- (\xb9-10,\yb9-10) -- (\xb9,\yb9);
 \draw (\xb9,\yb9) -- (\xc0,\yc0);
 \draw (\xc0,\yc0) -- (\xc2,\yc2);
 \draw (\xc2,\yc2) -- (\xc3,\yc3);
 \draw (\xc3,\yc3) -- (\xc4,\yc4);
 \draw (\xc4,\yc4)[rounded corners=4] -- (\xc4,\yc5-6) -- (\xc5+6,\yc5-6);
 \draw (\xc5,\yc5) -- (\xc6,\yc6);
 \draw (\xb1,\yb1) -- (\xc7,\yc7);
 \draw (\xc7,\yc7) -- (\xd2,\yd2);
 \draw (\xd2,\yd2)[rounded corners=4] -- (\xc8,\yc9-6) -- (\xc9+6,\yc9-6);
 \draw (\xc6,\yc6) -- (\xc9,\yc9);
 \draw (\xa4,\ya4) -- (\xd0,\yd0);
 \draw (\xd0,\yd0)[rounded corners=4] -- (\xd0,\yc9-6) -- (\xc9-6,\yc9-6);
 \draw (\xa1,\ya1) -- (\xd1,\yd1);
 \draw[white,line width=2.0] (\xd1,\yc5-16) -- (\xc5-6,\yc5-16);
 \draw (\xd1,\yd1)[rounded corners=4] -- (\xd1,\yc5-16) -- (\xc5-6,\yc5-16)
   -- (\xc5-6,\yc5-6);
 \draw (\xc6,\yc6) -- (\xc6,\yc6+15);
 \draw[dash pattern=on 1 off 1] (\xb0,\yb0) -- (\xa7-8,\ya7+12);
 \begin{scope}[shift={(\xa7-10,\ya7+10)},rotate=45]
  \putRubberBand(0,0)
 \end{scope}
 \boardBoundary(\xc3+95,\yc3)(\xa3-15,\ya3)
 \boardBoundary(\xc7-20,\yc7)(\xa6+100,\ya6)
 \boardBoundary(\xb3-20,\yb3)(\xc0+65,\yc0)
 \boardBoundary(\xb6-10,\yb6)(\xb7+10,\yb7)
 \putPar(\xa1,\ya1)
 \putLeftDiode(\xa2,\ya2)
 \putLowerSocket(\xa3,\ya3)
 \putPar(\xa4,\ya4)
 \putLeftDiode(\xa5,\ya5)
 \putLowerSocket(\xa6,\ya6)
 \putPar(\xa7,\ya7)
 \putRightDiode(\xa8,\ya8)
 \putEliminator(\xa9,\ya9)
 \putPositiveTerminal(\xb0,\yb0)
 \putPar(\xb1,\yb1)
 \putLeftDiode(\xb2,\yb2)
 \putLowerSocket(\xb3,\yb3)
 \putPar(\xb4,\yb4)
 \putRightDiode(\xb5,\yb5)
 \putUpperSocket(\xb6,\yb6)
 \putLowerSocket(\xb7,\yb7)
 \putLeftDiode(\xb8,\yb8)
 \putPar(\xb9,\yb9)
 \putUpperSocket(\xc0,\yc0)
 \putUpperSocket(\xc1,\yc1)
 \putConcave(\xc2,\yc2)
 \putUpperSocket(\xc3,\yc3)
 \putConcave(\xc4,\yc4)
 \putDuplicator(\xc5,\yc5)
 \putTensor(\xc6,\yc6)
 \putUpperSocket(\xc7,\yc7)
 \putUpperSocket(\xc8,\yc8)
 \putDuplicator(\xc9,\yc9)
 \putUpperSocket(\xd0,\yd0)
 \putConvex(\xd1,\yd1)
 \putConcave(\xd2,\yd2)
 \putDimple(\xc1-20,\yc1)
 \putDimple(\xc3-20,\yc3)
 \putDimple(\xc8-20,\yc8)
 \draw[>=latex,->] (\xa1-10,\ya1-24) -- (\xa1-10,\ya1-25);
 \draw[>=latex,->] (\xa1-4,\ya1-4) -- (\xa1-3.5,\ya1-3.5);
 \draw[>=latex,->] (\xa1+4,\ya1-4) -- (\xa1+3.5,\ya1-3.5);
 \draw[>=latex,->] (\xa2+3,\ya2) -- (\xa2+2.5,\ya2);
 \draw[>=latex,->] (\xa3,\ya3+5.5) -- (\xa3,\ya3+5);
 \draw[>=latex,->] (\xa4-4,\ya4-4) -- (\xa4-3.5,\ya4-3.5);
 \draw[>=latex,->] (\xa4+4,\ya4-4) -- (\xa4+3.5,\ya4-3.5);
 \draw[>=latex,->] (\xa5+3,\ya5) -- (\xa5+2.5,\ya5);
 \draw[>=latex,->] (\xa6,\ya6+5.5) -- (\xa6,\ya6+5);
 \draw[>=latex,->] (\xa7-4,\ya7+4) -- (\xa7-3.5,\ya7+3.5);
 \draw[>=latex,->] (\xa7+4,\ya7+4) -- (\xa7+3.5,\ya7+3.5);
 \draw[>=latex,->] (\xa8-3,\ya8) -- (\xa8-2.5,\ya8);
 \draw[>=latex,->] (\xa9,\ya9-6.5) -- (\xa9,\ya9-6);
 \draw[>=latex,->] (\xb1-4,\yb1-4) -- (\xb1-3.5,\yb1-3.5);
 \draw[>=latex,->] (\xb1+4,\yb1-4) -- (\xb1+3.5,\yb1-3.5);
 \draw[>=latex,->] (\xb2+3,\yb2) -- (\xb2+2.5,\yb2);
 \draw[>=latex,->] (\xb3,\yb3+5.5) -- (\xb3,\yb3+5);
 \draw[>=latex,->] (\xb4-4,\yb4+4) -- (\xb4-3.5,\yb4+3.5);
 \draw[>=latex,->] (\xb4+4,\yb4+4) -- (\xb4+3.5,\yb4+3.5);
 \draw[>=latex,->] (\xb5-3,\yb5) -- (\xb5-2.5,\yb5);
 \draw[>=latex,->] (\xb6,\yb6-5.5) -- (\xb6,\yb6-5);
 \draw[>=latex,->] (\xb7,\yb7+5.5) -- (\xb7,\yb7+5);
 \draw[>=latex,->] (\xb8+3,\yb8) -- (\xb8+2.5,\yb8);
 \draw[>=latex,->] (\xb9-4,\yb9-4) -- (\xb9-3.5,\yb9-3.5);
 \draw[>=latex,->] (\xb9+4,\yb9-4) -- (\xb9+3.5,\yb9-3.5);
 \draw[>=latex,->] (\xc0,\yc0-5.5) -- (\xc0,\yc0-5);
 \draw[>=latex,->] (\xc1,\yc1-5.5) -- (\xc1,\yc1-5);
 \draw[>=latex,->] (\xc2,\yc2-3) -- (\xc2,\yc2-2.5);
 \draw[>=latex,->] (\xc3,\yc3-5.5) -- (\xc3,\yc3-5);
 \draw[>=latex,->] (\xc4,\yc4-3) -- (\xc4,\yc4-2.5);
 \draw[>=latex,->] (\xc5+6.5,\yc5-6) -- (\xc5+6,\yc5-6);
 \draw[>=latex,->] (\xc5-6,\yc5-6.5) -- (\xc5-6,\yc5-6);
 \draw[>=latex,->] (\xc6+4.4,\yc6-3.5) -- (\xc6+3.9,\yc6-3.1);
 \draw[>=latex,->] (\xc6-4,\yc6-4) -- (\xc6-3.5,\yc6-3.5);
 \draw[>=latex,->] (\xc7,\yc7-5.5) -- (\xc7,\yc7-5);
 \draw[>=latex,->] (\xc8,\yc8-5.5) -- (\xc8,\yc8-5);
 \draw[>=latex,->] (\xc9+6.5,\yc9-6) -- (\xc9+6,\yc9-6);
 \draw[>=latex,->] (\xc9-6.5,\yc9-6) -- (\xc9-6,\yc9-6);
 \draw[>=latex,->] (\xd0,\yd0-6.5) -- (\xd0,\yd0-5);
 \draw[>=latex,->] (\xd1,\yd1-5.5) -- (\xd1,\yd1-4);
 \draw[>=latex,->] (\xc6,\yc6+14.5) -- (\xc6,\yc6+15);
 \filldraw[white] (\xb7+10-2,\yb7+10-4) rectangle (\xb7+10+2,\yb7+10+4);
 \node[right] at ($ .5*(\xb6,\yb6)+.5*(\xb6+4,\yb7) $) {$\scriptstyle v_{ijkl}$};
 \node[right] at ($ .5*(\xb4-10,\yb4)+.5*(\xb4-6,\yb9) $) {$\scriptstyle w_{jkl}$};
 \node[right] at (\xb1-8,\yb1-30) {$\scriptstyle z_{kl}$};
 \node[right] at (\xa3+2,\ya3+17) {$\scriptstyle y_{l}$};
 \node[left] at (\xa1-12,\ya1-15) {$\scriptstyle x$};
 \node[right] at (\xc6+2,\yc6+15) {$\scriptstyle
  (\{y_{l}\cdot\{z_{kl}\cdot\{\}_0\}_{k}\}_{l}+
  \{z_{kl}\cdot\{w_{jkl}\cdot\{v_{ijkl}\}_{i}\}_{j}\}_{kl})$};
 \node[right] at (\xc6+62,\yc6+6) {$\scriptstyle
  \cdot (\{x\cdot\{y_{l}\}_{l}\}_1+\{w_{jkl}\cdot\{v_{ijkl}\}_{i}
  \}_{jkl})$};
 \node at (\xc3+95,\ya3-4) {$\scriptstyle l$};
 \node at (\xc1+80,\ya6-4) {$\scriptstyle k$};
 \node at (\xc0+65,\yb3-4) {$\scriptstyle j$};
 \node at (\xb6+10,\yb7-4) {$\scriptstyle i$};
\end{tikzpicture}
$$
}
\caption{An example of a generic form}\label{smh13}
\end{figure}
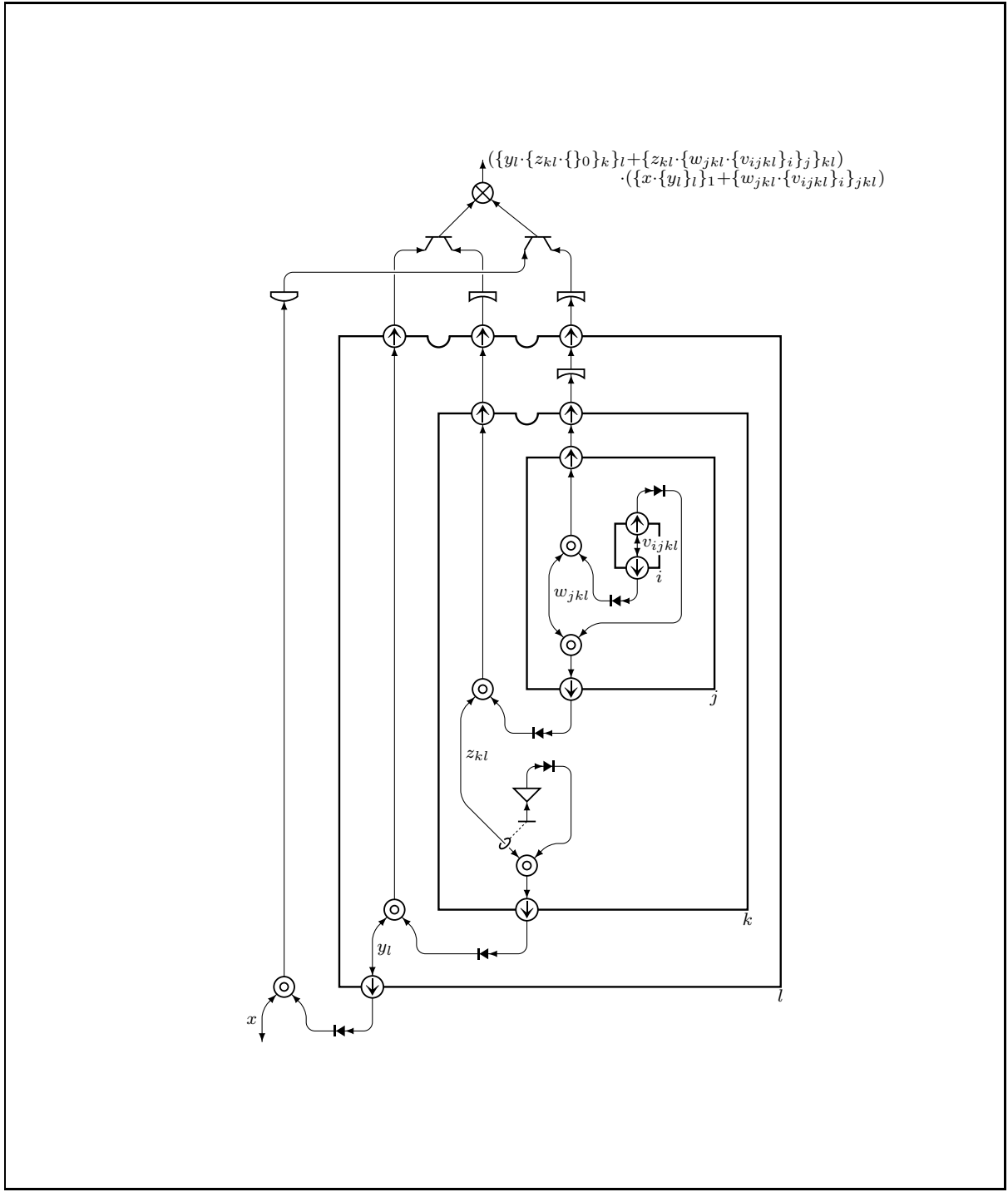

\noindent
Therefore, flows spring at the ridges of
 bioriented wires or unit parts, run down
 into the ocean through the outgoing wires, occasionally joining with one another.
If the generic forms assigned to the outgoing wires are
$$
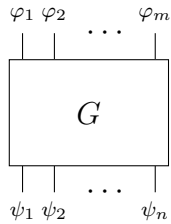
\begin{tikzpicture}[xscale=0.0352778, yscale=0.0352778, thin, inner sep=0]
  \def\p#1#2{%
    \ifcase #1
      \or \ifx#2x  5       \else  0 \fi
      \or \ifx#2x  \x1+12     \else  \y1 \fi
      \or \ifx#2x  \x1+50     \else  \y1 \fi
    \fi}
  \def\x#1{\p#1x}
  \def\y#1{\p#1y}
 \path[use as bounding box] (0,-20) rectangle (60,60);
 \draw (\x1,\y1-10) -- (\x1,\y1+50);
 \draw (\x2,\y2-10) -- (\x2,\y2+50);
 \draw (\x3,\y3-10) -- (\x3,\y3+50);
 \filldraw[fill=white] (0,0) rectangle (60,40);
 \node at (30,20) {$G$};
 \node at ($ .5*(\x3,\y3-10)+.5*(\x2,\y2-10) $) {$\cdots$};
 \node at ($ .5*(\x3,\y3+50)+.5*(\x2,\y2+50) $) {$\cdots$};
 \node at (\x1,\y1+57) {$\scriptstyle \varphi_1$};
 \node at (\x2,\y2+57) {$\scriptstyle \varphi_2$};
 \node at (\x3,\y3+57) {$\scriptstyle \varphi_m$};
 \node at (\x1,\y1-17) {$\scriptstyle \psi_1$};
 \node at (\x2,\y2-17) {$\scriptstyle \psi_2$};
 \node at (\x3,\y3-17) {$\scriptstyle \psi_n$};
\end{tikzpicture}
$$
 then the {\it generic form} $\varphi(G)$ is
 defined to be $(\varphi;\psi)$ where
 $\varphi=(\varphi_1,\varphi_2,\ldots,\varphi_m)$ and
 $\psi=(\psi_1,\psi_2,\ldots,\psi_n)$.
Figure~\ref{smh13} gives the generic form of the graph in Fig.~\ref{yes66}.

	\vskip2ex

\begin{lemma}\label{vqa10}
From the generic form $(\varphi;\psi)$ together
 with the types of the outgoing wires,
 we can restore the graph $G$ up to $\sim$.
\end{lemma}

\proof
We decompose the generic form one by one, restoring the graph.
Although the duality is not reflected by the generic form,
 its information is covered by the types.
For example, when $\varphi\cdot\psi$ is decomposed, the type
 of wires determines whether it comes from tensor or cotensor.
No information of the dotted lines of units/counits is reflected
 in generic forms, thus the restoration is up to $\sim$.
\endproof
\afterpage{\clearpage}

	\vskip2ex

\noindent
As seen in the lemma above, a generic form almost
 reflects the structure of a graph.
Simultaneously, the generic form behaves as a template
 the instances of which give information of the interpretation
 in $\mathscr{M}$.
An instance is provided by an assignment pair, which
 we explain shortly.

As $\mathop{\rm exp}A$ occurs
 in the interpretation of boards, and as
 boards can be nested, we have to manipulate
 iterative applications of the functor $\mathop{\rm exp}A$.
For instance, an element of $\mathop{\rm exp}(\mathop{\rm exp}A)$
 is a multiset of multisets, e.g., $\{\{a,b\},\{c,d,e\}\}$.
This is represented by the tree
$$
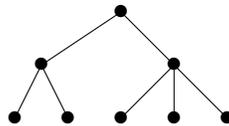
\begin{tikzpicture}[xscale=0.0352778, yscale=0.0352778, thin, inner sep=0]
  \def\pa#1#2{%
    \ifcase #1
          \ifx#2x  0           \else  0 \fi
      \or \ifx#2x  \xa0+20     \else  \ya0 \fi
      \or \ifx#2x  \xa0+40     \else  \ya0 \fi
      \or \ifx#2x  \xa0+60     \else  \ya0 \fi
      \or \ifx#2x  \xa0+80     \else  \ya0 \fi
      \or \ifx#2x  \xa0+10     \else  \ya0+20 \fi
      \or \ifx#2x  \xa0+60     \else  \ya0+20 \fi
      \or \ifx#2x  \xa0+40     \else  \ya0+40 \fi
    \fi}
  \def\xa#1{\pa#1x}
  \def\ya#1{\pa#1y}
  \draw (\xa7,\ya7) -- (\xa6,\ya6);
  \draw (\xa7,\ya7) -- (\xa5,\ya5);
  \draw (\xa6,\ya6) -- (\xa4,\ya4);
  \draw (\xa6,\ya6) -- (\xa3,\ya3);
  \draw (\xa6,\ya6) -- (\xa2,\ya2);
  \draw (\xa5,\ya5) -- (\xa1,\ya1);
  \draw (\xa5,\ya5) -- (\xa0,\ya0);
  \node at (\xa0,\ya0) {$\bullet$};
  \node at (\xa1,\ya1) {$\bullet$};
  \node at (\xa2,\ya2) {$\bullet$};
  \node at (\xa3,\ya3) {$\bullet$};
  \node at (\xa4,\ya4) {$\bullet$};
  \node at (\xa5,\ya5) {$\bullet$};
  \node at (\xa6,\ya6) {$\bullet$};
  \node at (\xa7,\ya7) {$\bullet$};
\end{tikzpicture}
$$
 where the leaves are annotated by the elements $a$ through $e$
 in this order.
We assign a sequence of
 non-negative integers to each board.
Suppose that $i$ is the identifier of the board.
The sequence is denoted by $\langle m_i(s)\,;\>s=1,2,\ldots,q_i\rangle$.
The length $q_i$ is determined recusrively
 as follows.
If $i$ is an outermost board, we set $q_i=1$.
If $i$ is the innermost board containing $j$,
 we set $q_j=\sum_{s=1}^{q_i}m_i(s)$.
The tree above is the case of $m_i(1)=2$ and
 $m_j(s)=2,3$ for $s=1,2$ respectively, for
 a graph that has two nested boards $i$ and $j$, where $i$ is outer,
The number $s$ eneumerates the nodes at a fixed height.
We identify $s$ with a finite sequence
 $(r_1,r_2,\ldots,r_n)$ of positive integers.
For instance, in the following tree
$$
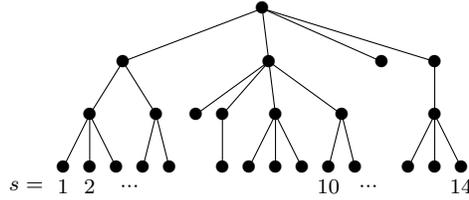
\begin{tikzpicture}[xscale=0.0352778, yscale=0.0352778, thin, inner sep=0]
  \def\pa#1#2{%
    \ifcase #1
          \ifx#2x  0           \else  0 \fi
      \or \ifx#2x  \xa0+10     \else  \ya0 \fi
      \or \ifx#2x  \xa0+20     \else  \ya0 \fi
      \or \ifx#2x  \xa0+30     \else  \ya0 \fi
      \or \ifx#2x  \xa0+40     \else  \ya0 \fi
      \or \ifx#2x  \xa0+50     \else  \ya0 \fi
      \or \ifx#2x  \xa0+60     \else  \ya0 \fi
      \or \ifx#2x  \xa0+70     \else  \ya0 \fi
      \or \ifx#2x  \xa0+80     \else  \ya0 \fi
      \or \ifx#2x  \xa0+90     \else  \ya0 \fi
      \or \ifx#2x  \xa0+100    \else  \ya0 \fi
    \fi}
  \def\xa#1{\pa#1x}
  \def\ya#1{\pa#1y}
  \def\pb#1#2{%
    \ifcase #1
          \ifx#2x  \xa0+100    \else  \ya0 \fi
      \or \ifx#2x  \xa0+110    \else  \ya0 \fi
      \or \ifx#2x  \xa0+120    \else  \ya0 \fi
      \or \ifx#2x  \xa0+130    \else  \ya0 \fi
      \or \ifx#2x  \xa0+140    \else  \ya0 \fi
      \or \ifx#2x  \xa0+150    \else  \ya0 \fi
      \or \ifx#2x  \xa0+10     \else  \ya0+20 \fi
      \or \ifx#2x  \xa0+35     \else  \ya0+20 \fi
      \or \ifx#2x  \xa0+50     \else  \ya0+20 \fi
      \or \ifx#2x  \xa0+60     \else  \ya0+20 \fi
    \fi}
  \def\xb#1{\pb#1x}
  \def\yb#1{\pb#1y}
  \def\pc#1#2{%
    \ifcase #1
          \ifx#2x  \xa0+80     \else  \ya0+20 \fi
      \or \ifx#2x  \xa0+105    \else  \ya0+20 \fi
      \or \ifx#2x  \xa0+140    \else  \ya0+20 \fi
      \or \ifx#2x  \xa0+22.5   \else  \ya0+40 \fi
      \or \ifx#2x  \xa0+77.5   \else  \ya0+40 \fi
      \or \ifx#2x  \xa0+120    \else  \ya0+40 \fi
      \or \ifx#2x  \xa0+140    \else  \ya0+40 \fi
      \or \ifx#2x  \xa0+75     \else  \ya0+60 \fi
    \fi}
  \def\xc#1{\pc#1x}
  \def\yc#1{\pc#1y}
  \draw (\xc7,\yc7) -- (\xc6,\yc6);
  \draw (\xc7,\yc7) -- (\xc5,\yc5);
  \draw (\xc7,\yc7) -- (\xc4,\yc4);
  \draw (\xc7,\yc7) -- (\xc3,\yc3);
  \draw (\xc6,\yc6) -- (\xc2,\yc2);
  \draw (\xc4,\yc4) -- (\xc1,\yc1);
  \draw (\xc4,\yc4) -- (\xc0,\yc0);
  \draw (\xc4,\yc4) -- (\xb9,\yb9);
  \draw (\xc4,\yc4) -- (\xb8,\yb8);
  \draw (\xc3,\yc3) -- (\xb7,\yb7);
  \draw (\xc3,\yc3) -- (\xb6,\yb6);
  \draw (\xc2,\yc2) -- (\xb5,\yb5);
  \draw (\xc2,\yc2) -- (\xb4,\yb4);
  \draw (\xc2,\yc2) -- (\xb3,\yb3);
  \draw (\xc1,\yc1) -- (\xb1,\yb1);
  \draw (\xc1,\yc1) -- (\xb0,\yb0);
  \draw (\xc0,\yc0) -- (\xa9,\ya9);
  \draw (\xc0,\yc0) -- (\xa8,\ya8);
  \draw (\xc0,\yc0) -- (\xa7,\ya7);
  \draw (\xb9,\yb9) -- (\xa6,\ya6);
  \draw (\xb7,\yb7) -- (\xa4,\ya4);
  \draw (\xb7,\yb7) -- (\xa3,\ya3);
  \draw (\xb6,\yb6) -- (\xa2,\ya2);
  \draw (\xb6,\yb6) -- (\xa1,\ya1);
  \draw (\xb6,\yb6) -- (\xa0,\ya0);
  \node at (\xa0,\ya0) {$\bullet$};
  \node at (\xa1,\ya1) {$\bullet$};
  \node at (\xa2,\ya2) {$\bullet$};
  \node at (\xa3,\ya3) {$\bullet$};
  \node at (\xa4,\ya4) {$\bullet$};
  \node at (\xa6,\ya6) {$\bullet$};
  \node at (\xa7,\ya7) {$\bullet$};
  \node at (\xa8,\ya8) {$\bullet$};
  \node at (\xa9,\ya9) {$\bullet$};
  \node at (\xb0,\yb0) {$\bullet$};
  \node at (\xb1,\yb1) {$\bullet$};
  \node at (\xb3,\yb3) {$\bullet$};
  \node at (\xb4,\yb4) {$\bullet$};
  \node at (\xb5,\yb5) {$\bullet$};
  \node at (\xb6,\yb6) {$\bullet$};
  \node at (\xb7,\yb7) {$\bullet$};
  \node at (\xb8,\yb8) {$\bullet$};
  \node at (\xb9,\yb9) {$\bullet$};
  \node at (\xc0,\yc0) {$\bullet$};
  \node at (\xc1,\yc1) {$\bullet$};
  \node at (\xc2,\yc2) {$\bullet$};
  \node at (\xc3,\yc3) {$\bullet$};
  \node at (\xc4,\yc4) {$\bullet$};
  \node at (\xc5,\yc5) {$\bullet$};
  \node at (\xc6,\yc6) {$\bullet$};
  \node at (\xc7,\yc7) {$\bullet$};
  \node at (\xa0-14,\ya0-7) {$\scriptstyle s\;={}$};
  \node at (\xa0,\ya0-7) {$\scriptstyle 1$};
  \node at (\xa1,\ya0-7) {$\scriptstyle 2$};
  \node at (\xa2+5,\ya0-7) {$\scriptstyle \cdots$};
  \node at (\xb0,\ya0-7) {$\scriptstyle 10$};
  \node at (\xb1+5,\ya0-7) {$\scriptstyle \cdots$};
  \node at (\xb5,\ya0-7) {$\scriptstyle 14$};
\end{tikzpicture}
$$
 the leaf $s=10$ is identified with $(2,4,1)$ since
 it is reached from the root by taking the second, the fourth, and the
 first children successively.
At the root level, $s=1$ is identified with an empty sequence $()$.
In general, suppose that $i_n,\ldots,i_2,i_1$ are the boards nesting from
 the outermost $i_n$ in this order.
For the sequence $(r_1,r_2,\ldots,r_n)$, the components range over
$$\vbox{\halign{$1\ \leq\ #$\hfil &${}\ \leq\ #$\hfil\cr
 r_1 & m_{i_n}()\cr 
 r_2 & m_{i_{n-1}}(r_1)\cr 
 r_3 & m_{i_{n-2}}(r_1,r_2)\cr 
 \omit\hfil $\vdots$.\cr}}$$
We note that $s$ enumerates the nodes of a fixed level.
If levels are different, the same number $s$ corresponds to
 sequences of different lengths.

Let $x_{i_1i_2\cdots i_n}$ be a variable for a wire of
 an atomic type $A$.
We recall that $i_1,i_2,\ldots,i_n$ are the boards containing
 the wire, listed from the innermost.
We consider an assignment $\eta(x_{i_1i_2\cdots i_n},s)$.
It assigns an element of the interpreation of $A$
 for each of $s=1,2,\ldots,q_{i_1}$.

An {\it assignment pair} is defined as $P=(\{m_i\}_i,\eta)$ where
 $i$ ranges over all boards of a graph.
Given $s=(r_1,r_2,\ldots,r_k)$,
 an assignment pair $P(s)$ is naturally induced.
It is the pair $(\{m'_{j}\}_{j},\eta')$
 where $m'_j(s)=m_j(s\mathbin{\widehat{\ }}s')$
 and $\eta'(x_{i_1i_2\cdots i_n},s')=\eta
(x_{i_1i_2\cdots i_n},s\mathbin{\widehat{\ }}s')$.
Here $s\mathbin{\widehat{\ }}s'$ denotes
 the concatenation of sequences.

We associate $|\varphi|_P$ with each generic
 form $\varphi$ and each assignment pair $P$.
If $\varphi$ annotates a wire of type $A$,
 then $|\varphi|_P$ is an element of the interpretation
 of the type $A$.
We define $|\varphi|_{P(s)}$ recursively.

\begingroup
	\vskip.5ex
        \hangafter0\hangindent2em
        \noindent
\llapem2{(i)}%
$|x_{i_1i_2\cdots i_n}|_{P(s)}=\eta(x_{i_1i_2\cdots i_n},s)$.

	\vskip0ex
        \noindent
\llapem2{(ii)}%
$|\mathord*|_{P(s)}$ is the unique element $*$ of ${\bf 1}$.

	\vskip0ex
        \noindent
\llapem2{(iii)}%
$|\varphi\cdot\psi|_{P(s)}$ is the pair
 of $|\varphi|_{P(s)}$ and $|\psi|_{P(s)}$.

	\vskip0ex
        \noindent
\llapem2{(iv)}%
$|\varphi+\psi|_{P(s)}$ is the multiset union
 of $|\varphi|_{P(s)}$ and $|\psi|_{P(s)}$.

	\vskip0ex
        \noindent
\llapem2{(v)}%
$|\{\}_0|_{P(s)}$ is the empty multiset $\emptyset$.

	\vskip0ex
        \noindent
\llapem2{(vi)}%
$|\{\varphi\}_1|_{P(s)}$ is the singleton the member
 of which is $|\varphi|_{P(s)}$.

	\vskip0ex
        \noindent
\llapem2{(vii)}%
$|\{\varphi\}_{i_1i_2\cdots i_l}|_{P(s)}$ is the multiset
 consisting of $|\varphi|_{P(s\mathbin{\widehat{\ }}s')}$
 where $s'=(r'_1,r'_2,\ldots,r'_l)$ ranges over the sequences
 in which the domain of $r'_j$ is determined by $\{m_i\}_i$.

	\vskip.5ex
        \hangafter0\hangindent0pt
\endgroup
        \noindent
Suppose that $\varphi(G)=(\varphi;\psi)$
 is the generic form of the graph $G$.
Then $(\alpha;\beta)$
 is an {\it intstance} of $\varphi(G)$
 if there is an assignment pair $P$ such
 that $\alpha=|\varphi|_P$ and $\beta=|\psi|_P$.

	\vskip2ex

\begin{remark}
A simple but important observation is the invariance
 of positive multiset occurrences.
The unions of multisets occur only in
 (iv) and (vii) above, that is,
 at the upper ends of a $d$-part and a $\delta$-part.
They are negative multisets.
In fact, the flows traverse upward in these parts.
Each positive multiset is created at the positive gate
 of a board.
It is never reformed afterward, thus appears
 exactly in the same shape at an outermost wire.
\end{remark}

	\vskip2ex

\begin{lemma}\label{tiz23}
Suppose that $G$ is the graph obtained from a morphism $f$ in
 normal form.
Then $M_f[\alpha;\beta]\neq\emptyset$
 if and only if $(\alpha;\beta)$ is an instance of $\varphi(G)$.
\end{lemma}

\proof
We consider the process
 $\pi(G)^{\alpha}_{\beta}$.
We show that the process succeeds if and only if
$(\alpha;\beta)$ is an instance.
Suppose that the process succeeds.
At a board, concurrent subprocesses are
 invoked.
We provide process numbers $s$ by appropriately ordering
 the subprocesses.
If the $s$-th subprocess finds another box $i$,
 the  multisets for its
 gates must have the same cardinality, say $n$.
We assign $m_i(s)=n$.
If the $s$-th process reaches a bioriented wire
 of $x_{i_1i_2\cdots i_k}$, its both ends must be
 annotated by the same element, say $a$.
We assign $\eta(x_{i_1i_2\cdots i_k},s)=a$.
Conversely, suppose that the given tuple is an instance.
Successful non-deterministic choices at $d$-parts, and $\delta$-parts
 and boards in the process are guided
 by the instance.
For example, if an element above a $\delta$-part
 is $\gamma=|\{\varphi\}_{i_1i_2\cdots i_n}|_{P(s)}$,
 it is decomposed into the multiset of multisets, $\{\gamma_1,\gamma_2,\ldots,
 \gamma_{m_{i_n}(s)}\}$, where $\gamma_j$ is given
 as $|\{\varphi\}_{i_1i_2\cdots i_{n-1}}|_{P(s\cdot j)}$.
Here $s\cdot j$ denotes the sequence obtained
 by adjoining $j$ to the right of the sequence $s$.
\endproof

	\vskip2ex

\begin{remark}
We can also use the generic forms to count the number 
$M_f[\alpha;\beta]$, not
 only to know whether it is non-zero.
To this end, we must modify the definitions
 of generic forms and instances
 to cope with the matter of asymmetry discussed in \S\ref{ypr76}.
In instances, we replace each positive
 multiset occurrence $\{\gamma_1,\gamma_2,\ldots,
\gamma_n\}$ by a list $[\gamma_1,\gamma_2,\ldots,
\gamma_n]$, the order of elements in which is fixed
 arbitrarily.
Negative multisets are not altered.
We let $\alpha'$ and $\beta'$ denote
 the result of such modification.
Accordingly, we replace the form $\{\varphi\}_i$ at a positive
 position with $[\varphi]_i$.
The instance $|[\varphi]_i|_{P}$ is redefined
 as a list, instead of a multiset.
Then $M_f[\alpha;\beta]$ is equal to
 the number of the assignment pairs that yield
 $(\alpha';\beta')$ as instances.
This method is not
 suitable to actually compute the matrix components,
 since we must screen the matching instances
 after generating all possible instances.
We prefer the method discussed in \S\ref{yhv89}.
In the proof of our main theorem, we are concerned
 only with the instances in which the difference between multisets and lists
 does not matter.
See Def.~\ref{ubo83}.
\end{remark}

\section{Confluence}

Confluence is the objective of this paper.
It is verified via the model $\mathscr{M}$.
We focus on special instances, called
 $p$-echo instances.
We show that the generic forms can be uniquely rebuilt
 from such special instances.
Suppose that a morphism has two normal forms, $f$ and $f'$,
 the graphs of which are denoted by $G$ and $G'$ respectively.
Then, the generic form $\varphi(G)$ must equal $\varphi(G')$, provided that
 they share a common $p$-echo instance.
Then $G\sim G'$ is concluded
 by Lem.~\ref{vqa10}.
Namely, the normal forms $f$ and $f'$ must be almost equal.
This establishes the confluence up to $\sim$.
A non-trivial part of this proof strategy is
 how to assure the existence of a common $p$-echo instance.
This is achieved via the elementary number theory
 applied to the interpretation in $\mathscr{M}$.

	\vskip2ex

\begin{definition}\label{ubo83}\rm
Let $p$ be a prime number.
A {\it $p$-echo} assignment pair
 of the generic form $\varphi(G)$ is defined
 as $P=(\{m_i\}_i,\eta)$ subject to the following conditions:

	\vskip.5ex
        \hangafter0\hangindent2em
        \noindent
\llapem2{(i)}%
The value of $m_i(s)$ does not depend on $s$,
 and $m_i(s)=p^{p^k}$ for some $k\geq 0$.
Moreover, if $p^{p^k}$ and $p^{p^l}$ are the values for
 distinct boards $i\neq j$, then $k\neq l$.

	\vskip0ex
        \noindent
\llapem2{(ii)}%
The assignment $\eta(x_{i_1i_2\cdots i_n},s)$
 does not depend on $s$.
Moreover, for distinct variables $x_{i_1i_2\cdots i_n}\neq
 y_{j_1j_2\cdots j_m}$, the assigned elements are different:
 $\eta(x_{i_1i_2\cdots i_n},s)\neq \eta(y_{j_1j_2\cdots j_m},s')$.

	\vskip.5ex
        \hangafter0\hangindent2em
        \noindent
An instance is {\it $p$-echo} if it is an instance produced by
 a $p$-echo assignment pair.
\end{definition}

	\vskip2ex

\noindent
The first halves of the conditions are called the {\it uniformity}
 condition:
 neither $m_i(s)=p^{p^k}$ nor $\eta(x_{i_1i_2\cdots i_n},s)$ does
 depend on $s$.
The second halves are called the {\it discernibility} condition:
$m_i(s)\neq m_j(s')$ for different boards
 and $\eta(x_{i_1i_2\cdots i_n},s)\neq \eta(y_{j_1j_2\cdots j_m},s')$
 for different variables.
A $p$-echo pair assigns the same number to the same board,
 and the same element to the same variable.
It assigns different things to different
 boards or different variables.
In a word, a $p$-echo assignment gives
 an ``echoing'' repercussion of
 the graph $G$.

	\vskip2ex

\begin{definition}\rm
A multiset is {\it homogeneous} if it
 is of the shape $\{\alpha,\alpha,\ldots,\alpha\}$.
We write $n\{\alpha\}$ regarding it as
 $\{\alpha\}+\{\alpha\}+\cdots +\{\alpha\}$
 ($n$ is the number of elements).
\end{definition}

	\vskip2ex

\noindent
The positive multisets occurring in a $p$-echo instance
 are all homogeneous.
We recall that positive multisets are never modified
 from the moments of creation at the positive gates of boards.
By the uniformity condition, the positive multisets are homogeneous
 at the positive gates.
Thus it remains so at the outermost wires.
The set of different positive multisets in the $p$-echo instance
 has a one-to-one correspondence to
 the set of boards in the graph.
It is an immediate consequence of the discernibility condition.
Different positive multisets have different numbers of elements.

We give an example of a $p$-echo instance of the generic form
 in Fig.~\ref{smh13}.
We assign $m_i(s)=p^{p^{k_1}},m_j(s)=p^{p^{k_2}},
m_k(s)=p^{p^{k_3}},m_l(s)=p^{p^{k_4}}$ from the innermost boards.
For the variables, we assign
$\eta(x,s)=a,\eta(y_{l},s)=b,
\eta(z_{kl},s)=c,\eta(w_{jkl},s)=d,\eta(v_{ijkl},s)=e$.
Then the $p$-echo instance
 is the pair of

	\vskip2ex
        \noindent\kern5em
$(p^{p^{k_1}}\{(b,\overline{p^{p^{k_2}}\{(c,\overline{\emptyset})\})}\}+
p^{p^{k_1}+p^{k_2}}\{(c,\overline{p^{p^{k_3}}\{(d,\overline{p^{p^{k_4}}\{e\})}\})}\},$
	\nopagebreak\par\nopagebreak\noindent\hfill
$p^0\{(a,\overline{p^{p^{k_1}}\{b\})}\}
 +p^{p^{k_1}+p^{k_2}+p^{k_3}}\{(d,\overline{p^{p^{k_4}}\{e\})}\})$

	\vskip2ex
        \noindent
 and $a$.
We add overbars to clarify the distinction between
 positive and negative.
The first component of the pair corresponds to the source
 of a morphism, thus regarded to be contravariant.
So its negative occurrences are positive.
The list of different positive multisets is
 $p^{p^{k_1}}\{b\}$, $p^{p^{k_2}}\{(c,\emptyset)\}$,
 $p^{p^{k_3}}\{(d,p^{p^{k_4}}\{e\})\}$, and $p^{p^{k_4}}\{e\}$.
These have a one-to-one correspondence to the boards.
Moreover, all are homogeneous.
We observe that the number of occurrences of these
 multisets are $p^0,\,p^{p^{k_1}},\,p^{p^{k_1}+p^{k_2}},\,
p^{p^{k_1}+p^{k_2}+p^{k_3}}$ respectively.

	\vskip2ex

\begin{definition}\rm
The {\it size} of a graph $G$ is defined
 as the number of wires.
It is denoted by $\mathop{\rm size}(G)$.
\end{definition}

	\vskip2ex

\noindent
For the following argument, indeed, a tighter bound works.
It suffices that $\mathop{\rm size}(G)$ is not less
 than the number of legs of any multi-duplicators,
 the number of boards, and the number of bioriented wires.
We choose this definition as is simple.

A $p$-echo instance directly reflects the structure of
 the underlying graph.
As naturally expected, thus, 
 the generic form can be restored from such an instance.
We show that this is true if the graph is not too big.

	\vskip2ex

\begin{lemma}\label{pwe08}
We consider the class of graphs $G$ satisfying that
 ${\rm size}(G)<p$.
If $(\alpha;\beta)$ is a $p$-echo instance of some
 generic form, then the form $\varphi(G)$ is uniquely
 determined up to the renaming of the identifies of boards
 and the variables.
\end{lemma}

\proof
Let $\delta_1,\delta_2,\ldots,\delta_m$ enumerate all
 different positive multisets in $(\alpha;\beta)$.
As explained above, they have one-to-one correspondence
 to the boards $i_1,i_2,\ldots,i_m$.
We put $m_{i_j}(s)=\smash{p^{p^{k_j}}}$.
The information of the nesting among boards is
 recovered from the number of occurrences of each $\delta_j$.
Let us suppose that $i_{u_1},i_{u_2},i_{u_3},\ldots$ nests
 from the outermost in this order.
Then, $\delta_{u_1}$ occurs exactly once, $\delta_{u_2}$
 occurs $p^{p^{k_{u_1}}}$ times, $\delta_{u_2}$ occurs
 $p^{p^{k_{u_1}}+p^{k_{u_2}}}$ times, and so on.
Namely, divisibility exactly corresponds
 to the nesting.
There is only one way to rewrite $p^l$
 into $p^{p^{k_{u_1}}+p^{k_{u_2}}+\cdots +p^{k_{u_t}}}$ except
 the order of the sum, since all $k_j$ are differenct
 by the discernibility condition.
By comparing the numbers by the divisibility,
 therefore, we can retain the complete information
 of nesting.
Here the hypothesis about the size is not needed.

We reconstruct the generic form from the
 instance, one by one, from
 the outside.
We replace a positive multiset $p^{p^{k_j}}\{\gamma'\}$ with
 $\{\gamma'\}_j$ and we continue from $\gamma'$.
For a negative multiset $\gamma$, we first decompose
 it into $\gamma=\gamma_1+\gamma_2+\cdots+\gamma_n$.
Here, $n<p$ and each $\gamma_i$ is a homogeneous multiset
 of $p^{l_i}$ elements for some $l_i\geq 0$.
As $n$ turns out to be the number of legs of a multi-duplicator,
 this constraint is justified as $n\leq\mathop{\rm size}(G)<p$.
We do not exclude the case $\gamma_i=\gamma_j$.
For example, if $p=3$ and $\gamma=\{\alpha,\alpha,\alpha,\alpha,\alpha,
 \alpha\}$, it is decomposed into
 $\gamma_1=\{\alpha,\alpha,\alpha\}=\gamma_2$.
It is disallowed to decompose $\gamma$ into the sum of six $\{\alpha\}$
 since $n<p$.
The decomposition of $\gamma$ is, thus, unique up
 to the shuffles of the sum.
We rewrite the number $p^{l_i}$ to
 $p^{p^{k_{u_1}}+p^{k_{u_2}}+\cdots +p^{k_{u_t}}}$.
Provided that the boards $i_{u_1},i_{u_2},\ldots,i_{u_t}$ nest
 from the outside in this order,
 we transform $\gamma_i=p^{l_i}\{\gamma'\}$ to
 $\{\gamma'\}_{i_{u_t}\cdots i_{u_2}i_{u_1}}$.
If $l_i=0$, we transform it to $\{\gamma'\}_1$.
Then we continue from $\gamma'$.
If $\gamma=\emptyset$, we transform it to
 $\{\}_0$.
We transform $(\alpha,\beta)$ to $\alpha\cdot\beta$
 and continue from each of $\alpha$ and $\beta$.
Finally, if we reach an atom, we translate it
 to a variable with appropriate indices.
We choose different variables for different atoms.
\endproof

	\vskip2ex

\noindent
This is not a deep result.
If we know that an instance is $p$-echo from the outset,
 we can rebuild the generic form.
What really matters is how to judge if
 a given instance is $p$-echo.
Before discussing this problem, we give
 an immediate consequence of the lemma.

	\vskip2ex

\begin{corollary}\label{pbm43}
Suppose that $\varphi(G)$ and $\varphi(G')$ shares
 a common $p$-echo
 instance.
If $\mathop{\rm size}(G),\mathop{\rm size}(G')<p$, then
 $G\sim G'$ holds.
\end{corollary}

\proof
By Lem.~\ref{pwe08} and \ref{vqa10}.
\endproof

	\vskip2ex

\noindent
We want to verify that, if $[\![f]\!]=[\![g]\!]$ in the
 model $\mathscr{M}$, then $\varphi(G_f)$ and $\varphi(G_g)$ shares
 a $p$-echo instance.
Here $G_f$ denotes the graph associated with $f$.
If $f$ and $g$ are normal forms of
 a common morphism, $[\![f]\!]=[\![g]\!]$ is
 true.
So we have $G_f\sim G_g$ by Cor.~\ref{pbm43}.
Namely $f\sim g$ is concluded.
This establishes confluence up to $\sim$.

To this end, we give a characterization of $p$-echo
 instances with no referrence to the assignment pairs.
We consider the following five conditions for a
 matrix $M=M_f$ in ${\bf Set}^{A\times B}$ and
 $(\alpha;\beta)\in A\times B$:

	\vskip.5ex
        \hangafter0\hangindent3em
        \noindent
\llapem3{($\star 1$)}%
$M[\alpha;\beta]\not\equiv 0\ {\rm mod}\ p$.

	\vskip0ex
        \noindent
\llapem3{($\star 2$)}%
Every positive multisets occurring in $(\alpha;\beta)$ is homogeneous
 and has $p^{p^k}$ elements for some $k\geq 0$.

	\vskip0ex
        \noindent
\llapem3{($\star 3$)}%
If two positive multisets $\gamma,\gamma'$ occurring in $(\alpha;\beta)$
 has the same $p^{p^k}$ elements, then $\gamma=\gamma'$.

	\vskip0ex
        \noindent
\llapem3{($\star 4$)}%
If a multiset occurs positively in $(\alpha;\beta)$,
 then it occurs positively exactly $p^l$ times for
 some $l\geq 0$.

	\vskip0ex
        \noindent
\llapem3{($\star 5$)}%
Each signed atom has exactly $p^l$ occurrences
in $(\alpha;\beta)$ for
 some $l\geq 0$.

	\vskip.5ex
        \hangafter0\hangindent0pt
        \noindent
The first condition is involved in the matrix.
The rest address only the elements $(\alpha;\beta)$.
We show that, for sufficiently large $p$, all
 $p$-echo instances satisfy the five conditions,
 and vice versa.

Whenever we mention a graph $G$ and a matrix $M$ in
 the following, we implicitly assume that
 they come from a morphism $f$ in normal form as $G=G_f$ and
 $M=M_f$.

	\vskip2ex

\begin{lemma}\label{twa45}
A $p$-echo instance $(\alpha;\beta)$ of $\varphi(G)$ satisfies
 $(\star 2)$ through $(\star 5$).
\rm
\end{lemma}

\noindent
Comment: we postpone $(\star 1)$ since it requires a
 sensitive argument and
 $p$ must be taken large.

\proof
Each positive multiset is the one created by the positive gate
 of a board.
Thus it must be homogeneous and contains $p^{p^k}$ elements.
This is the condition $(\star 2)$.
By the discernibility condition, $(\star 3)$ holds.
Moreover, if a board (a bioriented wire) occurs
 in boards $i_{u_1},i_{u_2},\ldots,i_{u_t}$, then
 the number of occurrences of the corresponding positive
 multisets (signed atoms) is of the form
 $p^{p^{k_{u_1}}+p^{k_{u_2}}+\cdots+p^{k_{u_t}}}$.
Hence ($\star 4$) and ($\star 5$) hold.
\endproof

	\vskip2ex

\noindent
Fermat's little theorem asserts that $a^p\equiv a\ {\rm mod}\ p$
 for prime $p$.
The theorem is extended to the following.

	\vskip2ex

\begin{lemma}\label{ihq88}
Let $p$ be a prime number and $l$ a positive integer.

	\vskip.5ex
        \hangafter0\hangindent2em
        \noindent
\llapem2{(i)}%
$a^{p^l}\equiv a\ {\rm mod}\ p$ holds.

	\vskip0ex
        \noindent
\llapem2{(ii)}%
Whenever $j\neq 0,p^l$, ${p^l\choose j}\equiv 0
 \ {\rm mod}\ p$ holds.

	\vskip0ex
\end{lemma}

	\vskip2ex
        \noindent
(i) is obtained by iterating Fermat's little theorem $l$ times.
Contrary to its name, Fermat's little theorem was first verified by
 Euler \cite[p.~63]{hawr}.
He proved (ii) first (for $l=1$), from which he derived the theorem.
Conversely, we can derive (ii) from Fermat's little theorem
 as follows.
We have $1+x^p\equiv 1+x\equiv (1+x)^p\ {\rm mod}\ p$ by the theorem.
Thus $(1+x)^p-(1+x^p)$ is constant $0$ as a function
 of $\mathbb{Z}/p\mathbb{Z}$.
Since the polynomial of degree $p-1$ has $p$ roots, we have
 $(1+x)^p=1+x^p$ as polynomials over
 $\mathbb{Z}/p\mathbb{Z}$.
Iterating twice, we have $(1+x)^{p^2}=(1+x^p)^p=1+x^{p^2}$.
In this way, we obtain $(1+x)^{p^l}=1+x^{p^l}$ by iteration.
Namely, all coefficients ${p^l\choose j}$ vanish to modulus $p$
 except in the constant and the leading term.
The statements (i) and (ii) of the lemma are
 essentially equivalent.

The next lemma is a key in all of the following arguments.

	\vskip2ex

\begin{lemma}\label{hsg59}
Let $G$ be a board.
We consider $\pi(G)^{\alpha_1,\alpha_2,\ldots,\alpha_m}_\beta$, where
 $\beta=\{b,b,\ldots,b\}$ is a homogeneous multiset having
 $p^{l}$ copies of $b$ for some $l\geq 1$.

	\vskip1ex
        \hangafter0\hangindent2em
        \noindent
\llapem2{(i)}%
$\pi(G)^{\alpha_1,\alpha_2,\ldots,\alpha_m}_\beta
\equiv 0\ {\rm mod}\ p$ holds
 unless all $\alpha_i$ are homogeneous multisets having $p^l$ elements.

	\vskip0ex
        \noindent
\llapem2{(ii)}%
$\pi(G)^{\alpha_1,\alpha_2,\ldots,\alpha_m}_\beta
 \equiv \pi(G')^{a_1,a_2,\ldots,a_m}_b
 \ {\rm mod}\ p$ holds
 if all $\alpha_i$ are homogeneous multisets $\{a_i,a_i,\ldots,a_i\}$
 having $p^l$ elements.
Here $G'$ denotes the graph inside the board.
	\vskip0ex
\end{lemma}

	\vskip1ex

\proof
We put $n=p^l$.
By definition, $\pi(G)^{\alpha_1\alpha_2\cdots\alpha_m}_\beta$
 equals $0$ unless all $\alpha_i$ have $n$ elements.
If all have $n$ elements, we have

	\vskip2ex
        \noindent\kern5em
$\displaystyle
\pi(G)^{\alpha_1\alpha_2\cdots\alpha_m}_\beta\ =
\ \sum \pi(G')^{a'_{11}a'_{21}\cdots a'_{m1}}_{b}
\pi(G')^{a'_{12}a'_{22}\cdots a'_{m2}}_{b}\>\cdots\>
\pi(G')^{a'_{1n}a'_{2n}\cdots a'_{mn}}_{b}$.
	\vskip2ex
        \noindent
 since $\beta=\{b,b,\ldots,b\}$.
The summation ranges over different $m\times n$
 matrices
$$\matrix{%
 a'_{11} & a'_{12} & \cdots & a'_{1n} \cr
 a'_{21} & a'_{22} & \cdots & a'_{2n} \cr
 \vdots & \vdots &         & \vdots \cr
 a'_{m1} & a'_{m2} & \cdots & a'_{mn} \cr
}$$
 where the $i$-th row is a linear disposition of
 $\alpha_i$.
If all $\alpha_i$ are homogeneous, there is
 only one matrix
$$\matrix{%
 a_{1} & a_{1} & \cdots & a_{1} \cr
 a_{2} & a_{2} & \cdots & a_{2} \cr
 \vdots & \vdots &         & \vdots \cr
 a_{m} & a_{m} & \cdots & a_{m} \cr
}$$
 thus $\pi(G)^{\alpha_1\alpha_2\cdots\alpha_m}_\beta
 =(\pi(G')^{a_1a_2\cdots a_m}_b)^n$.
Since $n=p^l$, Lem.~\ref{ihq88},(i) implies
 $\pi(G)^{\alpha_1\alpha_2\cdots\alpha_m}_\beta
 \equiv \pi(G')^{a_1a_2\cdots a_m}_b\ {\rm mod}\ p$.
This proves (ii).
If any $\alpha_i$ is not homogeneous,
 the matrix has different column vectors.
We note that a shuffle of column vectors do
 not change the value $\pi(G')^{a'_{11}a'_{21}\cdots a'_{m1}}_{b}
\pi(G')^{a'_{12}a'_{22}\cdots a'_{m2}}_{b}\cdots
\pi(G')^{a'_{1n}a'_{2n}\cdots a'_{mn}}_{b}$.
The number of different matrices obtained by shuffles
 is given by a multinomial coefficient ${n\choose
 k_1,k_2,\ldots,k_q}$ where $q\geq 2$ is the number of different
 column vectors, $k_j\geq 1$, and $k_1+k_2+\cdots +k_q=n$.
The multinomial coefficient is factored by  ${n\choose k_1}$,
 thus it vanishes to modulus $p$ by Lem.~\ref{ihq88}, (ii).
Therefore $\pi(G)^{\alpha_1\alpha_2\cdots\alpha_m}_\beta\equiv
 0\ {\rm mod}\ p$.
This ends the proof of (i).
\endproof

	\vskip2ex

\noindent
Under the hypothesis that all positive multisets are
 homogeneous,
 Lem.~\ref{hsg59} signifies
 that, as long as we count $\pi(G)$
 to modulus $p$, we can ignore non-homogeneous
 instances on a board $G$.
Moreover, if they are all homogeneous, it suffices to
 consider the single subprocess $\pi(G')^{a_1a_2\cdots a_m}_b$.

	\vskip2ex

\begin{lemma}\label{dtd53}
If $M$ and $(\alpha;\beta)$ satisfy $(\star 1)$ and $(\star 2)$,
 then $(\alpha;\beta)$ is the instance of $\varphi(G)$ yielded by
 an assignment pair subject to the uniformity condition.
\end{lemma}

	\noindent
Comment:
The lemma asserts that there is an assignment pair satisfying
 the uniformity condition.
It does not negate the existence of the assignment pairs
 that breach the condition but yield the same
 $(\alpha;\beta)$.

\proof
Since $M[\alpha;\beta]\neq 0$, the process $\pi(G)^\alpha_\beta$
 succeeds.
Namely, there is at least one non-deterministic branch that succeeds.
We show that a successful branch keeps
 the condition $\pi(G')^{\alpha'}_{\beta'}\not\equiv
 0\ {\rm mod}\ p$, and that we can extract an assignment pair
 satisfying the uniformity condition from the branch.

For the tensor introduction $\pi(G)_{(\alpha,\beta)}=\pi(G_1)_\alpha\cdot
 \pi(G_2)_\beta$, if the left-hand side is not congruent to $0$
 to modulus $p$, then neither of $\pi(G_1)_\alpha$ or
 $\pi(G_2)_\beta$ is congruent to $0$.
The case of cotensor elimination is similar.
For the duplicator case $\pi(G)^\gamma =\sum
 \pi(G')^{\gamma_1,\gamma_2}$, there is decomposition
 $\gamma=\gamma_1+\gamma_2$ such that $\pi(G')^{\gamma_1,\gamma_2}$
 is not congruent to $0$.
The case of a $\delta$-part is similar.
Finally, we consider the case of a board.
By the condition $(\star 2)$, the positive gate is annotated by
 a homogeneous multiset of $p^{p^k}$ elements.
As we count the numbers to modulus $p$, we can assume that
 the negative gates are also associated with homogeneous
 multisets of $p^{p^k}$ elements
 by Lem.~\ref{hsg59}.
Moreover, $\pi(G)^{\alpha_1\alpha_2\cdots\alpha_m}_\beta=
\pi(G')^{a_1a_2\cdots a_m}_b\ {\rm mod}\ p$ implies that
 the right-hand side is not congruent to $0$ to modulus $p$.

We show that there is an assingment
$\eta(x_{i_1i_2\cdots i_n},s)$ that does not depend on
$s=(r_1,r_2,\ldots,r_n)$.
We walk along a flow, departing from
 the bioriented wire marked by $x_{i_1i_2\cdots i_n}$.
We choose either of the two directions.
Eventually, we reach one of the positive and negative
 gates of the board $i_1$.
Since the instance associated with the gate is a homogeneous
 multiset,
 we have an assignment that does not depend on
 $r_n$ (the suffixes are reversed,
 since $i_1,i_2,\ldots$ are from the innermost while
 $r_1,r_2,\ldots$ are from the outermost).
If we proceed further, we reach a gate of the board $i_2$.
Since it is associated with a homogeneous multiset,
 we have an assignment that does not depend on $r_{n-1}$.
Repeating this, we conclude that
 the value is irrelevant of $s$.
The uniformity for $m_i$ is similarly verified.
In this case, we start from the positive gate of the board $i$
 and proceed.
\endproof
 
	\vskip2ex

\begin{lemma}\label{phr52}
Suppose that $\mathop{\rm size}(G)<p$ holds.
If an assignment pair subject to the uniformity condition
 yields the instance $(\alpha;\beta)$ that fulfills $(\star 3)$, $(\star 4)$,
 and $(\star 5)$,
 then the assignment pair
 satisfies the discernibility condition.
\end{lemma}

\proof
By the uniformity condition, $m_i=m_i(s)$ does not depend
 on $s$.
Let
 $\delta_i$ denote the multiset associated with
 the positive gate of the board $i$.
Assume that $m_i$ is equal to another $m_j$.
Then $\delta_i$ equals $\delta_j$ by the condition ($\star 3$).
The number of positive occurrences of each $\delta_j$ is
 of the shape $p^{l}=p^{p^{k_1}+p^{k_2}+\cdots +p^{k_t}}$
 in the instance subject to the uniformity condition,
 as mentioned in the proof of Lem.~\ref{twa45}.
Hence, $\delta_i$ occurs $p^{l_1}+p^{l_2}+\cdots
 +p^{l_n}$ times positively for some $n\geq 2$.
This violates the condition ($\star 4$),
 since $n\leq\mathop{\rm size}(G)<p$ holds
 as $n$ is bounded by the number of boards.
Namely, all $m_i$ are distinct.
Likewise, we have different assignments to
 different variables.
We use that the number of bioriented
 wires is not greater than $\mathop{\rm size}(G)$.
\endproof

	\vskip2ex

\noindent
By two lemmata above, one direction
 of the implication is completed.
Namely, we have
 the following proposition.

\begin{proposition}\label{gxe81}
Suppose that $\mathop{\rm size}(G)<p$ holds.
If $M$ and $(\alpha;\beta)$ satisfy $(\star 1)$
 through $(\star 5)$, then $(\alpha;\beta)$ is
 a $p$-echo instance of $\varphi(G)$.
\end{proposition}

	\vskip2ex

\noindent
What remains is to prove that $p$-echo instances
 satisfy the condition ($\star 1$).
We must take $p$ large.
The size of $G$ is not sufficient.

We introduce the notion of {\it duplication scale}
$d(G)$.
The definition follows the reverse of sequentialization.
We consult the reader to the construction of
 the process $\pi(G)$ in \S\ref{yhv89}, and reuse the symbols therein.
If the graph is decomposed into two graphs $G_1$ and $G_2$ as
 in the case of tensor introduction, we
 set $d(G)=d(G_1)d(G_2)$.
For duplicators, we use the multi-duplicators
 integrating successive duplicators as much as possible.
If $G$ equals
$$
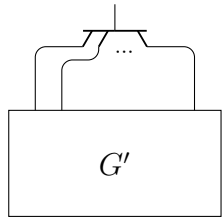
\vcenter{\hbox{%
\begin{tikzpicture}[xscale=0.0352778, yscale=0.0352778, thin, inner sep=0]
  \def\p#1#2{%
    \ifcase #1
      \or \ifx#2x  30         \else  70 \fi
      \or \ifx#2x  -10        \else  0 \fi
      \or \ifx#2x  \x2+10     \else  \y2+40 \fi
      \or \ifx#2x  \x2+20     \else  \y2+40 \fi
      \or \ifx#2x  \x2+70     \else  \y2+40 \fi
    \fi}
  \def\x#1{\p#1x}
  \def\y#1{\p#1y}
 \filldraw[fill=white] (\x2,\y2) rectangle (70,40);
 \node at (30,20) {$G'$};
 \begin{scope}[line width=.7]
  \draw (\x1,\y1) -- ++(12,0);
  \draw (\x1,\y1) -- ++(-12,0);
  \draw (\x1+8.6,\y1) -- ++(3.5,-6.0);
  \draw (\x1-8.6,\y1) -- ++(-3.5,-6.0);
  \draw (\x1-2.6,\y1) -- ++(-3.5,-6.0);
 \end{scope}
 \draw (\x3,\y3)[rounded corners=4] -- (\x3,\y1-6) -- (\x1-12.1,\y1-6);
 \draw (\x4,\y4)[rounded corners=4] -- (\x4,\y1-11) -- (\x1-6.1,\y1-11) -- (\x1-6.1,\y1-6);
 \draw (\x5,\y5)[rounded corners=4] -- (\x5,\y1-6) -- (\x1+12.1,\y1-6);
 \draw (\x1,\y1) -- (\x1,\y1+10);
 \node at (\x1+3.5,\y1-8) {$\scriptstyle\cdots$};
\end{tikzpicture}%
}}%
$$
 we set $d(G)=n\mathord!\cdot d(G')$.
Here $n$ is the number of legs of the multi-duplicator.
For all other parts, we set $d(G)=d(G')$.
In particular, if $G$ is a board containing a subgraph $G'$,
 we have $d(G)=d(G')$.
If we reach a bioriented wire, we set $d(G)=1$.

	\vskip2ex

\begin{lemma}\label{vcl75}
Suppose $\mathop{\rm max}\{\mathop{\rm size}(G),d(G)\}<p$ holds.
If $(\alpha;\beta)$ is a $p$-echo instance of $\varphi(G)$,
 then $M[\alpha;\beta]$ is congruent to one of $1,2,\ldots,d(G)$
 to modulus $p$.
\end{lemma}

\proof
We explore all non-deterministic branches of
 the process $\pi(G)^\alpha_\beta$ to modulus $p$.
We can ignore branches that vanish
 to modulus $p$.
Viewing Lem.~\ref{hsg59}, thus, we can assume that
 the instance at a gate of a board is a homogeneous multiset
 of $p^{p^k}$ elements.
Accordingly, the instance above a $\delta$-part
 is a homogeneous multiset of $p^l$ elements.
We verify that $p$ is so large that no carry-over happens
 when counting in base $p$.
Non-deterministic branches may occur at $d$-parts, $\delta$-parts
 and boards.
By Lem.~\ref{hsg59}, the value of $\pi(G)$ at a board
 is congruent to $\pi(G')$, where $G'$ is the inside of the board,
 to modulus $p$.
At a $\delta$-part, the
 decomposition is uniquely determined
 by divisibility relation, as shown in Lem.\ref{pwe08}.
Only $d$-parts have non-deterministic choices.
The multiset $\gamma$ associated with the upper side of a multi-duplicator
 is decomposed into a sum $\gamma=\gamma_1+\gamma_2+\cdots+\gamma_n$
 of homogeneous multisets where $n$ is the number of legs.
Each $\gamma_i$ has $p^{l_i}$ elements.
Since $n<p$, the collection of $\gamma_1,\gamma_2,\ldots,\gamma_n$
 is unique up to permutation, as discussed in the proof of Lem.~\ref{pwe08}.
For generic forms, the permutation does not matter,
 since $\varphi+\psi=\psi+\varphi$.
For the process $\pi(G)$, however, $\alpha+\beta$ and $\beta+\alpha$
 are separately counted unless $\alpha=\beta$.
Hence the shuffles of $\gamma_1,\gamma_2,\ldots,\gamma_n$ matter.
Not all of the shuffles lead to the successful processes in general.
In the worst case, however, all of $n\mathord!$ shuffles may succeed.
Recall $d(G)=n\mathord!\cdot d(G')$.
Hence, provided that each subprocess
 $\pi(G')$ returns a value up to $d(G')$ to modulus $p$,
 the value of $\pi(G)$ is up to $d(G)$.
Finally, the value of $\pi(G)$ is
 not congruent to $0$, since the $p$-echo
 assignment pair leads to a successful branch.
\endproof

	\vskip2ex
\noindent
The following gives an example where $M[\alpha;\beta]\neq 1$.

$$
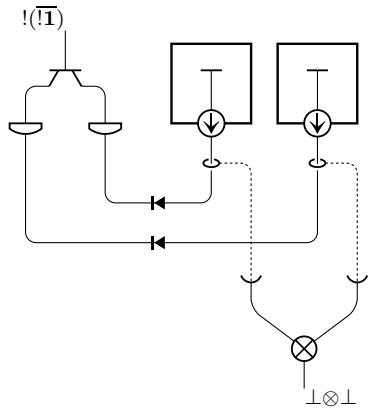
\vcenter{\hbox{%
\begin{tikzpicture}[xscale=0.0352778, yscale=0.0352778, thin, inner sep=0]
  \def\pa#1#2{%
    \ifcase #1
          \ifx#2x  0           \else  0 \fi
      \or \ifx#2x  \xa0-15     \else  \ya0-20 \fi
      \or \ifx#2x  \xa0+15     \else  \ya1 \fi
      \or \ifx#2x  \xa2+40     \else  \ya0 \fi
      \or \ifx#2x  \xa3+40     \else  \ya0 \fi
      \or \ifx#2x  \xa3        \else  \ya1 \fi
      \or \ifx#2x  \xa4        \else  \ya1 \fi
      \or \ifx#2x  \xa2+20     \else  \ya2-30 \fi
      \or \ifx#2x  \xa7        \else  \ya7-15 \fi
      \or \ifx#2x  \xa3+15     \else  \ya8-15 \fi
    \fi}
  \def\xa#1{\pa#1x}
  \def\ya#1{\pa#1y}
  \def\pb#1#2{%
    \ifcase #1
          \ifx#2x  \xa4+15     \else  \ya9 \fi
      \or \ifx#2x  \xa3+35     \else  \yb0-25 \fi
      \or \ifx#2x  \xa5        \else  \ya5-15 \fi
      \or \ifx#2x  \xa6        \else  \yb2 \fi
    \fi}
  \def\xb#1{\pb#1x}
  \def\yb#1{\pb#1y}
 \path[use as bounding box] (\xa0-21,\ya0+24) rectangle (\xb0+4,\yb1-22);
 \draw (\xa0+6,\ya0-6)[rounded corners=4] -- (\xa2,\ya0-6) -- (\xa2,\ya7)
   -- (\xa5,\ya7) -- (\xa5,\ya3);
 \draw (\xa0-6,\ya0-6)[rounded corners=4] -- (\xa1,\ya0-6) -- (\xa1,\ya8)
   -- (\xa4,\ya8) -- (\xa4,\ya4);
 \draw (\xa9,\ya9)[rounded corners=4] -- (\xa9,\ya9-10) -- (\xb1,\yb1);
 \draw (\xb0,\yb0)[rounded corners=4] -- (\xb0,\yb0-10) -- (\xb1,\yb1);
 \draw[dash pattern=on 1 off 1] (\xb2+3,\yb2)[rounded corners=4] -- (\xa9,\yb2)
   -- (\xa9,\ya9);
 \draw[dash pattern=on 1 off 1] (\xb3+3,\yb3)[rounded corners=4] -- (\xb0,\yb3)
   -- (\xb0,\yb0);
 \draw (\xb1,\yb1) -- (\xb1,\yb1-15);
 \draw (\xa0,\ya0) -- (\xa0,\ya0+15);
 \putDuplicator(\xa0,\ya0)
 \putConvex(\xa1,\ya1)
 \putConvex(\xa2,\ya2)
 \boardBoundary(\xa3-15,\ya3+10)(\xa5+15,\ya5)
 \boardBoundary(\xa4-15,\ya3+10)(\xa6+15,\ya5)
 \putPositiveTerminal(\xa3,\ya3)
 \putPositiveTerminal(\xa4,\ya4)
 \putLowerSocket(\xa5,\ya5)
 \putLowerSocket(\xa6,\ya6)
 \putLeftDiode(\xa7,\ya7)
 \putLeftDiode(\xa8,\ya8)
 \putUpperNegativeTerminal(\xa9,\ya9)
 \putUpperNegativeTerminal(\xb0,\yb0)
 \putTensor(\xb1,\yb1)
 \putRubberBand(\xb2,\yb2)
 \putRubberBand(\xb3,\yb3)
 \node[above left] at (\xa0,\ya0+15) {$\scriptstyle \mathord!(\overline{\mathord!{\bf 1}})$};
 \node[below right] at (\xb1,\yb1-15) {$\scriptstyle \bot\otimes\bot$};
\end{tikzpicture}}}
$$
Let $\alpha$ be $\{p^{p^k}\{*\},p^{p^l}\{*\}\}$ with $k\neq l$
 and let $\beta$ be $(*,*)$.
$\alpha$ is a two-point multiset.
At the $d$-part, both $\{p^{p^k}\{*\}\}+\{p^{p^l}\{*\}\}$ and
 $\{p^{p^l}\{*\}\}+\{p^{p^k}\{*\}\}$ succeed.
Therefore $M[\alpha;\beta]$ equal $2\mathord!=2$.
This example generalize to $n\mathord!$ and
 tells why we need the duplication rate $d(G)$.

Now we have the other direction of implications, provided
that $p$ is sufficiently large:

	\vskip2ex

\begin{proposition}\label{nbi47}
Suppose that $\mathop{\rm max}\{\mathop{\rm size}(G),d(G)\}<p$ holds.
If $(\alpha;\beta)$ is a $p$-echo instance of $\varphi(G)$, then
all of the conditions $(\star 1)$ through $(\star 5)$ hold.
\end{proposition}

\proof
$(\star 1)$ is a consequence of Lem.~\ref{vcl75}.
All others come from Lem.~\ref{twa45}.
\endproof

	\vskip2ex

\begin{theorem}
If $f$ and $f'$ are normal forms satisfying $[\![f]\!]
=[\![f']\!]$, then $f\sim f'$ holds.
\end{theorem}

\proof
Let $G$ and $G'$ denote the graphs obtained from $f$ and $f'$.
We take a prime number so
 that $\mathop{\rm max}\{\mathop{\rm size}(G),
\mathop{\rm size}(G'),d(G)\}<p$ holds.
We take a $p$-echo instance $(\alpha;\beta)$
 of $\varphi(G)$.
It satisfies five conditions by
 Prop.~\ref{nbi47}.
Since $[\![f]\!]=[\![f']\!]$ holds, the assocaited matrices $M$
 are the same.
The five conditions are unaltered.
Hence $(\alpha;\beta)$ is a $p$-echo instance of $\varphi(G')$
 by Prop.~\ref{gxe81}.
Therefore $G\sim G'$ by Cor.~\ref{pbm43}.
\endproof

	\vskip2ex

\noindent
So we have the confluence of our system up to the equivalence:

	\vskip2ex

\begin{corollary}
The normal form of a morphism is determined uniquely up
 to $\sim$.
\end{corollary}

	\vskip2ex

\noindent
Let us write $f=g$ if these morphisms are equal in the free linear
 category in the ordinary sense, that is, when
 all defining diagrams are understood to be commutative
 as usual, rather than rewriting.

	\vskip2ex

\begin{corollary}
If $f=g$ in the free classical linear category, 
 their normal forms are almost equal.
\end{corollary}

	\vskip2ex

\noindent
In a previous paper, we have verified a termination
 property \cite{hase2}.
Although it is weak termination, we have a certain strategy
 leading to normal forms definitely.
Together with the result in this paper, we can derive a
 procedure to determine whether given two morphisms in
 the free linear category are equal up to
 the equivalence $\sim$.
Namely, we first transfer them to normal forms,
 and just check the equivalence.
In \cite{bcst}, the existence of a procedure to
 determine the equality between morphisms is regarded as
 the generalization of coherence.
We have shown that the classical linear
 category satisfies the generalized coherence up to
 the equivalence $\sim$.

The current result remains partial.
We have not succeeded in handling the
 isomorphisms related to the units appropriately.
Seemingly, the structure involved in the units
 in the $*$-autonomous category is more intricate
 than one naturally imagines.
Two units, ${\bf 1}$ and $\bot$, interact delicately.
We leave the problem related to the units to
 future work.
The linear normal functor model $\mathscr{M}$ is not
 enough to disentangle the intricacy caused by
 the existence of two distinct units.

\section*{Acknowledgements}
We are grateful to Kazushige Terui to make us notice
 \cite{cato}.
The author is supported by JSPS Kakenhi Grant Number JP15500003.

\begin{references*}

\bibitem{aspe}
A.~Asperti, Linear logic, comonads and optimal reductions,
{\it Fundamenta Informaticae}, 22(1-2):3--22, 1995.

\bibitem{besh}
 U.~Berger and H.~Schwichtenberg,
 \newblock An inverse of the evaluation functional for typed lambda calculus,
 \newblock {\em Proceedings of
 the Sixth Annual Symposium on Logic in Computer Science, LICS '91},
 Amsterdam, The Netherlands, July, 1991, pages 203--211, IEEE, 1991.

\bibitem{bcst}
R.~F.~Blute, J.~R.~B.~Cockett, R.~A.~G.~Seely, and T.~H.~Trimble,
{\rm Natural deduction and coherence for weakly distributive
 categories}, {\it Journal of Pure and Applied Algebra}, 113(3):229--296,
 1996.

\bibitem{carv1}
D.~de~Carvalho, Execution time of $\lambda$-terms via denotational
 semantics and intersection types, {\it Mathematical Structures in
 Computer Science}, 28(7):1169--1203, 2006.

\bibitem{carv2}
D.~de~Carvalho, Taylor expansion in linear logic is invertible,
 {\it Logical Methods in Computer Science}, 14(4:21):1-73, 2018.

\bibitem{cato}
D.~de~Carvalho and L.~Tortora~de~Falco,
 The relational model is injective for Multiplicative Exponential
 Linear Logic (without weakenings), preprint,
 Computing Research Repository - CORR, 2010.

\bibitem{cose}
J.~R.~B.~Cockett and R.~A.~G.~Seely, {\rm Weakly distributive categories},
 {\it Journal of Pure and Applied Algebra}, 114(2):133--173, 1997.

\bibitem{dare}
V.~Danos and L.~Regnier, The structure of multiplicatives,
{\it Archive for Mathematical Logic} 28:181--203. 1989.

\bibitem{gira1}
J.-Y. Girard, The system $F$ of variable
 types, fifteen years later, {\it Theoretical
 Computer Science}, 45:159--192, 1986.

\bibitem{gira2}
J.-Y.~Girard, {\rm Linear logic}, {\it Theoretical Computer Science},
 50(1):1--101, 1987.

\bibitem{gptf}
G.~Guerrieri, L.~Pellissier and L.~Tortora de Falco,
 Computing connected proof(-structure)s from their Taylor
 expansion, In D. Kesner, B. Pientka (eds.), {\it 1st
 International Conference on Formal Structures
 for Computation and Deduction, FSCD 2016}, Vol.~52,
 LIPIcs : Leibniz International Proceedings in Informatics, 2016.

\bibitem{hawr}
G.~H.~Hardy and E.~M.~Wright, {\it An Introduction to the
 Theory of Numbers}, Fifth Edition, Oxford Science Publications, 1979.

\bibitem{hase1}
R.~Hasegawa, Two applications of analytic functors,
{\it Theoretical Computer Science}, 272(1-2):113--175, 2002.

\bibitem{hase2}
R.~Hasegawa, A categorical reduction system for linear logic,
{\it Theory and Applications of Categories}, 35:1833-1870, 2020.

\bibitem{heho}
W.~Heijltjes and R.~Houston, {\rm No proof nets for MLL with units:
 Proof equivalence in MLL is PSPACE complete}, in {\it 
Proceedings of the Joint Meeting of the Twenty-Third EACSL
 Annual Conference on Computer Science Logic 
 and the Twenty-Ninth Annual ACM/IEEE Symposium on Logic in Computer Science,
 CSL-LICS '14}, ACM, 2014.

\bibitem{hugh}
D.~J.~D.~Hughes, {\rm Simple free star-autonomous categories
 and full coherence}, {\it Journal of Pure and Applied Algebra},
 216(11):2386--2410, 2012.

\bibitem{jaco}
B.~Jacobs, {\it Categorical Logic and Type theory}, Elsevier, 2001.

\bibitem{joya}
A.~Joyal, Foncteurs analytiques et esp\`eces de structures,
 {\it Combinatoire Enum\'e\-rative}, Proceedings, Montreal, Qu\'ebec,
 Canada, 1985, G.~Labelle, P.~Leroux, eds., Lecture Notes in
 Mathematics 1234, (Springer, 1986) pp. 126--159.

\bibitem{jost}
A.~Joyal and R.~Street, Braided tensor categories,
 {\it Advances in Mathematics}, 102:20-78, 1993.

\bibitem{kell}
G.~M.~Kelly, {\rm On MacLane's conditions for coherence of natural
 associativities, commutativities, etc.}, {\it Journal of
 Algebra}, 1(4):397--402, 1964.

\bibitem{lasc}
J.~Lambek and P.~J.~Scott, {\it Introduction to Higher-Order Categorical
 Logic}, Cambridge Universit Press, 1988.

\bibitem{lamp}
J. Lamping. An algorithm for optimal
 lambda calculus reductions. {\it Proceedings of the 17th Symposium on
 Principles of Programming Languages (POPL 90)}, San Francisco, C.A., 
U.S.A., 1990, pp.16--30.

\bibitem{mmpr}
M.~E.~Maietti, P.~Maneggia, V.~de~Paiva, and E.~Ritter,
 {\rm Relating categorical semantics for intuitionsitic linear logic},
 {\it Applied Categorical Structures}, 13(1):1--36, 2005.

\bibitem{sefl}
R.~Sedgewick and P.~Flajolet, {\it An Introduction to the
 Analysis of Algorithms}, Addison-Wesley Publishing Company,
 1996.

\bibitem{sips}
M.~Sipser, {\it Introduction to the Theory of Computation},
PWS Publishing Company, 1997.

\end{references*}


\end{document}